\documentclass[]{amsart}
\usepackage{amsmath,amsfonts,amssymb,amsthm,float,verbatim, xspace}
\usepackage{epsfig}
\begin{document}

\newtheorem{thm}{Theorem}[subsection]
\newtheorem{prop}[thm]{Proposition}
\newtheorem{lemma}[thm]{Lemma}
\newtheorem{cor}[thm]{Corollary}
\newtheorem{dfn}[thm]{Definition}
\newtheorem{axiom}[thm]{Axiom}
\restylefloat{figure}

\numberwithin{figure}{section}
\numberwithin{equation}{section}
\renewcommand{\thesubsubsection}{}
\renewcommand{\thefootnote}{\fnsymbol{footnote}}

\newcommand{\no}{\noindent}
\newcommand{\vv}{\vspace{2mm}}
\newcommand{\hh}{\hspace{2mm}}
\newcommand{\h}{\ensuremath{\epsilon}\xspace}
\newcommand{\fg}{\mathfrak{g}}
\newcommand{\fgs}{{\mathfrak{g}}^{\star}}
\newcommand{\ug}{\mathcal{U}(\fg)}
\newcommand{\uge}{\mathcal{U}(\fg_{\epsilon})}
\newcommand{\ech}{\ensuremath{\mathcal{H}}\xspace}
\newcommand{\bech}{\ensuremath{\bar{\mathcal{H}}}\xspace}
\newcommand{\st}[1]{{\star}_{#1}}
\newcommand{\stu}{\circledast}
\newcommand{\aro}{\epsfig{file=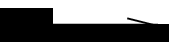, height=1.9mm}}
\newcommand{\arow}{\epsfig{file=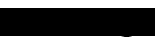, height=0.4mm, width=8mm}}
\newcommand{\bc}{\mathbb C}
\newcommand{\br}{\mathbb R}
\newcommand{\bz}{\mathbb Z}
\newcommand{\bn}{\mathbb N} 
\newcommand{\bp}{\mathbb P}
\newcommand{\Bn}[1]{{\hat{B}}_{#1}}
\newcommand{\Bc}[1]{{\bar{B}}_{#1}}
\newcommand{\Be}[1]{{\mathcal{B}}_{#1}}
\newcommand{\Bb}[1]{{\hat{\mathcal{B}}}_{#1}}
\newcommand{\la}[2]{\ensuremath{\lambda_#1^#2}}
\newcommand{\lak}[1]{\ensuremath{\lambda_k^#1}}
\newcommand{\Pol}[2]{{\mathcal{P}}_{#1}({#2})}
\newcommand{\vac}[1]{\underline{#1}}

\def\bea{\begin{eqnarray*}}
\def\eea{\end{eqnarray*}}
\def\eg{e.g.\ }
\def\ie{i.e.\ }
\def\stp{$\star$-product\xspace}

\newcommand{\pt}[2]{\frac{\partial #1}{\partial #2}}

\newcommand{\biv}[2]{\overleftarrow{\frac{\partial}{\partial #1}}
\overrightarrow{\frac{\partial}{\partial #2}}}

\newcommand{\bii}[2]{\overleftarrow{{\partial}_#1}\overrightarrow{{\partial}_#2}}
\newcommand{\pois}{\alpha^{ij} \overleftarrow{{\partial}_i}\overrightarrow{{\partial}_j}}

\newcommand{\del}[1]{\ensuremath{\partial_#1}}
\newcommand{\delt}[2]{\ensuremath{{\partial_#1}(#2)}}
\newcommand{\dell}[1]{\overleftarrow{{\partial}_#1}}
\newcommand{\delr}[1]{\overrightarrow{{\partial}_#1}}
\newcommand{\delle}{\overleftarrow{\partial}}
\newcommand{\delre}{\overrightarrow{\partial}}

\newcommand{\matr}[4]{\ensuremath{\left[ \begin{array}{cc}#1&#2\\#3&#4\end{array}\right] } }
\newcommand{\dfequal}{\stackrel{\mbox {\tiny {def}}}{=}}
\newcommand{\smalcirc}{\mbox{\tiny{$\circ $}}}

\newcommand{\ctn}{\bc^{2n}}
\newcommand{\rtn}{\br^{2n}}
\newcommand{\rn}{\br^{n}}
\newcommand{\rd}{\br^{d}}
\newcommand{\rtm}{\br^{2m}}

\newcommand{\cala}{{\mathcal A}}
\newcommand{\smo}{\sum_{j=0}^{\infty}{\epsilon}^{j}}
\newcommand{\smj}[1]{\sum_{#1}}

\title[Universal Formulae for Deformation Quantization]{Kontsevich's Universal Formula for Deformation Quantization and the Campbell-Baker-Hausdorff Formula, I}
\author{Vinay Kathotia}
\address{Department of Mathematics \\ University of California \\
Davis, CA 95616}
\email{kathotia@math.ucdavis.edu}
\date{November 30, 1998}
\begin{abstract}
\no We relate a universal formula for the deformation quantization
of  Poisson structures ({\stp}s) on $\rd$ proposed by Maxim Kontsevich to the Campbell-Baker-Hausdorff formula. Our basic thesis is that 
exponentiating a suitable deformation of the Poisson structure provides
a prototype for such formulae. For 
the dual of a Lie algebra, the \stp given by the universal enveloping algebra via symmetrization is shown to be of this type. In fact this \stp is essentially given by the Campbell-Baker-Hausdorff (CBH) formula. We call this the CBH-quantization. Next we limn Kontsevich's  construction using a graphical representation for differential calculus. We outline a structure theory for the weighted graphs which encode
bi-differential operators in his formula and compute certain  weights. We then establish that the Kontsevich and CBH quantizations are identical  for the duals of nilpotent Lie algebras. Consequently part of Kontsevich's \stp is determined by the CBH formula. Working the other way, we have a graphical encoding for the CBH formula. We end with some questions that could help clarify the overall picture. 
\end{abstract}
\maketitle
\tableofcontents

\setcounter{section}{-1}

\section{Introduction}
\label{chap0}
\no {\bf Overview.} Deformation quantization \cite{bf} aims to provide a mathematical model for quantum mechanics in accordance with Dirac's correspondence principle (the principle is outlined in the chapter on {\em quantum conditions} in \cite{di}). One starts with Poisson manifolds as models for classical mechanics and {\em quantizes} them by introducing a new multiplication on the space of functions (observables). The correspondence principle is addressed by requiring that the multiplication be {\em associative} and that it be a deformation of pointwise multiplication in the direction of the Poisson bracket. For an overview of the field please refer to \cite{we}.

\vv

\no In this paper (which was essentially the author's Ph.D. dissertation at Berkeley) we
investigate certain universal formulae for (formal and bi-differential) deformation quantizations of $\rd$. By a universal formula on $\rd$  we mean one that works for a whole class of Poisson structures. \eg , the Moyal product \cite{mo}, the motivating example for the theory,  is a universal deformation quantization formula for {\em constant} Poisson structures. One can also interpret the quantization for the dual of a Lie algebra arising from its universal 
enveloping algebra as  a universal formula for linear Poisson structures.

\vv

\no  In September 1997 
Maxim Kontsevich  \cite{ko:defq} presented a scheme for the  deformation 
quantization of all Poisson manifolds. As a first step he constructed a universal deformation quantization for any open domain in $\rd$. His method is motivated by topological quantum 
field theory and is based on a graphical representation for 
bi-differential operators. 
\vv

\no The aim of this paper is to relate the above formulae - Kontsevich's quantization and the one coming from the universal enveloping algebra. Our
 guiding principle, motivated by the Moyal product, is that one can 
obtain universal formulae by {\em exponentiating}  
 the Poisson bracket (or a suitable deformation of it).

\vv

\no We will show that both the Kontsevich and the universal enveloping algebra quantizations can be interpreted as being obtained 
via the exponentiation of a deformation of the Poisson bracket. In the Lie 
algebra case the exponential setting is evinced by the occurrence of the Campbell-Baker-Hausdorff 
(CBH) formula. Recall that for $X, Y$ in some Lie algebra, and for $\exp(X) \exp(Y) = \exp(C)$, the CBH formula provides a determination of $C$ as a formal series whose terms are elements in the Lie algebra generated by $X$ and $Y$.

\vv

\no In addition to establishing a connection at the conceptual level, we are able to 
show that the CBH formula is {\em contained} in Kontsevich's construction.
This allows us to  encode the CBH formula graphically and provides a new scheme for computing its coefficients.

\vv

\no{\bf Layout (or how to skip most of the material).} The layout of the sections is as follows:

\vv

\no In Section 1 we define deformation quantization and present some examples, namely the Moyal product and a deformation quantization for linear Poisson structures (duals of Lie algebras) based on the universal enveloping algebra.

\vv

\no The next section goes into the details of quantizing the linear Poisson structures . An approach based on the CBH formula is presented. It is shown that this leads to the same quantization as the one obtained using  the universal enveloping algebra. The occurrence of the Bernoulli numbers in both approaches is presented. We assemble some identities satisfied by the Bernoulli numbers for later use.

\vv

\no Section 3 introduces a graphical representation for differential calculus. This provides a graphical shorthand for addressing quantization questions and 
paves the way for describing Kontsevich's construction. We also show how the same graphical set-up yields an encoding for Lie algebras and for the CBH formula, yielding a much simplified framework for  the CBH quantization.

\vv

\no Section 4 delves into Kontsevich's construction. We discuss the nature of his {\em admissible} graphs, set up his weight integrals, and go about computingthem for a small subfamily of tractable graphs. Though these computations will not completely determine the Kontsevich quantization even for linear Poisson structures, they contain enough data to allow for a comparison with the CBH quantization. The main computation is presented in Section  \ref{4wc} where, not surprisingly, the Bernoulli numbers resurface.

\vv

\no In the last section we try to pull together the various strands. We explicitly relate
 Kontsevich's quantization to the one obtained from the CBH formula. We establish that both yield the same quantization for nilpotent Lie algebras. We also point out how Kontsevich's formula is of the form \hh `exp(suitable deformation of the Poisson structure)'. We end by briefly outlining some natural questions that seem to arise from the analysis presented here. 

\vv

\no {\bf Note.} Those primarily interested in Kontsevich's construction can skip most of Sections 1 and 2. The only results there that are used later in the weight computations are certain identities for the Bernoulli numbers that are summarized in Section \ref{secber}, and an isomorphism criterion for deformation quantizations (Section \ref{c2s1}) that is used in Section 5. Also, part of Section 3 can be viewed as low-tech motivation for Kontsevich's formula.

\vv

\no {\bf Plans.}  In a forthcoming  sequel to this paper we will study the efficacy of using  Kontsevich's graphs and weights to determine the CBH formula.

\vv

\no {\bf Acknowledgements.}  I am grateful to  Felix Alexandrovich Berezin and to Maxim Kontsevich for the ideas underlying this work. That these ideas have not remained inchoate is mainly due  to the encouragement and guidance  of Alan Weinstein.
{\bf Thank you.   }



\section{Deformation Quantization}
\label{chap1}
\subsection{The Definition}
\label{1.1}
\no Given our self-imposed restriction to $\rd$ we will not make an effort to state results in their general setting . Working in the formal differentiable category means that the basic objects of interest will be formal power series of smooth bi-differential operators. 

\vv

\no {\bf Notation:} Unless stated otherwise, $\cala$ will denote the algebra of $\br$ valued smooth functions on $\rd$. \\
Furthermore $\cala$ will be equipped with a Poisson structure $\alpha$ with bracket $\{ a , b \} $ for $ a,b \in \cala$. $\h$ will denote a formal parameter \footnote{We have chosen $\h$ as our deformation parameter instead of the more commonly used $\hbar$  in order to avoid introducing complex numbers. Setting $\h = i\hbar$ in Definition \ref{dfn-dq} would yield the standard formulation of deformation quantization.}.

\vv

\no We will follow the Einstein summation convention of summing over repeated indices over the appropriate range. At times we will introduce the summation sign $\sum$ if it helps to clarify matters.
 
\vv

\begin{dfn}
\label{dfn-dq}

A deformation quantization of $\rd$ in the direction of $\alpha$ is an $\br[[\h]]$ linear product on $\cala[[\h]]$, defined for $a, b \in \cala$ by
$$ a \star b = \smo \Pi_{j}(a,b) \; ,$$ 
where  the $\Pi_{j}$ are bi-differential operators satisfying
\begin{align}
(D1) \qquad&\smj{j+k=n} \Pi_{j}( \Pi_{k}( a , b) , c) = \smj{j+k=n} \Pi_{j}(a, \Pi_{k}( b , c ) )  & &\quad\text{(Associativity)} \notag \\
(D2) \qquad&\Pi_{0}(a,b) = ab & &\quad\text{(Classical Limit)} \notag \\
(D3) \qquad&\Pi_{1}(a,b) - \Pi_{1}(b,a) = \{ a , b \} & &\quad\text{(Semi-Classical Limit)} \notag
\end{align}

\no and $\star$ extends naturally to $\cala[[\h]]$.
\end{dfn}

\vv

\no Sometimes the $\Pi_{j}$ are not required to be bi-differential operators and the term $\star$-product is used for bi-differential deformation quantizations. We will use the two terms interchangeably. We now discuss some examples.

\subsection{The Moyal product}
\label{secmoy}
\no  The motivating example for the entire theory is the Moyal product for a constant Poisson structure $\alpha$ on $\rd$. 

\no Let
\begin{equation}
\label{podef}
\alpha = \frac{1}{2}{\alpha}^{ij} \partial_{i} \wedge \partial_{j} \hh , \hh {\alpha}^{ij} = - {\alpha}^{ji} \in \br
\end{equation}

\no where we are using the Einstein summation convention, $ \partial_{i} = \pt{}{x_{i}}$, and the $x_{i}$ are coordinates on $\rd,  i = 1, \ldots , d$.

\vv

\no  Note that viewed as a bi-differential operator the Poisson bracket associated to  $\alpha$ is given by

\begin{equation}
\label{podiff}
\hat{\alpha} (a,b) = \alpha^{ij} \bii{i}{j} (a,b)  \; ,
\end{equation}

\vv

\no where the $\dell{i}$ act on $a$ and the $\delr{j}$ act on $b$.
The Moyal $\star$-product is then given by simply exponentiating this Poisson operator. Symbolically,

\begin{equation}
\label{moya1}
a \star b = \exp(\frac{\h}{2} \alpha^{ij} \bii{i}{j}) (a,b)  \; .
\end{equation}

\vv

\no Note that  all the elements in the argument of the $\exp$ commute, and the zeroth power of any bi-differential operator stands for usual pointwise multiplication. 

\no Here are the first few terms of the expanded product:

\begin{equation}
\label{moya2}
 a \star b = ab + \frac{\h}{2} \alpha^{ij} \delt{i}{a} \delt{j}{b} + \frac{\h^{2}}{2^{2}2!} \alpha^{ij} \alpha^{kl} \delt{i}{\delt{k}{a}} \delt{j}{\delt{l}{b}} + \ldots \; .
\end{equation}

\vv

\no At this stage we will not prove that equation \eqref{moya1} yields a \stp. It will follow as a corollary of more general results that will be proved in the next section. 

\vv

\no The basic problem in trying to generalize the exponentiation idea to non-constant Poisson structures is that the $\dell{i} , \delr{j}$ no longer commute with the $\alpha^{ij}$ (though they still commute with each other). Invoking the imagery of exponentiating non-commuting variables is the perfect lead-in to the Campbell-Baker-Hausdorff (CBH) formula. The CBH formula does yield a \stp, but we defer its study to the next section. We now move on to our second example of deformation quantization.

\subsection{Linear Poisson Structures}
\label{secug}
\no Considering a linear Poisson structure on $\rd$ is equivalent to studying the dual of a Lie algebra, $\fgs$, with the Poisson structure being given by the
Lie bracket on $\fg$. Let $\{ X^{i} , \; i = 1, \ldots , d\}$ be a basis for $\fg$. Simultaneously we will think of the $X^{i}$ as giving coordinates on $\fgs$. Given the structure constants $c_k^{ij}$, the bracket
\begin{equation}
\label{cijk}
[ X^{i} , X^{j} ] = c_k^{ij} X^{k}  
\end{equation}

\vv

\no determines the Poisson structure
\begin{equation}
\label{linpo}
 \alpha = \frac{1}{2} c_k^{ij} X^{k} \del{i} \wedge \del{j} \; .
\end{equation}

\vv 

\no We construct a deformation quantization for $\fgs$ using the universal enveloping algebra of $\fg$. 

\subsubsection{The $\mathcal{U}$-quantization}
\label{ucal}
\no Recall that the universal enveloping algebra $\ug$ of a Lie algebra $\fg$ is the associative algebra obtained via
\begin{equation}
\label{unive}
 \ug  = \mathcal{T}(\fg) / \mathcal{I} \hh , \hh \mathcal{I} = <X{\otimes}Y - Y{\otimes}X - [ X , Y ]> \; ,
\end{equation}

\vv

\no where $\mathcal{T}(\fg)$ is the tensor algebra of $\fg$ and one quotients out the ideal $\mathcal{I}$ generated by terms of the type $(X \otimes Y - Y \otimes X - [ X , Y ]) \hh , \hh X,Y \in \fg.$ One can modify the construction slightly to obtain $\uge$,
\begin{equation}
\label{univd}
 \uge = \mathcal{T}(\fg)[[\h]] / {\mathcal{I}}_{\h} \hh , \hh {\mathcal{I}}_{\h} = <X{\otimes}Y - Y{\otimes}X - \h [ X , Y ]> \; .
\end{equation}

\vv

\no Notice that all we have done is re-scale the bracket by the formal parameter $\h$. 

\no $\uge$ provides a model for the quantized Poisson algebra as follows. 

\vv

\no Let  $S(\fg)$ be the symmetric tensor algebra of $\fg$, \ie the space of polynomial functions on $\fgs$. The Poincare-Birkhoff-Witt theorem leads to a vector space  isomorphism between $S(\fg)[[\h]]$ and $\uge$. The isomorphism $\sigma : S(\fg)[[\h]] \to \uge$ is given by the symmetrization of monomials,

\begin{equation}
\label{symdef}
\sigma(X^{1}X^{2} \ldots X^{n}) = \frac{1}{n!}\smj{\pi \in \mathfrak{S}_n}
X^{\pi(1)} \circ X^{\pi(2)} \circ \ldots \circ X^{\pi(n)} \; ,
\end{equation}

\vv

\no where $\mathfrak{S}_n$ is the permutation group of order $n$ and $\circ$ denotes the multiplication in $\uge$.
Now for $ p,q \in S(\fg)[[\h]]$  we define their $\stu$-product by

\begin{equation}
\label{uea:sp}
p \stu q = \sigma^{-1} ( \sigma(p) \circ \sigma(q) ) \; .
\end{equation}

\vv

\no Associativity of $\stu$ follows from the associativity of $\circ$. The properties $(D2)$ and $(D3)$ of Definition \ref{dfn-dq} are built into the definition of $\uge$. Furthermore $\stu$ is given by bi-differential operators. This follows from the fact that 
$  \sigma^{-1} (X^1 \circ X^2 \circ \ldots \circ X^n)$ is a sum of monomials whose factors are the $X^i$ and/or terms obtained from multiple bracketing of various $X^i$.
Being given by bi-differential operators, $\stu$   extends to a $\star$-product on all the smooth functions.

\vv

\no We wish to point out that in the above construction the introduction of $\h$ seems to be extraneous, with $\h = 1$ representing the {\em natural} theory. One way to think of this is that we are  working with the Poisson structure $\h \alpha$ instead of $\alpha$. More important, one should see the $\h$ factor in a \stp expansion, or rather its exponent, as simply coding the number of commutator (or Poisson) brackets contributing to each term. 
\vv

\no We now proceed to explore the $\fgs$ quantization in greater detail.

\section{Quantizing \protect{$\fgs$}}
\label{chap2} 
\subsection{Basic Set-up}
\label{c2s1}
\no  Throughout this section we will be working with $\fgs \cong\rd$ equipped with the Poisson structure \ $\frac{1}{2} c_k^{ij} X^{k} \del{i} \wedge \del{j}$ \ inherited from the Lie bracket on $\fg$. Polynomials and bi-differential operators will be allowed to have coefficients in $\br[[\h]]$. We plan to study and relate two quantizations for $\fgs$, the $\mathcal{U}$-quantization introduced in Section \ref{secug} and the quantization based on the CBH formula (CBH-quantization)
that was promised at the end of Section \ref{secmoy}.
Before we discuss the individual {\stp}s we will explore their general nature. We will establish the following facts. 

\vv

\newcounter{b}
\begin{list}
{(P\arabic{b})}{\usecounter{b}}
\item For any two polynomials $p_{n}$ and $q_{m}$ on $\fgs$ (or 
equivalently $p_{n}, q_{m} \in S(\fg)$) of degree
$n$ and $m$ respectively, both the $\star$-products mentioned above satisfy
\begin{equation}
\label{stform}
p_{n} \star q_{m} = p_{n} q_{m} + r_{m+n-1} \; ,
\end{equation}
\vv

\no where $r_{m+n-1}$ is a polynomial of degree $m+n-1$. \ie, modulo the pointwise product, $\star$ reduces total degree.
\item Given property (P1), the associativity of $\star$ 
implies that the $\star$-product of any two polynomials is determined 
by knowing $(X)^n \star Y$, where $n \in \bn$ and $X$ and $Y$ are arbitrary elements of $\fg$ and $(X)^n$ stands for the $n$th power of $X$ in $S(\fg)$.
\end{list}

\vv

\no {\bf Note:} To avoid any confusion between indexed elements $ X^n \in \fg$ and powers $(X)^n$ in $S^n(\fg)$ we shall always represent the latter with parentheses preceding the exponent as in  $(X)^n$.

\vv

\no Property (P1) follows for the $\mathcal{U}$-quantization because terms 
other than the pointwise product in $\stu$ are obtained by multiple 
bracketing, and the bracketing by the  Poisson structure, 
$\frac{1}{2} c_k^{ij} X^{k} \del{i} \wedge \del{j}$, reduces total degree by one. Later we shall see that the same reasoning gives us property (P1) for the CBH-quantization.

\vv

\no The proof of (P2) is essentially an exercise in induction and a
polarization identity and is outlined below. First we set up the 
following framework. 

\subsubsection{$\star$-monomials}
\no Given a monomial $p \in S(\fg),\ p = X^{i_{1}} X^{i_{2}} \ldots 
X^{i_{m}}$ we define $p^{\star}$ to be the $\star$-monomial 
$X^{i_{1}} \star X^{i_{2}} \star \ldots  \star X^{i_{m}}$. (Remember that the superscripts here are indices, not exponents!)

\vv

\no The basic result here is that the $\star$-monomials 
generate $S(\fg)[[\h]]$. 

\vv

\begin{lemma}
\label{Lemma1} Let $\mathcal{B}_{n}$ be a basis for 
$\sum_{j=0}^n S^{j} (\fg)$, the polynomials of degree $\leq n$, consisting of monomials ordered by decreasing degree,
$$ \mathcal{B}_{n} = \{ p_{1}, p_{2}, \ldots , p_{n_{d}} \} \; . $$

\vv
\no Here $n_d$ is the dimension of $\sum_{j=0}^n S^{j} (\fg)$. If $\star$ satisfies property (P1), the corresponding set 
$$\mathcal{B}^{*}_{n}  = \{ p_{1}^{\star}, p_{2}^{\star}, \ldots , 
p_{n_{d}}^{\star} \}  $$

\vv
\no is a spanning set for $S^{n}(\fg)[[\h]]$ over $\br[[\h]]$.
\end{lemma}

\begin{proof} Property (P1) implies that the transition matrix between the two sets 
$\mathcal{B}_{n}$ and $\mathcal{B}^{*}_{n}$ is given by
an $n_d \times  n_d$ lower triangular matrix with 1's on the diagonal.
\end{proof}

\subsubsection{Polarization} 
\no The aim here is to show that given $Y \in \fg$ 
and $p_{m}$  any homogeneous polynomial of degree $m$, $p_{m} \star Y$   is 
determined by $(X)^m \star Y$, where $X$ runs over $\fg$.
 This follows immediately from Lemma \ref{Lemma2} below.

\begin{lemma}
\label{Lemma2} Given a basis $\{ X^{1},  X^{2}, \ldots ,
X^{d} \}$ for $\fg$, $Y \in \fg$, 

\no $X^{i_{1}} X^{i_{2}} \ldots 
X^{i_{m}} \star Y$ can be obtained from the values of 
\begin{equation}
\label{polar}
 (t_{1} X^{1} + t_{2} X^{2} + \ldots  + t_{d} X^{d})^{m} \star Y
\end{equation}

\no as $(t_1, t_2, \ldots, t_d)$ range over $\rd$.
\end{lemma}

\begin{proof} We differentiate \eqref{polar} above with respect to $t_{i_{1}}, 
t_{i_{2}}, \ldots, t_{i_{m}}$ and then set all the $t_{k}$ equal to 
zero.
\end{proof}

\vv

\no We are now in a position to state our main theorem on the isomorphism of $\star$ products on $\fgs$.

\begin{thm}
A $\star$ product on $\fgs$ satisfying property (P1) is entirely determined by 
 $(X)^{n} \star Y$, where $X$ and $Y$ are arbitrary elements of $\fg$, $n 
\in \bn$.
\end{thm}

\begin{proof} It is adequate to show that the product is determined for arbitrary polynomials. To determine $p_{n} \star q_{m}$ we proceed by 
induction on the total degree of the expression, \ie $(n+m)$. 

\vv 

\no Assume that we have proved our result for products of polynomials with 
total degree $k$. Let us now consider $p_{n} \star q_{m}$ with $n+m = k+1$.  
Further let $q_{m} = Y^1 \ldots Y^{m-1} Y^m $. Then
Lemma \ref{Lemma1} tells us that, modulo polynomials of degree $k$,
\begin{align}
\label{eqipf}
p_{n} \star q_{m} & = \qquad p_{n} Y^1 \ldots Y^{m-1}Y^m & &\quad\text{(modulo degree $k$ terms)} \notag \\
& = \qquad p_{n} Y^1 \ldots Y^{m-1} \star Y^m \; . & &\quad\text{(modulo degree $k$ terms)}
\end{align}

\vv

\no But Lemma \ref{Lemma2} then implies that the right hand side of Equation \eqref{eqipf} is determined by the $(X)^k \star Y^m$. Note that it is safe to ignore the terms of degree $\leq k$ since they are determined via the induction hypothesis and Lemma \ref{Lemma1}.

\no So all we need to show is that we can start the induction, but that follows from the tautology $X \star Y = X \star Y$. 
\end{proof}

\vv

\no We now proceed to introduce the CBH-quantization. Subsequently we will show that the CBH-quantization and the $\mathcal{U}$-quantization are identical because their $(X)^{n} \star Y$ terms are equal. Moreover these terms have the Bernoulli numbers as coefficients, and we will try to explain the occurrence of the Bernoulli numbers in both settings. 

\subsection{The CBH-quantization}
\label{cbhsec}
\no The CBH-quantization is obtained in a fairly straightforward manner from the CBH formula. We would like to mention that though the derivation that follows is original, the CBH-quantization is implicit in Berezin's work \cite{be}. In fact most of this section can be viewed as a low-tech approach to the results outlined in his paper. We begin with an introduction to the CBH formula.

\subsubsection{The CBH Formula}
\no The Campbell-Baker-Hausdorff (CBH) formula states that for $X,Y \in \fg$ the formal series, $ H(X,Y) = \ln(\exp(X) \exp(Y))$ , is a Lie series.

\vv

\no Here $\ln$ and $\exp$ denote the usual series for those functions, and by a Lie series we mean a series whose terms are elements in the Lie algebra generated by $X$ and $Y$. H(X,Y) is sometimes called the {\em Hausdorff series}. Another way of stating the formula is that given the product of two group elements 
\begin{equation}
\label{cbh:df}
\exp(X) \centerdot \exp(Y) = \exp( H(X,Y) ) \; ,
\end{equation}

\no $ H(X,Y)$ is a formal series whose terms lie in the Lie algebra spanned by $X$ and $Y$. 

\vv

\no Here are its first few terms.  
\begin{equation}
\label{cbh:ex}
H(X,Y) = X + Y \; + \; \frac{1}{2} [X,Y] \; + \; \frac{1}{12} ([X,[X,Y]] + [[X,Y],Y]) \; + \ldots \; .
\end{equation}

\vv

\no Most proofs of the CBH formula, including the original (independent) proofs of Campbell, Baker, and Hausdorff use an iterative or recursive scheme to 
establish the formula. The invariable starting point is some identity in non-commutative calculus. None of these are of much use for computing the coefficients for the various commutators. We will not delve too much into the history and import of the formula. The reader is referred to Reutenauer's text \cite{re} and the references therein for information on pure mathematical aspects of the formula.

\vv

\no We would like to mention that, to the best of our knowledge, there is no reasonable explicit expression for the Hausdorff series. There is a formula due to Dynkin (see \cite{re}) which is of closed form, but to obtain individual terms one has to calculate sums which increase very rapidly in complexity. A much more manageable formula is presented in Newman and Thompson \cite{nt} where they build on earlier work of Goldberg \cite{go}. But here too, given a generic commutator expression, one has to sum an increasing number of coefficients to obtain the required coefficient. An additional complication arises from the fact that 
there are dependencies between various commutator expressions as a result of skew-symmetry and the Jacobi identity. This spanning set dependence of the formula makes it quite challenging to compare various approaches. We will probe these matters in the sequel to this paper.

\vv

\no We now return to our quantization problem.

\subsubsection{Group Multiplication}
\no Recall that we are simultaneously using $\{ X^{i} , \; i = 1, \ldots , d\}$ as a basis for $\fg$ as well as for coordinates on $\fgs$. This allows us to establish, at least for small $t \in \br$, a correspondence between group elements $\exp(tX)$, $X \in \fg$, and functions of the form $\exp(tX)$ on $\fgs$. Pulling back group multiplication to the functions via this
correspondence leads to the CBH-quantization. 

\vv

\no Starting with the CBH formula we will introduce a bi-differential operator $\hat{D}$ on $\fgs$. We will then show that $\hat{D}$ provides a \stp for functions of the form $\exp(tX)$, $X \in \fg$. But since these functions separate points on $\fgs$ the \stp extends to all of $\cala$.

\subsubsection{Symbols for {\stp}s}

\no We will try to define the bi-differential operator $\hat{D}$ as explicitly as possible. On occasion we will write $\hat{D}$ in the following expanded form
\begin{equation}
\label{Ddef}
\hat{D}(f,g) \dfequal \hat{D}(\vac{X},\hat{\vac{c}},(\dell{1}, \dell{2}, \ldots, \dell{d}) , 
(\delr{1}, \delr{2}, \ldots, \delr{d})) (f,g) \; .
\end{equation}

\vv

\no This is meant to reinforce the dependence of $\hat{D}$ on the coordinates functions $X^i$ (indicated here by the $\vac{X}$), on the structure constants $c_k^{ij}$ (indicated here by the $\hat{\vac{c}}$), and on the basic differential operators $\del{i}$.

\vv

\begin{dfn} The symbol of $\hat{D}$ as defined in Equation \eqref{Ddef} is the multivariable expression
\begin{equation}
\notag
D(\vac{X},\hat{\vac{c}},(s_1, s_{2}, \ldots, s_{d}) , 
(t_{1}, t_{2}, \ldots, t_{d}))
\end{equation}

\no obtained by replacing the $\dell{i}$ in $\hat{D}$ by $s_i$ and the $\delr{j}$ by $t_j$, where $s_i, t_j, \hh i,j \in \{1, \ldots, d\}$ are real-valued commuting variables.
\end{dfn}

\vv

\no For notational ease we will also use $\vac{s} = (s_1, s_2, \dots, s_d)$ and $\vac{t} = (t_1, t_2, \ldots, t_d)$. The symbol then being given by
\begin{equation}
\label{Dsdef}
{D}(\vac{X},\hat{\vac{c}},\vac{s} , 
\vac{t}) \; .
\end{equation}

\vv

\no It will turn out that ${D}(\vac{X},\hat{\vac{c}},\vac{s}, 
\vac{t})$ will depend linearly on $\vac{X}$ and is best considered as a formal power series in the $s_i$ and $t_j$ with (non constant) coefficients given by multiple products of the $c_k^{ij}$ with a single $X^l$. We now proceed to construct the symbol $D$ of $\hat{D}$ using the CBH formula.

\vv

\no Whereas $\exp(X)\centerdot\exp(Y)$ is given by exponentiating the Hausdorff series $H(X,Y)$ (Equation \eqref{cbh:df}), let us consider exponentiating only the terms that arise from bracketing in $H(X,Y)$. \ie, $H(X,Y) - X - Y$. Following Equation \eqref{cbh:ex} we have
\begin{equation}
\label{symbcb1}
\exp\bigl(H(X,Y) - X - Y\bigr) \; = \; \exp\bigl(\frac{1}{2} [X,Y] \; + \; \frac{1}{12} ([X,[X,Y]] + [[X,Y],Y]) \; + \ldots \bigr) \; .
\end{equation}

\vv

\no We now substitute $X = s_i X^i$ and $Y = t_j X^j$ in Equation \eqref{symbcb1}. Remember that we are using the Einstein summation convention. The resulting expression is precisely the symbol $D$ we seek.
\begin{equation}
\label{symbcb2}
\begin{split}
{D}(\vac{X},\hat{\vac{c}},\vac{s} , 
\vac{t}) \dfequal \exp\bigl(\; &\frac{1}{2} [s_i X^i,t_j X^j] \; + \\
& \frac{1}{12} ([s_i X^i,[s_l X^l,t_j X^j]] + [[s_i X^i,t_j X^j],t_l X^l]) \; + \ldots \bigr) \; .\\
\end{split} 
\end{equation}

\vv

\no Since the $s_i, t_j, \ldots $ are scalars we can extract them from the brackets. Further, $[X^i,X^j]$ translates to  $X^kc_k^{ij}$. This yields the following simplification(?),
\begin{equation}
\label{symbcb3}
D(\vac{X},\hat{\vac{c}},\vac{s} , 
\vac{t}) \; \dfequal \; \exp\bigl(\frac{1}{2} X^k c_k^{ij} s_i t_j \; + \; \frac{1}{12} (
X^k c_k^{im} c_m^{lj} s_i s_l t_j\; +\; X^k c_k^{ml} c_m^{ij} s_i t_j t_l) \; + \ldots \bigr) \; .
\end{equation}

\vv

\no We hope that the nature of the subsequent terms in $D(\vac{X},\hat{\vac{c}},\vac{s} , 
\vac{t})$ is clear from Equations \eqref{symbcb2} and  \eqref{symbcb3}. There should be no confusion about the exponentiation since for the purposes of the symbol, and for $\hat{D}$, the $X^i$ are being interpreted as coordinates on $\fgs$ and commute with one another. All the other terms, the $\hat{\vac{c}}, \vac{s}, \vac{t}$, commute freely. One now obtains the bi-differential operator $\hat{D}$ by 
replacing the $s_i$ by the $\dell{i}$ and the $t_j$ by the $\delr{j}$. This yields

\begin{equation}
\label{symbcb4}
\hat{D} = I +  \frac{1}{2} X^k c_k^{ij} \dell{i} \delr{j} \; +  \ldots 
\end{equation}

\vv

\no where we remind ourselves that the identity operator $I$ stands for pointwise multiplication. Equation \eqref{symbcb4} says that the \stp $\st{1}$ defined by 
\begin{equation}
\label{symbcb5}
f \st{1} g \dfequal \hat{D}(f,g)
\end{equation}

\vv

\no does give a deformation of pointwise multiplication in the direction of the Poisson bracket. The missing deformation parameter $\h$ can be introduced by
placing $\h c_k^{ij}$ instead of $c_k^{ij}$ in the definition of $\hat{D}$. We still need to show that this product is associative. We start by showing that it is associative on functions of the form $\exp(X), \; X \in \fg$. Consider $ \exp(s_i X^i) \st{1} \exp(t_j X^j) $. It follows from  Equation \eqref{Ddef} and from the rule for differentiating exponentials that
\begin{align}
& \exp(s_i X^i) \st{1} \exp(t_j X^j) \notag  \\
& \qquad =  \hat{D}\bigl(\vac{X},\hat{\vac{c}},(\dell{1}, \dell{2}, \ldots, \dell{d}) , 
(\delr{1}, \delr{2}, \ldots, \delr{d})\bigr) \bigl( \exp(s_i X^i), \exp(t_j X^j) \bigr) \notag \\
& \qquad =  D\bigl(\vac{X},\hat{\vac{c}},(s_{1}, s_{2}, \ldots, s_{d}) , 
(t_{1}, t_{2},\ldots, t_{d})\bigr) \exp(s_i X^i) \exp(t_j X^j) \; . \notag
\end{align}

\vv

\no The last line being a consequence of $\frac{d}{dX} \exp(tX) = t\exp(tX)$. Now using  Equation \eqref{symbcb2} we have
\begin{align}
& \exp(s_i X^i) \st{1} \exp(t_j X^j) \notag  \\
& \qquad = \exp\bigl(H(s_i X^i, t_j X^j) - s_i X^i - t_j X^j\bigr) \exp(s_i X^i) \exp(t_j X^j) \label{cbfinale} \\
& \qquad = \exp\bigl(H(s_i X^i, t_j X^j)\bigr) \; . \label{cbfinl}
\end{align}

\vv

\no But Equation \eqref{cbfinl} says that on functions of the form $\exp(X)$ the product $\st{1}$ corresponds to group multiplication via the CBH formula (compare Equation \eqref{cbfinl} with Equation \eqref{cbh:df}). This yields associativity for this family of functions. But since these functions separate all the points of $\fgs$ we obtain a bi-differential \stp on $\fgs$.

\subsubsection{Properties of the CBH-quantization}
\no Given the results of Section \ref{c2s1} we will concentrate on understanding the form of $(X)^n \st{1} Y, \; n \in \bn$ (since they determine the \stp). Thus we need not consider $\exp(s_i X^i) \st{1} \exp(t_j X^j) $ as we did in Equation \eqref{cbfinl},  instead we can simply work with $\exp(sX) \st{1} \exp(tY) \;,\; s,t \in \br$.

\vv

\no Now ignoring terms of $H(X,Y)$ which are of quadratic and higher degree in  $Y$, the CBH formula can be written as
\begin{equation}
\label{bern}
\exp(X)\centerdot\exp(Y) \cong \exp\bigl( X + Y + \frac{1}{2}[X,Y] + \sum_{k=2}^{\infty}\frac{B_k}{k!} \underbrace{[X,[X,\ldots,[X,[X,Y]]\ldots]]}_{k X's} \bigr) \; ,
\end{equation} 

\no where \hh $\cong$ \hh denotes equality modulo terms of quadratic and higher degree in $Y$. Alternately,
\begin{equation}
\label{bern1}
\exp(X)\centerdot\exp(Y) \cong  \exp\bigl(X + Y + \frac{1}{2}[X,Y] + \sum_{k=2}^{\infty}\frac{B_k}{k!} ({\rm ad}_X)^k(Y) \bigr) \; ,
\end{equation}

\vv

\no Here the $B_k$ are the Bernoulli numbers. We postpone the definition of the Bernoulli numbers and reasons for their arising to the end of this section. One thing worth pointing out at this stage is that 
\begin{equation}
\label{bink}
B_1 = -(\frac{1}{2}) \quad \mbox{and all the other $B_k$ for {\bf odd} $k$ are zero.}
\end{equation}

\vv

\no Given this information the $\frac{1}{2}[X,Y]$ term in Equation \eqref{bern1} can also be absorbed into the indexed sum if we work with $(-1)^k B_k$ instead of $B_k$. As we shall see, this is a more  natural variant of the Bernoulli numbers than appearances may suggest.

\vv

\no We return to  computing $\exp(sX) \st{1} \exp(tY)$. We will do so by expanding the corresponding version of  Equation  \eqref{cbfinale}. Notice that in the present case we can use the simpler form of $H(X,Y)$ as expressed in Equation \eqref{bern1}. Doing so we obtain,
\begin{align}
\label{berns}
&\exp(sX) \st{1} \exp(tY) \notag \\
&\quad \cong \exp\bigl( \frac{1}{2}[X,Y]st + \sum_{k=2}^{\infty}\frac{B_k}{k!}(({\rm ad}_X)^k(Y)) s^k t + \ldots \bigr) \exp(sX) \exp(tY) \; .
\end{align}

\vv

\no This yields the following simpler symbol for $\hat{D}$ restricted to products of the form $(X)^n \st{1} Y$.
\begin{equation}
\label{bersym}
 \exp\bigl(  \frac{1}{2}[X,Y]st + \sum_{k=2}^{\infty}\frac{B_k}{k!}(({\rm ad}_X)^k(Y)) s^k t \bigr) \; .
\end{equation}

\vv

\no Remember that $s$ represents $\dell{X}$ and $t$ represents $\delr{Y}$. Therefore in computing $(X)^n \st{1} Y$ we can safely ignore terms that are quadratic or higher in $t$. This allows us to write the relevant symbol in the still simpler form,
\begin{equation}
\label{bersym1}
I +  \frac{1}{2}[X,Y]st + \sum_{k=2}^{\infty}\frac{B_k}{k!}(({\rm ad}_X)^k(Y)) s^k t \; .
\end{equation}

\vv

\no where the $I$ as usual represents pointwise multiplication. Finally the above symbol gives us
\begin{align}
\label{sofin}
(X)^n \st{1} Y &= \bigl( I + \frac{1}{2}[X,Y]\dell{X}\delr{Y} + \sum_{k=2}^{\infty}\frac{B_k}{k!}(({\rm ad}_X)^k(Y)) (\dell{X})^k \delr{Y} \bigr) \hh ((X)^n,Y) \notag \\
&= (X)^n Y + \frac{1}{2}[X,Y]n(X)^{n-1} \notag \\
& \quad + \sum_{k=2}^{n}\frac{B_k}{k!}(({\rm ad}_X)^k(Y))
 n(n-1)\ldots(n-k+1)(X)^{n-k} \; .
\end{align}

\vv

\no The last sum is truncated at $n$ since subsequent terms include more than $n$ derivatives of $(X)^n$. Written in slightly condensed form,
\begin{equation}
\label{sofing}
(X)^n \st{1} Y = (X)^n Y + \frac{1}{2}n(X)^{n-1}[X,Y] + \sum_{k=2}^{n}\frac{n! B_k}{k! (n-k)!} (X)^{n-k}(({\rm ad}_X)^k(Y)) \; .
\end{equation}

\vv

\no As a check observe that the first and second terms on the right hand side of Equation \eqref{sofing} give the pointwise product and one half the Poisson bracket respectively. We also wish to point out that though we have avoided introducing \h in these computations the \h version of the \stp is easily obtained by scaling the Poisson bracket by \h.
\begin{equation}
\label{sofinh}
(X)^n \st{\epsilon} Y \dfequal (X)^n Y + \frac{\h}{2}n(X)^{n-1}[X,Y] + \sum_{k=2}^{n}\frac{\h^k n! B_k}{k! (n-k)!} (X)^{n-k}(({\rm ad}_X)^k(Y)) \; .
\end{equation}

\vv

\no We now return to the $\mathcal{U}$-quantization introduced in Section \ref{ucal} to show that it yields an identical result for $(X)^n \star Y$.

\subsection{The $\mathcal{U}$-quantization Revisited}
\label{secuq}
\no Once again we work with the $\h = 1$ theory since our arguments carry over to the more general case by simply rescaling the Lie bracket.  To compute $(X)^n \stu Y$ we resort to some straightforward combinatorial manipulation. Recall (Equation \eqref{uea:sp}) that $(X)^n \stu Y$ is defined by 
\begin{equation}
\label{chp2sp}
(X)^n \stu Y = \sigma^{-1} ( \sigma((X)^n) \circ \sigma(Y) ) \; ,
\end{equation}

\vv

\no where $\sigma$ is the symmetrization map defined in Equation \eqref{symdef}. Now Equation \eqref{chp2sp} readily simplifies to 
\begin{equation}
\label{chp2s2}
(X)^n \stu Y = \sigma^{-1} (X \circ X \circ \ldots \circ X \circ Y) \; . 
\end{equation}

\no Alternately,
\begin{equation}
\label{chp2s3}
\sigma ((X)^n \stu Y) = X \circ X \circ \ldots \circ X \circ Y \; . 
\end{equation}

\vv 

\no The only monomials in $S(\fg)$ that could make a non zero contribution to the left hand side of Equation \eqref{chp2s3} are those monomials which are built using  $n \; X$'s and a $Y$. But these are spanned by
\begin{equation}
\label{chp2s4}
(X)^{n-k}(({\rm ad}_X)^k(Y)) \quad \hspace{2mm} k \in \{0,1, \ldots, n\} \; .
\end{equation}

\vv

\no Having obtained all the terms that could possibly contribute to $(X)^n \stu Y$ we first observe that they match exactly with the terms that arose in the
CBH-quantization (Equation \eqref{sofing}). Therefore we collect terms using coefficients which reflect those in Equation \eqref{sofing}. \ie, let 
\begin{equation}
\label{chp2s5}
(X)^n \stu Y = \sum_{k=0}^{n}\frac{n! \Bn{k}}{k! (n-k)!}(X)^{n-k}({\rm ad}_X)^k(Y) \; ,
\end{equation}

\vv

\no where we still need to compute the $\Bn{k}$. 

\no For this we match the coefficients of $\underbrace{X \circ X \circ \ldots \circ X}_{n X's} \circ Y$ on both sides of Equation \eqref{chp2s3}. The right hand side coefficient is 1. For computing the left hand side let us first observe that 
\begin{equation}
\label{chp2s6}
\sigma(\underbrace{XX{\ldots}X}_{m X's}Z) = \frac{(X{\circ\ldots\circ}X{\circ}Z) \; +\;  (X{\circ\ldots\circ}X{\circ}Z{\circ}X) \; + \ldots +\; (Z{\circ}X{\circ\ldots\circ}X) }{m+1} \; .
\end{equation}

\vv

\no Now every term in the right hand side of Equation \eqref{chp2s5} is of the form $XX\ldots XZ$ where $Z =(({\rm ad}_X)^k(Y)), \; k = 0,1,\dots , n$. Therefore from Equations \eqref{chp2s5} and \eqref{chp2s6} it follows that the coefficient of $\underbrace{X \circ X \circ \ldots \circ X}_{n X's} \circ Y$ on the left hand side of Equation \eqref{chp2s3} is 
\begin{equation}
\label{chp2s7}
\sum_{k=0}^{n}\frac{n! \Bn{k}}{k! (n-k)!}\frac{1}{n-k+1} \; .
\end{equation}

\vv

\no This leads to the identity
\begin{equation}
\label{chp2s8}
\sum_{k=0}^{n}\frac{n! \Bn{k}}{k! (n-k)!}\frac{1}{n-k+1} = 1 \; .
\end{equation}

\vv

\no One can similarly compare the coefficients of  $Y \circ \underbrace{X \circ X \circ \ldots \circ X}_{n X's}$ on both sides of  Equation \eqref{chp2s3}. This time one has no contribution from the right hand side. But on the left hand side, as $k$ increases, the contribution of $({\rm ad}_X)^k(Y)$ alternates in sign. This yields the identity
\begin{equation}
\label{chp2s9}
\sum_{k=0}^{n}\frac{(-1)^k n! \Bn{k}}{k! (n-k)!}\frac{1}{n-k+1} = 0 \; .
\end{equation}

\vv

\no The identities \eqref{chp2s8} and \eqref{chp2s9} are quite similar to the most commonly given identity for the Bernoulli numbers. In fact as we shall see towards the end of this  section,
\begin{equation}
\label{chp2s10}
\Bn{k} = (-1)^kB_k \; .
\end{equation}

\vv

\no This means that $\Bn{0} = 1, \Bn{1} = \frac{1}{2}$, and the other $\Bn{k}$ are exactly the Bernoulli numbers (recall the statements immediately after Equation \eqref{bink}). But then $(X)^n \star Y$ is identical for the CBH and $\mathcal{U}$-quantizations (refer to Equations \eqref{sofing} and \eqref{chp2s5}). Hence the two quantizations are the same.

\subsection{Consolidation}
\label{consol}
\no We would like to consolidate some of the results obtained this far (pending the last subsection).

\newcounter{c}
\begin{list}
{(\arabic{c})}{\usecounter{c}}
\item The quantizations for $\fgs$ obtained via the universal enveloping algebra and via the CBH formula are equal and are determined by $(X)^n \star Y$ where$X$ and $Y$ are arbitrary elements in $\fg, n \in \bn$.
\item An oft repeated statement in the literature on the CBH formula is that all the coefficients in the CBH formula are determined by the Bernoulli numbers in that they are equal to sums of products of the Bernoulli numbers. The fact that the CBH-quantization is determined by $(X)^n \star Y$ is another restatement of this.
\item In our introduction we mentioned that one could consider obtaining deformation quantizations by exponentiating deformations of the Poisson bracket. We can now make this more precise in light of the CBH-quantization. In this case the deformation of the Poisson structure (or rather its symbol) is the argument of the $\exp$ in Equation \eqref{symbcb2} which is essentially given by the Hausdorff series. The other key point that needs clarification is  what we mean by exponentiation since one is dealing with bi-differential operators with non-constant coefficients. We used the transition to the symbols in order to do this. This boils down to saying that we multiply two bi-differential operators using the following simple rule
\begin{equation}
\label{con1}
(f \; \dell{A} \; \delr{B}) (g \; \dell{C} \; \delr{D}) = fg \; \dell{A}\dell{C} \; \delr{B}\delr{D} \; .
\end{equation}

\vv

\no Another way of saying it is that there is some deformation of the Poisson structure such that exponentiating the deformation by treating it as a constant coefficient bi-differential operator yields a \stp. 

\vv

\no Note that we have not characterized this magical deformation of the Poisson structure nor indicated how one can go about finding it. All we have shown is that in the $\fgs$ setting the Hausdorff series leads us to the required deformation. Later we will show that Kontsevich's product is of the same form. In fact the imagery suggested by his paper yields a simpler way of visualizing the above product and exponentiation for bi-differential operators.
\item The Moyal product introduced in Section \ref{secmoy} can be interpreted as a special case of the CBH-quantization. Given $\rd$ with a constant Poisson structure $\alpha$ (as defined in Equation \eqref{podef}) consider the 2-step nilpotent, $(d+1)$-dimensional Lie algebra with basis $\{ X^1, X^2, \ldots, X^d, H \}$, and bracket relations given by
\begin{equation}
\label{cbhmoy}
[X^i,X^j] = \alpha^{ij} H \; ,
\end{equation}

\vv

\no all other brackets being zero (\ie, $H$ commutes with everything). It follows that in this case the CBH-quantization is given by the bi-differential operator
\begin{equation}
\label{cbhmoy2}
\hat{D} \; = \; \exp(\frac{H}{2}{\alpha^{ij}}\dell{{X^{i}}}\delr{{X^{j}}}) \; .
\end{equation}

\vv

\no Since $H$ commutes with everything we can treat it as a scalar, in fact setting $H = \h \;$ gives the Moyal product as defined in Equation \eqref{moya1}.
\end{list}

\vv

\no We end this section with a discussion of the Bernoulli numbers.

\subsection{The Bernoulli Numbers}
\label{secber}
\no The Bernoulli numbers were discovered by Jacob Bernoulli (1654-1705) in the study of a universal formula for the following family of sums.
\begin{equation}
\label{2b1}
{S_n}(x) = 1^n + 2^n + \ldots + (x-1)^n \; .
\end{equation}

\vv

\no The information we shall gather here is quite standard (see \cite{gkp} for example). The simplest way of introducing the Bernoulli numbers is as the coefficients in the Taylor series 
\begin{equation}
\label{2b2}
\phi(x) \; \dfequal \; \frac{x}{ \exp (x) - 1} \; = \; \sum_{k=0}^{\infty} \frac{B_k}{k!}x^k \; .
\end{equation}

\vv

\no The most commonly given recurrence relation satisfied by the Bernoulli numbers is
\begin{equation}
\label{2b3}
\sum_{k=0}^{n} \binom{n+1}{k} B_k \; = \; 0 \quad n \geq 1 \quad ; \quad B_0 = 1 \; .
\end{equation}

\vv

\no It can be easily verified that $B_1 = - \frac{1}{2}$ and that for all the other  odd $k$ the $B_k$ are zero. 
If one defines $\Bc{k}$ by
\begin{equation}
\label{2b4}
\phi(-x) \; = \; \frac{-x}{ \exp (-x) - 1} \; \dfequal \; \sum_{k=0}^{\infty} \frac{\Bc{k}}{k!}x^k 
\end{equation}

\no then by the very definition of $\phi$
\begin{equation}
\label{2b4n}
B_k \; = \; (-1)^{k}\Bc{k} 
\end{equation}

\no and  given the vanishing of the  $B_k$ for $k$ odd, $k \neq 1$, it follows
that $\Bc{k} = B_k , k \neq 1$ \; and \; $\Bc{1} = \frac{1}{2}.$

\vv

\no But then adding $\binom{n+1}{1} (\Bc{1} - B_1)$ to Equation \eqref{2b3} yields
\begin{equation} 
\sum_{k=0}^{n} \binom{n+1}{k} \Bc{k} = \binom{n+1}{1} (\Bc{1} - B_1) \notag
\end{equation}
\begin{equation}
\mbox{or,} \quad \sum_{k=0}^{n} \binom{n+1}{k} \Bc{k} = n+1 \notag
\end{equation}
\begin{equation}
\mbox{or,} \quad \frac{1}{n+1}\sum_{k=o}^{n} \binom{n+1}{k} \Bc{k} = 1 \notag
\end{equation}
\begin{equation}
\mbox{or,} \quad \sum_{k=0}^{n}\frac{n! \Bc{k}}{k! (n+1-k)!} = 1 \; .
\end{equation}

\vv

\no But this is precisely the identity satisfied by the $\Bn{k}$ in Equation \eqref{chp2s8}. Similarly replacing $B_k$ by $(-1)^k \Bc{k}$ in Equation \eqref{2b3} and comparing it with Equation \eqref{chp2s9} also yields that the coefficients $\Bn{k}$ introduced in section \ref{secuq} are the same as the $\Bc{k}$. This proves the contention of Equation \eqref{chp2s10} that $\Bn{k} = (-1)^k B_k$.

\vv

\no Just as we have related the occurrence of the Bernoulli numbers in the $\mathcal{U}$-quantization to the identity \eqref{2b3}, one can relate the occurrence of the Bernoulli numbers in the CBH-quantization to their being the Taylor coefficients of $\phi(x)$. It turns out that computing the terms in the CBH formula (for $\exp(X) \centerdot \exp(Y))$ of the type $({\rm ad}_X)^n(Y)$ corresponds to expanding 
\begin{equation}
\label{2b6}
\frac{{\rm ad}_X}{1 - \exp(-{\rm ad}_X)} \; .
\end{equation}

\vv

\no \cite{re} contains a clear explanation for this or it is illuminating to  work out the simplest non trivial example, that of the 2-dimensional solvable Lie algebra. Now comparing Equations \eqref{2b4} and \eqref{2b6} explains the occurrence of the $\Bn{k}$ in the CBH formula.

\vv

\no We close this section with the definition of the Bernoulli polynomials, and a variant based on the $\Bn{k}$, as they will be needed later. 

\begin{dfn}
The $m$th Bernoulli polynomial $\Be{m}(x)$ is defined by

$$\label{bpold} \Be{m}(x) = \sum_{k=0}^m \binom{m}{k} B_{m-k} x^{k} \; .$$
\end{dfn}
\begin{dfn}

The $m$th $\Bn{}$ernoulli polynomial $\Bb{m}(x)$ is defined by

$$\label{bpold2} \Bb{m}(x) = \sum_{k=0}^m \binom{m}{k} \Bn{m-k} x^{k} \; .$$
\end{dfn}

\section{A Graphical Representation for Differential Calculus}
\label{chp3}
\subsection{ A Brief Overview}
\no The following graphical representation for differential calculus will be quite useful in understanding Kontsevich's work (\cite{ko:defq}). We start by listing the essential ingredients and then fill in details and examples.
\newcounter{g}
\begin{list}
{G-\arabic{g}.}{\usecounter{g}}
\item Functions are to be represented by points, or by {\em colored} points.
\item The sum of two functions will be represented by their disjoint union.
\item The product of two functions will be represented by placing them {\em indistinguishably} close (separable under {\em magnification}).
\item A vector field will be represented by a {\em short} arrow.  
\end{list}

\no Since the notions of scale that arise in the above are ill-defined at best,
we will follow the following conventions. Sums, i.e. disjoint unions, will be represented by placing the appropriate points in parentheses. On occasion we may introduce `+' signs if it helps clarify matters. Some of this would make a lot more sense if one thinks of addition as `{\em or}', and of multiplication as `{\em and}'.

\vv

\begin{figure}[H]
\centering{\epsfig{file=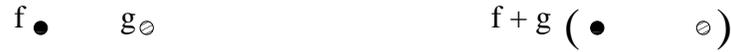 , height=6.1mm , clip= }}
\caption{Addition of functions}
\label{3p1}
\end{figure}

\no Products will be represented by placing points closer, and within a circle,
to denote that they should be seen as indistinguishably close.

\vv

\begin{figure}[H]
\centering{\epsfig{file=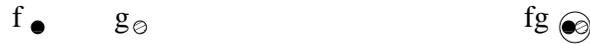 , height=5.3mm , clip= }}
\caption{Multiplication of functions}
\label{3p2}
\end{figure}

\no A vector field of the form $f\pt{}{x}$ will be denoted by $\aro$. Note that the $\bullet$ denotes the function $f$ and the $\arow$ the operator 
$\pt{}{x}$. More precisely, the vector field is obtained by a suitable `coloring' of $\aro$, where by a coloring we mean the following pair of maps.

\begin{dfn}
\label{colorv}

A {\bf coloring} of a point, $\bullet$, is the assignment of a function to the $\bullet$
$$ C_V : \{ \bullet \} \to  \cala \; .  \hspace{8mm} \mbox{(recall} \; \cala = C^{\infty} (\rd) \;) $$
\end{dfn}

\begin{dfn}
\label{colore}

Given a set of coordinates $\{x_1, x_2,  \ldots, x_d\}$ on $\rd$ a {\bf coloring} of an arrow, $\searrow$, is a map which assigns some $\del{i}$ to the $\searrow$ and is therefore given by a map $C_E$
$$ C_E :  \{ \searrow \} \to  \{ 1,2, \ldots, d \} \; .$$

\end{dfn}

\vv

\no The notion of coloring can be extended to all graphs made up of $\bullet$'s
and $\searrow$'s. In particular bi-differential operators of the form $\pois$
 will be called  {\bf wedges} and 
are coded by
\begin{figure}[H]
\centering{\epsfig{file=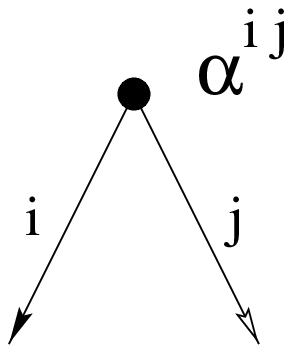 , height=16mm , clip= }}
\caption{A wedge}
\label{3p3}
\end{figure}

\no where we sum over colorings following the Einstein summation convention. The power of this method is that using a single wedge one can represent any Poisson tensor, without having to worry about the dimension $d$ of $\rd$ and the particular nature of the $\alpha^{ij}$. This makes it an ideal tool for studying universal formulae. 

\vv

\no Further, if one were
working over the polynomials, one could provide more detail by using
\begin{figure}[H]
\centering{\epsfig{file=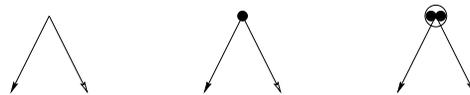 , height=12mm , clip= }}
\caption{More wedges}
\label{3p4}
\end{figure}

\no to encode constant, linear, and quadratic bi-vector fields respectively. We will  not do so and use $\bullet$'s to represent a generic function. We now display the sum and product rules for differentiation.

\begin{figure}[H]
\centering{\epsfig{file=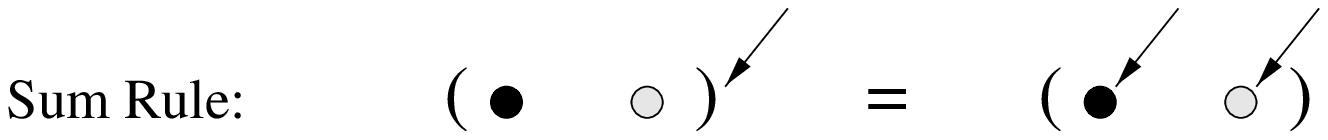 , height=8mm , clip= }}
\caption{Sum rule}
\label{3p5}
\end{figure}

\begin{figure}[H]
\centering{\epsfig{file=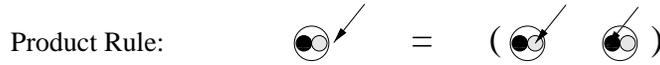 , height=8mm , clip= }}
\caption{Product rule}
\label{3p6}
\end{figure}

\no The product rule has the following visual interpretation. On looking closely (`zooming in') to see where an arrow has landed (on two proximal points) one finds it resting on one point {\em or} on the other. This is the basic principle at play when one uses graphs to encode the composition of differential operators, something we will be doing often when checking associativity for a $\star$-product.    
    
\vv
        
\no {\bf Note:} Technically speaking we should place the arrows and wedges inside circles but we will not always do that in order to reduce clutter. We will use some circling of operators when we compose two differential operators as this involves different levels of magnification. Hopefully this point will be clarified in subsequent examples.

\subsection{Iterative Symmetrizing}
\label{secflat}
\no Before we outline Kontsevich's construction of a universal \stp \ let us try to construct one using our rudimentary graphical calculus. We start with a product $\st{\alpha}$ given by 
\begin{equation}
\label{stal} 
f \st{\alpha} g = fg + \frac{1}{2}\{f,g\} 
\end{equation}

\vv

\no and try to modify the product by adding terms as mandated by associativity. As this is our first assay into the graphical world  we will (for a while)  write the equivalent algebraic expression below a graphical equation. Furthermore $\h$ will not be mentioned but one can think of it as implicitly scaling the Poisson structure. We also need to introduce a graphical representative for the \stp. We do this by using a circumscribing box to denote the \stp. \ie, $f \st{\alpha} g =$ \fbox{$f g$}. As a warm-up we restate the definition of $\st{\alpha}$ (Equation \eqref{stal}).

\begin{figure}[H]
\centering{\epsfig{file=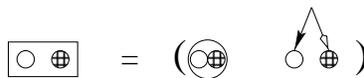 , height=1cm , clip= }}
\caption{$f \st{\alpha} g \hspace{5mm} = \hspace{2mm} fg \; +\; \frac{1}{2}\{f,g\}$}
\label{3p7}
\end{figure}

\no where {\epsfig{file=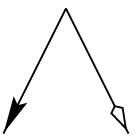 , height=6mm , clip= }} represents  $\frac{1}{2}\{\ ,\ \}$ and this is attained by the following coloring

\begin{figure}[H]
\centering{\epsfig{file=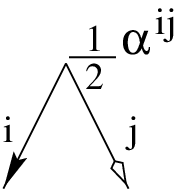 , height=14mm , clip= }}
\caption{$\frac{1}{2}\{\ ,\ \}$}
\label{3p7b}
\end{figure}

\no Next we consider $(f \st{\alpha}g) \st{\alpha} h$. 

\begin{figure}[H]
\centering{\epsfig{file=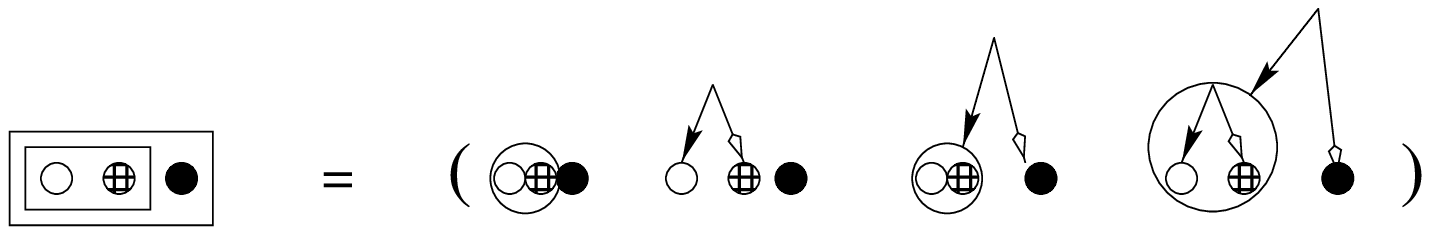 , height=16mm , clip= }}
\caption[$(f \st{\alpha}g) \st{\alpha} h$ - level I]{$(f \st{\alpha}g) \st{\alpha} h \quad= \; fgh \;+\; \frac{1}{2}\{f,g\}h \;+\; \frac{1}{2}\{fg,h\} \;+\; \frac{1}{4}\{ \{f,g\},h\}$}
\label{3p8}
\end{figure}

\no There should be a couple more bounding circles (to indicate multiplication) in Figure \ref{3p8}. We have omitted
them in the interest of clarity, simply retaining circles for terms that need to be expanded. Before we move on to the next step we will drop some more detail from our pictures. Given the ordering of the functions $(f,g,h)$ there should be no confusion even if we drop the distinguishing features of the various points and arrowheads. On the other hand it would make the pictures a lot easier to analyze. Here is what Figure \ref{3p8} looks like in this simplified setting.

\begin{figure}[H]
\centering{\epsfig{file=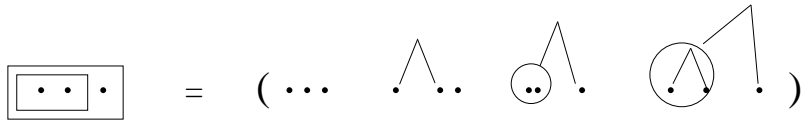 , height=16mm , clip= }}
\caption{Simpler version of Figure \ref{3p8}}
\label{3p9}
\end{figure}

\no We now expand the circled expressions ({\bf zoom in!}) using the product rule (see Figure \ref{3p6}). In words, the graph for $\frac{1}{2}\{fg,h\}$ will lead to two new graphs, one when the left `leg' of the wedge lands on $f$, and the other when it lands on $g$. Similarly one gets three new graphs from the $ \frac{1}{4}\{ \{f,g\},h\}$ term.

\begin{figure}[H]
\centering{\epsfig{file=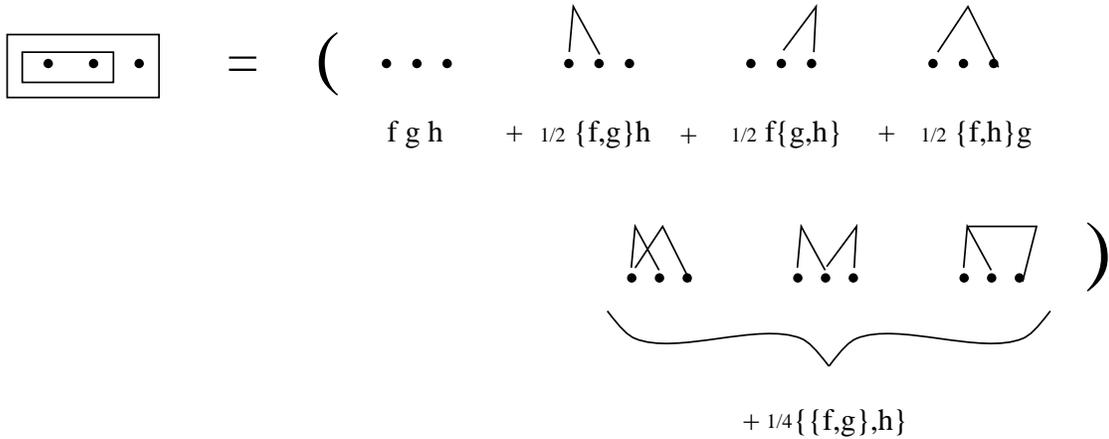 , height=58mm , clip= }}
\caption[$(f \st{\alpha}g) \st{\alpha} h$ - level II]{$(f \st{\alpha}g) \st{\alpha} h = fgh + \frac{1}{2}\{f,g\}h + f\frac{1}{2}\{g,h\} + \frac{1}{2}\{f,h\}g + \frac{1}{4}\{ \{f,g\},h\}$}
\label{3p10}
\end{figure}

\no Notice how there is no succinct way to expand the algebraic $ \frac{1}{4}\{ \{f,g\},h\}$, whereas things work just fine with the graphs. To simplify life, we now restrict ourselves to a constant Poisson structure. Let us call the corresponding product $\st{\beta}$. This allows us to drop the last term \epsfig{file=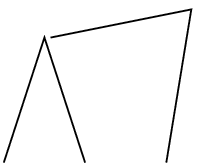, height=6mm, clip= } in Figure \ref{3p10}. We also stop writing the algebraic expansions.

\begin{figure}[H]
\centering{\epsfig{file=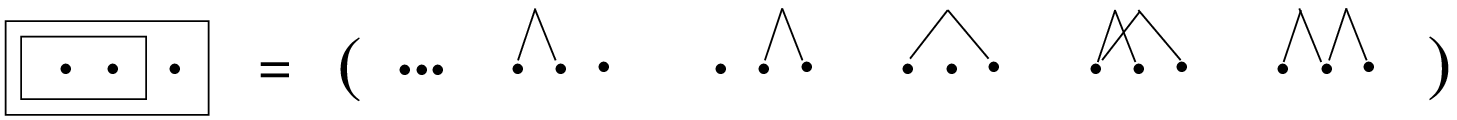 , height=1.2cm , clip= }}
\caption{$(f \st{\beta} g) \st{\beta} h$}
\label{3p11}
\end{figure}

\no Notice that if the graphs for $(f \st{\beta} g) \st{\beta} h$ in the right hand side of Figure \ref{3p11} were all present along with their mirror images (mirror images about a vertical line through the middle point representing the function $g$) then we would have associativity. To reinforce this point let us display  $f \st{\beta} (g \st{\beta} h)$.

\begin{figure}[H]
\centering{\epsfig{file=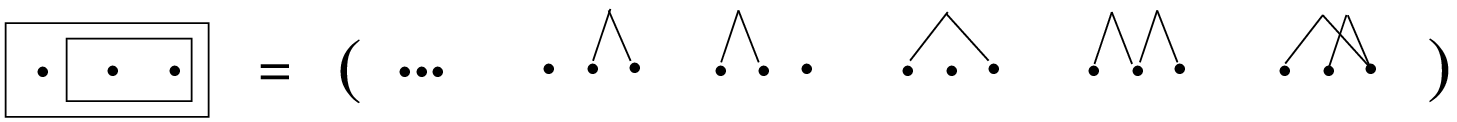 , height=1.2cm , clip= }}
\caption{$f \st{\beta} (g \st{\beta} h)$}
\label{3p12}
\end{figure}

\no Being able to introduce the term  \epsfig{file=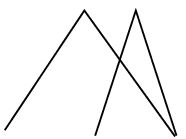, height=6mm, clip= }  in Figure \ref{3p11}  and the corresponding  \epsfig{file=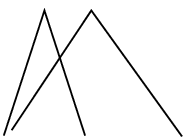, height=6mm, clip= }  in Figure \ref{3p12} would get us associativity up to the $\h^2$ level. Note that these graphs represent tri-differential operators whereas we can only modify our $\star$ product by adding bi-differential operators. There is a simple way of determining the bi-differential operator that would do the job. All we do is mege two of the vertices of the tri-differential operator/graph to obtain the bi-differential operator that needs to be added., \ie, being desirous of  \epsfig{file=tweed1.eps, height=6mm, clip= } means that we need to form (and add to our product) the following bi-differential operator.

\vv
\begin{figure}[H]
\centering{\epsfig{file=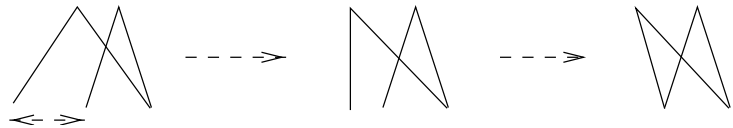 , height=16mm , clip= }}
\caption{Merging vertices - I}
\label{3p13}
\end{figure}

\no Note that we do not merge the two feet of a single wedge since the the skew-symmetry of the Poisson structure would cause the corresponding operator to vanish. Similarly from  \epsfig{file=tweed2.eps, height=6mm, clip= } we obtain

\vv

\begin{figure}[H]
\centering{\epsfig{file=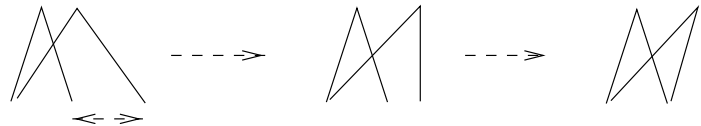 , height=16mm , clip= }}
\caption{Merging vertices - II}
\label{3p14}
\end{figure}

\no This gives us the right correction up to a scaling factor. To show this, let us analyze the effect of adding this bi-differential operator to our \stp. At the $\h^2$ level, new terms will arise due to the composition of pointwise multiplication and the \epsfig{file=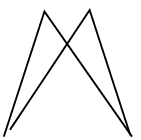, height=6mm, clip= } term that has been introduced. We consider the possible new terms for $(f \st{\beta} g) \st{\beta} h$.

\vv

\begin{figure}[H]
\centering{\epsfig{file=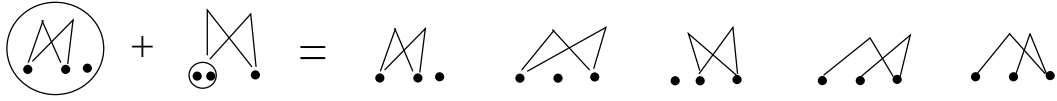 , height=13mm , clip= }}
\caption{Correction analysis}
\label{3p14aa}
\end{figure}

\no The left hand side of Figure \ref{3p14aa} shows the two ways in which \epsfig{file=tweed3.eps, height=6mm, clip= } can be composed with pointwise multiplication in $(f \st{\beta} g) \st{\beta} h$. Notice that as a result of this we get twice the term mandated by associativity (the last two terms on the right hand side of Figure \ref{3p14aa}) and some additional terms (the first three terms on the right hand side). The additional terms that are introduced provide no obstruction to associativity as they arise in  mirror symmetric (about the vertical axis) pairs.

\vv

\no Thus we add $\frac{1}{2}$(\epsfig{file=tweed3.eps, height=6mm, clip= }) to the product $\st{\beta}$. The corresponding colored graph is

\begin{figure}[H]
\centering{\epsfig{file=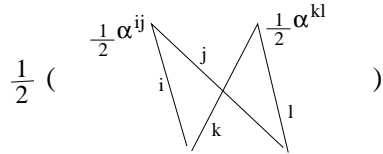 , height=2cm , clip= }}
\caption{Quadratic ($\h^2$) correction term}
\label{3p15}
\end{figure}

\no and in coordinates it is given by 
\begin{equation}
\frac{1}{2} \frac{\alpha^{ij}}{2} \frac{\alpha^{kl}}{2} \dell{i}\dell{k} \; \delr{j}\delr{l} \; .
\label{3p15eq}
\end{equation}

\vv

\no In the constant coefficient case it is not hard to see that at the nth stage one will need to introduce 

\begin{figure}[H]
\centering{\epsfig{file=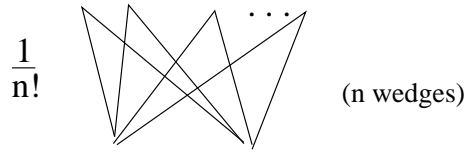 , height=2cm , clip= }}
\caption{General correction term}
\label{3p16}
\end{figure}

\no This gives us further insight regarding the occurrence of the exponential function in the Moyal product. Had we been working with a non-constant Poisson structure we would also need to `symmetrize' the \epsfig{file=tweedy.eps, height=6mm, clip= } term in $(f \st{\alpha}g) \st{\alpha} h$ (Figure \ref{3p10}).

As before, this requires that we introduce the bi-differential operator that is obtained by merging the vertices of \epsfig{file=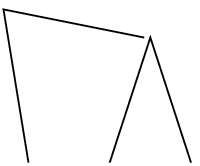, height=6mm, clip= } (the mirror image of \epsfig{file=tweedy.eps, height=6mm, clip= }).

\vv

\begin{figure}[H]
\centering{\epsfig{file=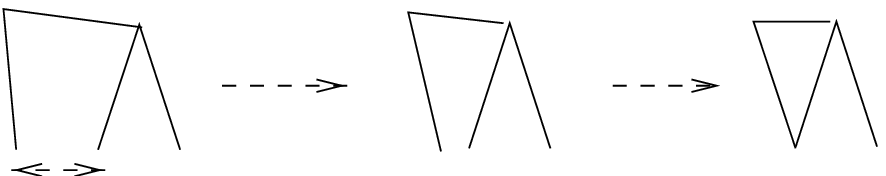 , height=2cm , clip= }}
\caption{Merging vertices - III}
\label{3p17}
\end{figure}

\no Similarly the corresponding term  \epsfig{file=tweedym.eps, height=6mm, clip= } in $f \st{\alpha} (g \st{\alpha} h)$ requires the introduction of \epsfig{file=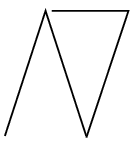, height=6mm, clip= }.

\no The net result is that we have to add the following correction term.

\vv

\begin{figure}[H]
\centering{\epsfig{file=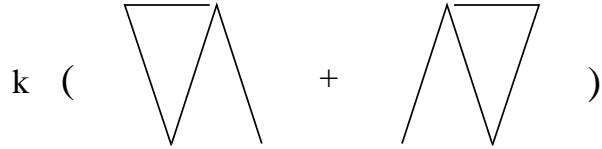 , height=2cm , clip= }}
\caption{Quadratic correction - non-constant case}
\label{3p18}
\end{figure}

\no The scaling factor $k$ still needs to be determined. One has to take the Jacobi identity into account, and it leads to a factor of $\frac{1}{3}$. There is also a determination of sign to be made here which is accounted for in the coloring. We will not go through the details of the computation (a couple of lines of diagrams). Suitably colored, the above correction term is given by

\vv

\begin{figure}[H]
\centering{\epsfig{file=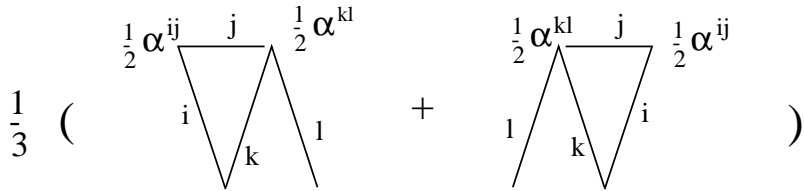 , height=25mm , clip= }}
\caption{Colored quadratic correction - non-constant case}
\label{3p19}
\end{figure}

\no which, in coordinates, translates to the following bi-differential operator.
\begin{equation}
\label{flat8}
 \frac{1}{12} (\alpha^{ij} \del{j}(\alpha^{kl}) \dell{i}\dell{k} \; \delr{j}) \; + \; \frac{1}{12} (\alpha^{ij} \del{j}(\alpha^{kl}) \dell{l} \; \delr{i}\delr{k}) \; .
\end{equation}

\vv

\no {\bf Some Observations.}

\no {\bf 1.} The correction term of Figure \ref{3p19} matches exactly with the term in the CBH-quantization that comes from  $ \frac{1}{12} ([X,[X,Y]] + [[X,Y],Y])$. The remainder of this paper is in some sense devoted to showing that this correspondence (the occurrence of $\exp$ and the matching of terms) between the CBH-quantization and `iterative symmetrizing' is not just coincidental. The CBH formula offers {\em the} prototype for universal deformation quantization formulae.

\vv

\no {\bf 2.} Subsequent direct computations for correction terms are cumbersome. On the other hand prior knowledge of some the coefficients can guide one towards more efficient computations. The advantage of Kontsevich's treatment (in \cite{ko:defq}) is that he simply bypasses any such iterative approaches though still relying on a graphical approach. We will briefly present some motivation for his method at the end of this section and then move on to defining and computing his product in Section 4.

\vv

\no{\bf 3.} Since the `wedges' allow us to code Poisson structures, in particular we can use them to represent Lie algebras. We describe that in some detail now.

\subsection{Graphical Encoding of Lie Algebras}
\label{freeg}
\no The graphical methods introduced in this section provide a natural framework for encoding Lie algebras. This can be done  at two levels. One can simply use wedges to denote bracketing, by using {\epsfig{file=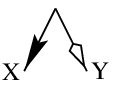 , height=8mm , clip= }} to represent $[X,Y]$, or one can think of the  wedge as representing a bi-differential operator colored by the linear Poisson structure. This duality gives us a simple way of explaining the CBH-quantization construction. The transition we effected in Section \ref{cbhsec} (using the notion of symbols) from the Hausdorff series of Lie elements to a bi-differential operator boils down to viewing the same collection of graphs in the two different ways indicated above. This subsection is designed to clarify this issue and to familiarize ourselves further with using graphs to represent bi-differential operators. Later this will be useful in relating the CBH and Kontsevich quantizations.

\subsubsection{The Free Lie Algebra on Two Generators}
\no To begin with, let us consider the free Lie algebra $\fg_F$ generated by {\bf two} elements, $X$ and $Y$. It is fairly straightforward to show that well-formed expressions in  $\fg_F$ can be identified with binary, planar, rooted {\em trees} with leaves labelled by $X$ and $Y$ (see \cite{re} for example). 

\no A brief note on the terminology. In our figures, roots correspond to {\em free} vertices (\ie, vertices with no wedges landing on them). Furthermore, our graphs are shown as descending from the roots, and `leaves' in the standard literature correspond to the feet at the bottom of our graphs.
 
\no The identification between Lie elements and rooted trees is based on the following principle. For $A$ and $B$ well-formed expressions in $\fg_F$, the tree for $[A,B]$ is obtained by taking a new root, with immediate left subtree given by (the tree for) $A$ and immediate right subtree given by $B$.

\vv

\begin{figure}[H]
\centering{\epsfig{file=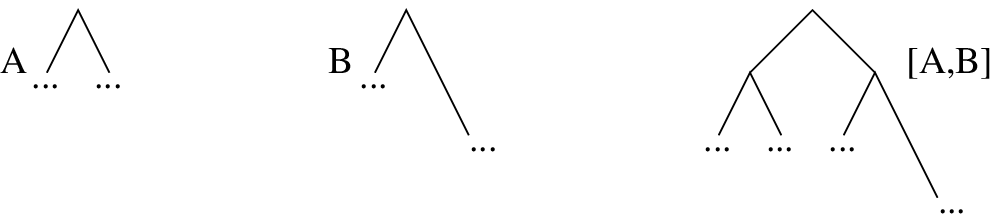 , height=2cm , clip= }}
\caption{Tree for $[A,B]$}
\label{cbg0}
\end{figure}

\no This is how one would represent [[X,[X,Y]],X].

\vv

\begin{figure}[H]
\centering{\epsfig{file=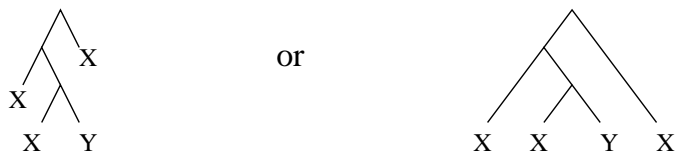 , height=2cm , clip= }}
\caption{Tree for $[[X,[X,Y]],X]$}
\label{cbg01}
\end{figure}

\no We carry this picture over to our wedge setting with one significant alteration. We place all the X leaves to the left and all the Y leaves to the right. We are able to do this as we can order the feet of our wedges by decorating them, and are not constrained to keep the order of bracketing consistent with the left-right ordering of the leaves. To make this clear we will encode $[X,Y]$ using the following {\em ordered}  wedge.

\vv

\begin{figure}[H]
\centering{\epsfig{file=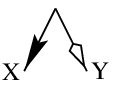, height=1cm , clip= }}
\caption{$[X,Y]$}
\label{cbg1}
\end{figure}

\no Here by ordered wedge we mean that there is an ordering for the feet of the wedge, set by declaring  \epsfig{file=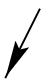 , height=8mm , clip= } to be the first leg (note that we do not necessarily require it to be the left leg). This determines the ordering of the bracketed terms. Given this setting, $[[X,[X,Y]],X]$ would be represented by 

\vv

\begin{figure}[H]
\centering{\epsfig{file=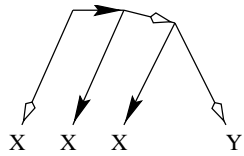 , height=2cm , clip= }}
\caption{Wedge based tree for $[[X,[X,Y]],X]$}
\label{cbg02}
\end{figure}

\begin{figure}[H]
\centering{\epsfig{file=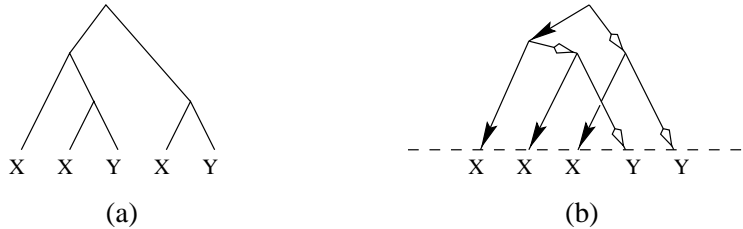 , height=3cm , clip= }}
\caption{Trees for $[[X,[X,Y]] , [X,Y]]$}
\label{cbg03}
\end{figure}

\no In preparation for Kontsevich's construction, we would like to specialize our representation even further. As a first step, we require that our graphs lie entirely in the closed upper half plane $\bar{\ech}$, with the $X$'s and $Y$'s being aligned along the real axis. Figure \ref{cbg03}(a) shows $[[X,[X,Y]] , [X,Y]]$ using the original tree representation displayed in Figure \ref{cbg0} whereas Figure \ref{cbg03}(b) uses the convention outlined in this paragraph.

\no We wish to point out that the crossing  \epsfig{file=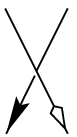 , height=8mm , clip= } in Figure \ref{cbg03}(b) is not a new vertex in the graph but simply a crossing introduced by the loss of planarity in the new representation.\footnote{Truly speaking this particular graph can still be represented by a planar version, but  $[[[X,[X,Y]],[X,Y]],[X,Y]]$ is an example for which one has to use a non-planar graph.}

\vv

\no We now make the final move to fix our graphical presentation for Lie algebras. Given that the $X$'s and $Y$'s have been segregated we can, without any loss of information, merge all the $X$'s into one point and all the $Y$'s into another point. Here is what  $[[X,[X,Y]] , [X,Y]]$ looks like after this procedure has been  applied to the graph in Figure \ref{cbg03}(b)).

\begin{figure}[H]
\centering{\epsfig{file=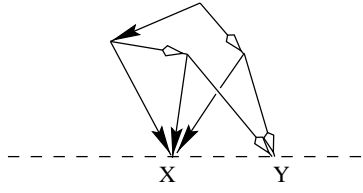 , height=25mm , clip= }}
\caption{Final form of $[[X,[X,Y]] , [X,Y]]$}
\label{cbg05}
\end{figure}

\no Let us summarize the method of assigning graphs to Lie elements in $\fg_F$  that we have just outlined. We are abusing nomenclature slightly at this point since by Lie element we mean elements of $\fg_F \smallsetminus \{X \cup Y\}$, \ie, all elements other than the generators $X$ and $Y$. 

\newcounter{l}
\begin{list}
{\arabic{l}.}{\usecounter{l}}
\item $[A,B]$ is to be represented by {\epsfig{file=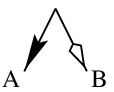 , height=10mm , clip= }}.
\item  Given a Lie element in $\fg_F$ we glue ordered wedges to form a rooted binary tree that represents that element.
\item We then pull all the $X$ leaves into one point and all the $Y$ leaves into another point.
\item The resulting graph is placed in the upper half plane with the $X$'s at 0 and the $Y$'s at 1 (we add this last directive in anticipation of the next section).
\end{list}

\begin{dfn}
A graph assigned to a Lie element $A \in \fg_F $ using the procedure outlined in
steps 1 through 4 above will be called the {\bf L-graph} of $A$ and will be denoted by $\Gamma_A$. Graphs which are L-graphs will be called {\bf Lie-admissible} graphs.
\end{dfn}

\no  We also fix some terminology that will be needed later. Vertices of L-graphs other than $0$ and $1$ will be called {\bf aerial} vertices. Wedges that land entirely on aerial vertices will be called aerial wedges. $0$ and $1$, often denoted by $X$ and $Y$ respectively, will be called the {\bf ground} vertices. Wedges that land on $X$ {\bf or} $Y$ will be called grounded wedges, those landing on   $X$ {\bf and} $Y$ being termed totally grounded wedges.

\no {\bf Note:} Aerial vertices are required to lie in the strict upper half plane \ech. \ie they cannot lie on $\br \subset \bar{\ech}$.

\vv 

\no The reason for introducing these more specialized(?)  graphs (instead of the traditional rooted trees) to represent Lie elements is simple. Not only do they represent Lie elements but they also  represent bi-differential operators when we color the graphs as in Figure \ref{3p7b} in Section \ref{secflat}. Recall that in Section \ref{cbhsec} when we were constructing the CBH-quantization we assigned  a bi-differential operator (using the notion of symbols)   to each bracketed Lie element in the Hausdorff series\footnote{the Hausdorff series we introduced in Section \ref{cbhsec} is best interpreted as belonging in $\fg_F$, though at that time we did not mention it.}. That correspondence can now be made transparent. It is nothing other than interpreting a given L-graph as Lie bracketing on one hand and as a bi-differential operator (colored by the given Poisson structure) on the other.

\vv

\no In the CBH-quantization there was one more step in addition to the replacement of the Hausdorff series by bi-differential operators - taking the exponential. This required taking products of bi-differential operators, but in a special way. As stated in Equation \eqref{con1} in Section \ref{consol}, this involved multiplying them as if they were constant coefficient operators. There is a simple transcription of this multiplication in the L-graph setting and we now describe that.

\vv

\begin{dfn}
\label{defspo} If $\Gamma_A$ and $\Gamma_B$ are two Lie-admissible graphs representing Lie algebra elements $A$ and $B$ respectively, then the composite graph formed by identifying their $X$ and $Y$ vertices is  denoted by $\Gamma_{AB}$
 and called the product of  $\Gamma_A$ and $\Gamma_B$.
\end{dfn}

\no Pictorially,

\begin{figure}[H]
\centering{\epsfig{file=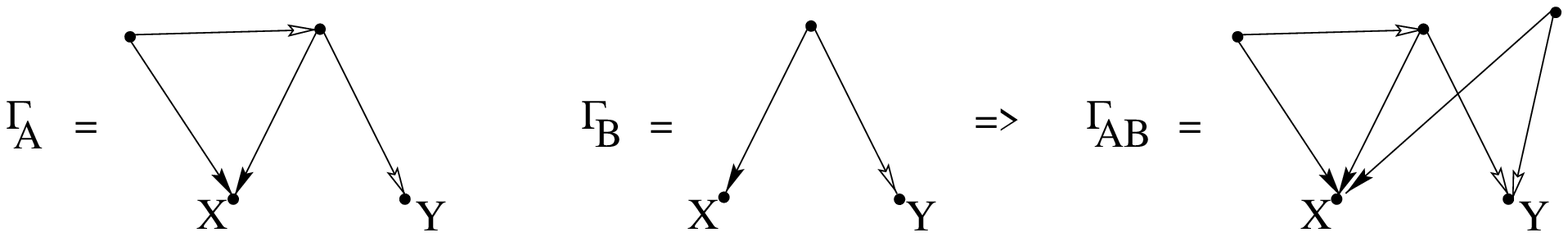 , height=2cm , clip= }}
\caption{L-graph multiplication}
\label{cbg4}
\end{figure}

\vv

\begin{figure}[H]
\centering{\epsfig{file=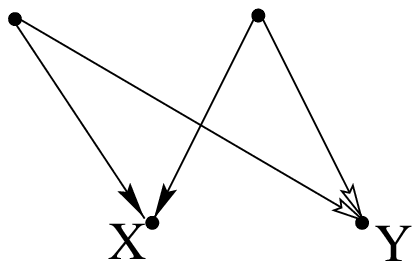 , height=2cm , clip= }}
\caption{$[X,Y][X,Y]$}
\label{cbg3m}
\end{figure}

\no  Apart from giving the required product for bi-differential operators, when we think of the L-graphs as Lie elements this product translates into the product structure on the symmetric algebra $S(\fg_F \smallsetminus \{X \cup Y\})$. This agrees with our imagery that two points, in this case two wedges, that are placed adjacent to each other are being multiplied. The product is clearly commutative and associative. As an example, the `element' encoded in Figure \ref{cbg3m} below is $[X,Y][X,Y] \in S(\fg_F)$.

\vv

\no One can now define the larger category of {\bf Sym-admissible graphs}, namely 
the semi-group generated by the Lie-admissible graphs under the product defined in Definition \ref{defspo}. The semi-group becomes a monoid when we introduce the identity element given by the graph with no edges and only the vertices $X$ and $Y$. 

\vv

\no There is a simple topological characterization
which determines the Sym-admissible graphs. Whereas the Lie-admissible graphs arose by taking a single rooted, binary tree and merging its leaves into two points, for the Sym-admissible graphs we are allowed to start with more than one such tree. Note that the Sym-admissible graphs account for all the terms that arise in the CBH-quantization.

\vv

\no We wish to point out that not every graph formed by superposing wedges yields a valid $S(\fg_F)$ element. Here are two instances of {\bf non} Sym-admissible graphs formed by superposing wedges.

\begin{figure}[H]
\centering{\epsfig{file=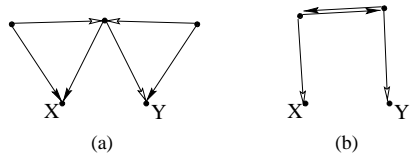 , height=2cm , clip= }}
\caption{Non Lie-admissible graphs}
\label{cbg3}
\end{figure}

\vv

\no Figures \ref{cbg3}(a) and \ref{cbg3}(b) are not Sym-admissible as they do not arise from binary trees. Whereas  Figure \ref{cbg3}(a) still arises from a tree (we do not consider the loops created by the collapsing of the $X$ and $Y$ leaves to be `true' loops),
 Figure \ref{cbg3}(b) has a non-trivial loop in it. As we will encounter such graphs later let us study them for a while.

\no  Figure \ref{cbg3}(a) suggests simultaneously bracketing more than one element with $[X,Y]$, \ie  something like $[ X, Y, [X,Y] ]$ which obviously does not make sense in the Lie algebra setting. Interpreted as a bi-differential operator this term {\em does} make sense, but it vanishes in the Lie algebra setting as we are taking two derivatives of a linear function (the coefficient of the linear Poisson structure) at the aerial vertex of $[X,Y]$. Thus terms such as this could arise in our iterative symmetrizing approach if we were dealing with a quadratic Poisson structure for example.

\vv

\no Figure \ref{cbg3}(b) is harder(?) to make sense of, though as a bi-differential operator it has a straightforward interpretation. Coloring the wedges with the Poisson structure, one finds that it represents the symmetric bi-differential operator that corresponds to the Killing form. A term such as this one could not have arisen in our iterative symmetrizing approach, as in that case we were only merging ground vertices of tri-differential operators to obtain new bi-differential operators, and that would never create an aerial loop. Since Kontsevich's approach allows for more general merging of vertices, such terms do arise in his product. 

\vv

\no Let us summarize the structure we have developed this far. We are interested in studying graphs that can be formed by gluing wedges together such that, suitably colored, they lead to bi-differential operators. 

\vv

\no One meaningful sub-category is that of Sym-admissible graphs. In relation to Lie algebras, they represent elements in $S(\fg_F)$. As bi-differential operators, they account for all the operators that arise in the CBH-quantization. They are obtained by merging leaves of binary rooted trees. The number of roots gives the degree of the corresponding symmetric polynomial in $S(\fg_F)$. Thus, Lie algebra elements are those with a single root. Two necessary conditions satisfied by Sym-admissible graphs that will help us recognize them are that 

(a) they cannot have more than one foot landing on an aerial vertex, and 

(b) they cannot contain loops. 

\no 

\vv {\bf Note:} As mentioned earlier, apparent `loops' at ground vertices are not considered loops. \eg, there is no loop in the following graph.

\begin{figure}[H]
\centering{\epsfig{file=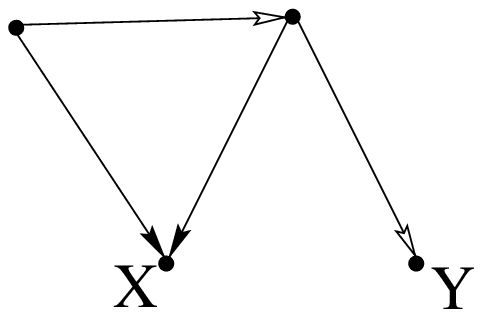 , height=2cm , clip= } }
\caption{$[X,[X,Y]]$}
\label{cbg2}
\end{figure}

\no There are two ways to reconcile this. One, we use the word loops to describe only those loops that may have existed before merging ground vertices. By definition  these  are ruled out for Sym-admissible graphs since we start with trees. Secondly, given that we are using directed graphs, the apparent loops are not valid loops when one takes direction into consideration. We are discussing these issues here since they will resurface in  Kontsevich's construction.

\subsection{Prolog to Kontsevich's Construction}

\no In the earlier subsection we have only identified some bi-differential operators that may arise in {\stp}s, without specifying the exact numerical coefficients that would guarantee associativity. In the CBH-quantization, the required numerical coefficients were given by the CBH formula. In our iterative symmetrizing in Section \ref{secflat}, we started with the requirement that we have a $\frac{1}{2}$ scaling the Poisson bracket and then iteratively got further coefficients. One can pose the following question. Is there a natural way of assigning numerical weights to graphs such that taking a weighted sum over all (suitably colored) graphs will yield a (bi-differential operator for a) \stp?

\vv

\no As stated this is an ill-posed problem since the numerical coefficients can be absorbed into the weights as well as into the coloring. Also, summing over {\em all} graphs could be an onerous task. One could address these issues by fixing a certain  class of graphs and a specific coloring scheme. Such a procedure is contained in Kontsevich's paper \cite{ko:defq}. We briefly outline his approach at this stage.

\vv

\no Kontsevich first introduces a class of admissible graphs that encode all bi-differential operators constructible from wedges. He also requires that the edges and (aerial) vertices be ordered. This determines the number of graphs that lead,  up to sign, to the same bi-differential operator. The coloring for an ordered wedge is the same as we had introduced in Section \ref{secflat} (Figure \ref{3p7b}). Namely one colors the aerial vertex of wedge with $\frac{1}{2}\alpha^{ij}$.\footnote{A fleeting perusal of \cite{ko:defq} may suggest that Kontsevich uses $\alpha^{ij}$ for his coloring. But his $\alpha^{ij}$ is our $\frac{1}{2}\alpha^{ij}$.} Having fixed all this, he defines  {\em natural} weights for all such graphs, leading to a \stp.

\vv

\no As an illustration, let us discuss the weight for the basic wedge 
\begin{figure}[H]
\centering{\epsfig{file=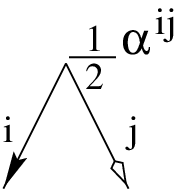, height=2cm , clip= }}
\caption{$\frac{1}{2}\{\ ,\ \}$}
\label{base}
\end{figure}

\no Despite the $\frac{1}{2}$ being absorbed into the coloring of the wedge, it turns out that in Kontsevich's scheme its weight is not $1$ but $\frac{1}{2}$. This is because there are {\bf two} ordered wedges whose contributions have to be added, and taking signs into consideration requires that the weight of {\epsfig{file=3p7a.eps , height=8mm , clip= }}  be $\frac{1}{2}$. We will not go into the details of the computation here, but postpone them to the next section.

\vv

\no The actual weights are given by integrals over certain configuration spaces
determined by  the graphs in question. The merging of points used in Section \ref{secflat} to obtain correction terms corresponds to moving to the boundary in these spaces. The vanishing of various `combinations' of weights (signed sums of various pairwise products of weights) yields associativity. 
Kontsevich uses Stokes' theorem to translate this vanishing condition to the vanishing of certain integrals on the boundary strata. He then does a case by case analysis to show that these integrals vanish. 

\vv

\no Having introduced Kontsevich's weights we mention that they are intimately
related to the coefficients arising in the CBH formula. Recall that all the bi-differential operators for the CBH-quantization arise from Sym-admissible graphs. The bi-differential operators corresponding to   Sym-admissible graphs in Kontsevich's scheme  match exactly with those in the CBH quantization. 
This is one of the key results of this paper. Stated alternately, Kontsevich's quantization applied to the linear Poisson case is none other than the CBH-quantization,  with some additional terms which are inessential for associativity in the linear setting. Given the universality of his formula (it works for {\em all} Poisson structures on ${\br}^d$) it is not surprising that Kontsevich's formulation may not  yield a `minimal' quantization for linear Poisson structures. 
We will address the correspondence between the Kontsevich and CBH  quantizations in Section 5.

\section{Kontsevich's Construction}
\no Kontsevich's \stp can be described in a fairly straightforward manner. The following subsection largely regurgitates pages 5 and 6 of \cite{ko:defq} with slight notational modification. Wherever appropriate we make some comparative remarks.
\subsection{The Formula} 

\no $G_n,\ n \geq 0$ will denote the class of all {\bf admissible} graphs. This is a special class of oriented, labeled graphs which accounts for all bi-differential operators that can be  constructed from $n$ wedges.

\begin{dfn}
\label{defff}
An oriented graph $\Gamma$ is a pair ($V_\Gamma, E_\Gamma$) of two finite sets such that $E_\Gamma$ is a subset of $V_\Gamma \times V_\Gamma$.
\end{dfn}
\no Elements of $V_\Gamma$ are the vertices of  $\Gamma$, and elements of $E_\Gamma$ are its edges. For $e = (v_1, v_2) \in E_\Gamma$ we say that $e$ starts at $v_1$ and ends at $v_2$.

\no We say that a labeled graph $\Gamma$ belongs to $G_n$ if
\newcounter{k}
\begin{list}
{(\arabic{k})}{\usecounter{k}}
\item $\Gamma$ has $n+2$ vertices and $2n$ edges.
\item The set of vertices $V_\Gamma$ is given by $\{1,2,\ldots,n\} \cup \{X, Y\}$ where $X$ and $Y$ are just two symbols.
\item For every $k \in \{1,2,\ldots,n\}$ there are two edges starting at $k$. These are ordered and are labeled by $e_k^1$ and $e_k^2$.
\item The edges emanating from $v$ can land on any vertex other than $v$ itself. \ie for every $v \in V_\Gamma$ the ordered pair $(v, v) \not\in  E_\Gamma$.
\item For $n \geq 1, \quad G_n$ has $(n(n+1))^n$ elements, and $G_0$ has a single element.
\end{list}

\vv

\no {\bf Note:} Any one of the  vertices  $\{1,2,\ldots,n\}$ and its two emanating edges correspond to our wedges. Following the terminology of Chapter \ref{chp3} the vertices  $\{1,2,\ldots,n\}$ will be called {\bf aerial} vertices. Kontsevich uses $L$ and $R$ (for Left and Right) instead of $X$ and $Y$ for his {\bf ground} vertices.
\vv 

\no Next, to each labeled graph $\Gamma \in G_n$ we associate a bi-differential operator $B_{\Gamma,\bar{\alpha}}$,

\begin{equation}
\label{kodef1}
B_{\Gamma,\bar{\alpha}} : \cala \times \cala \to \cala , \quad \cala = C^{\infty}(\rd) \; ,
\end{equation}

\no where $\bar{\alpha}$ is a bi-vector field on $\rd$, though not necessarily a Poisson one. Before we present the rule for constructing $B_{\Gamma,\bar{\alpha}}$ we wish to point out that it boils down to coloring the constituent wedges of $\Gamma$ with $\bar{\alpha}$ in all possible ways and then summing over them. Kontsevich uses $\alpha$ in his definition but since his 
$\alpha^{ij}$ corresponds to our $\frac{1}{2}\alpha^{ij}$ we have switched to $\bar{\alpha}$ in this section. Here is Kontsevich's example explaining the coloring scheme. Consider $\Gamma \in G_3$ with its edges given by
\begin{equation}
\label{kodef2}
(e_1^1, e_1^2, e_2^1, e_2^2, e_3^1, e_3^2) = ( (1,X),(1,Y),(2,Y),(2,3),(3,X),(3,Y) ) \; .
\end{equation}
\begin{figure}[H]
\centering{\epsfig{file=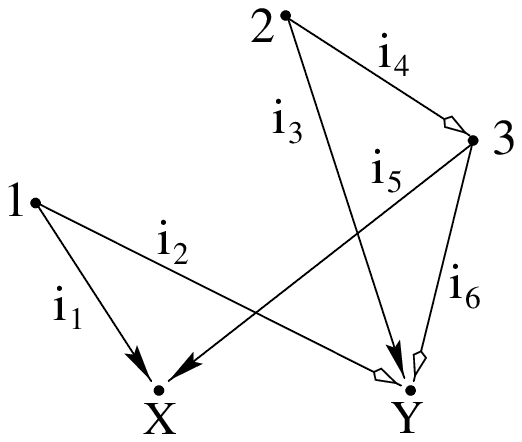 , height=4cm , clip= }}
\caption{$\Gamma \in G_3$}
\label{kopic1}
\end{figure}

\no In his picture Kontsevich uses the same arrowhead for all the edges whereas we have retained our earlier convention of using distinguishable arrowheads with the \epsfig{file=leg1.eps, height=6mm, clip= } denoting the first edge ($e_*^1$) for a given wedge. Independent indices $1\leq i_1,i_2,\ldots,i_6 \leq d$ have been placed along the edges instead of labels $e_m^n$ to facilitate the description of  $B_{\Gamma,\bar{\alpha}}$. It is given by 
\begin{equation}
\label{kodef3}
(f,g) \mapsto \sum_{i_1,\ldots,i_6} {\bar{\alpha}}^{{i_1}{i_2}} {\bar{\alpha}}^{{i_3}{i_4}}\del{{i_4}}( {\bar{\alpha}}^{{i_5}{i_6}})\del{{i_1}}\del{{i_5}}(f)\del{{i_2}}\del{{i_3}}\del{{i_6}}(g) \; .
\end{equation}
\no The general formula for  $B_{\Gamma,\bar{\alpha}}$ is as follows 
\begin{equation}
\label{kodef4}
\begin{split}
B_{\Gamma,\bar{\alpha}}(f,g) \dfequal &\sum_{I:E_\Gamma \to \{1,\ldots,d\}} 
\Biggl[\prod_{k=1}^n\biggl(\prod_{e\in E_\Gamma, e=(*,k)} \del{I(e)}\biggr){\bar{\alpha}}^{{I(e_k^1)}{I(e_k^2)}}\Biggr] \times \\
&\times \biggl(\prod_{e\in E_\Gamma, e=(*,X)} \del{I(e)}\biggr)f \times
\biggl(\prod_{e\in E_\Gamma, e=(*,Y)} \del{I(e)}\biggr)g \; . \\
\end{split}
\end{equation}
\no The map $I:E_\Gamma \to \{1,2,\ldots,d\}$ is the same as our coloring map for edges. Since we are only using  the ${\bar{\alpha}}^{ij}$ to color vertices, the map $I$ also helps define the coloring for the vertices. 

\vv

\no The next step consists of associating a numerical weight $w_K(\Gamma)$ with each graph $\Gamma \in G_n$. Kontsevich uses the notation $w_\Gamma$ but as we will be introducing some additional modified weights we choose  $w_K(\Gamma)$ for his weights. The weights are constructed using a special angular measure defined on a punctured upper half-plane.

\vv

\begin{figure}[H]
\centering{\epsfig{file=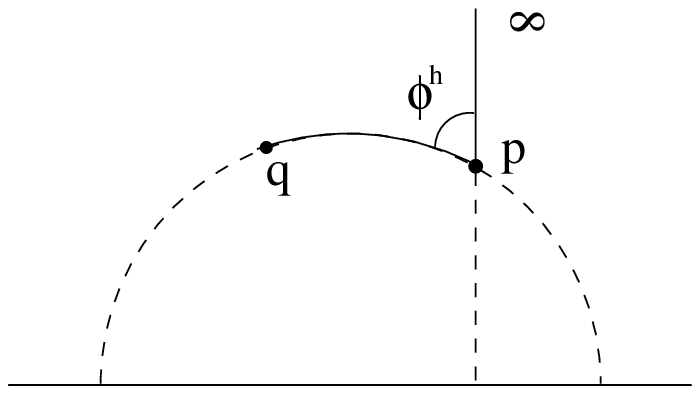 , height=3cm , clip= }}
\caption{The harmonic angle $\phi^h$}
\label{kopic2}
\end{figure}

\no Let \ech denote the standard upper half plane with its constant negative curvature metric. Recall that geodesics in \ech are given by vertical straight lines and semi-circles centered on $\br \subset \bech$. Given $p,q \in \ech, \quad p \ne q,\quad $ we define $\; \phi^h(p,q) \in \br/2\pi\bz \quad $ to be the angle between the two directed geodesics $l(p,\infty)$ and $l(p,q)$. As displayed in Figure \ref{kopic2} $\quad l(p,\infty)$ is the vertical line through $p$ and the angle is measured counterclockwise from $l(p,\infty)$ to $l(p,q)$.
 
\no Though based at $p$ in Figure \ref{kopic2} the angle $\phi^h(p,q)$ provides
the notion of an angle subtended by $p$ at $q$. 
\begin{figure}[H]
\centering{\epsfig{file=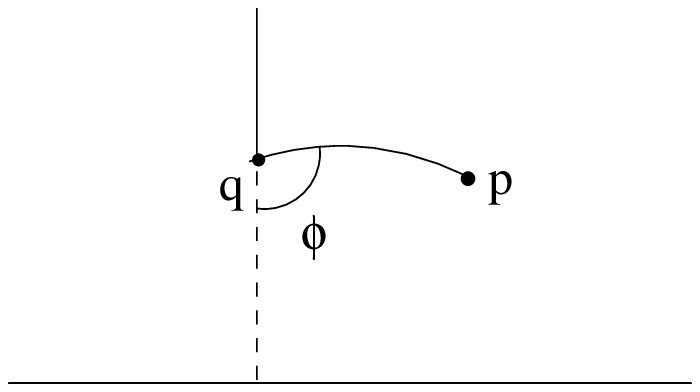 , height=3cm , clip= }}
\caption{The angle subtended by $p$ at $q$}
\label{kopica}
\end{figure}
 
\no This will become clear when we study the level curves of $\phi^h(p,q)$ in the next subsection. The $h$ in $\phi^h$ stands for {\em harmonic}. Observe that  $\phi^h(p,q)$ extends (by continuity) to $p,q \in \bech$. Furthermore  $\phi^h(x,q) = 0$ for $x \in \br$. \ie, points on $\br \subset \ech$ subtend a zero angle at any point $q \in \ech$. It is this property that makes this angle measure quite different from the standard angle measure on the plane. 

\no Next, by $\ech_n$ we denote the space of configurations of $n$ numbered pairwise distinct points on $\ech$.
\begin{equation}
\label{kodef5}
\ech_n = \{(p_1,p_2,\ldots,p_n): p_k \in \ech, p_k \ne p_l \; \mbox{for} \; k \ne l \} \; .
\end{equation}
\no $\ech_n \subset \bc^n$ is a non-compact smooth $2n$-dimensional manifold and inherits an orientation via its natural complex structure. Now given a graph $\Gamma \in G_n$ and a configuration of points $(p_1,p_2,\ldots,p_n) \in \ech_n$ we can draw a representative for $\Gamma$ in $\ech$ by placing the $k$th vertex at $p_k$, $1 \leq k \leq n$, the vertex $X$ at $0$, $Y$ at $1$,
and using directed geodesics as edges. \eg for $\Gamma =$\epsfig{file=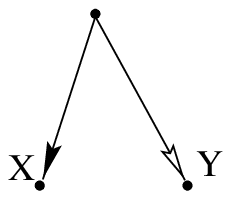 , height=8mm , clip= }, and $p$ any point in \ech we get a figure of the form
\begin{figure}[H]
\centering{\epsfig{file=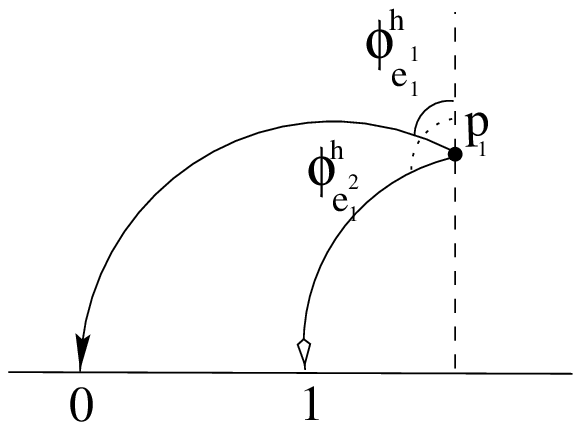 , height=3cm , clip= }}
\caption{A graph in \ech}
\label{kopic3}
\end{figure}
 
\no For the sake of clarity we will revert to drawing graphs with straight lines despite the hyperbolic picture. As the $p_i$ vary over $\ech_n$, each edge $e$ joining two vertices, say  $p$ and $q$, defines a function $\phi_e^h:\ech_n \to \br/2\pi\bz$, obtained in the following way. We identify $p$ and $q$ with their images under the natural projection from $\ech_n$ to $\ech$ and $\phi_e^h$ is then the pullback of $\phi^h(p,q)$. The weight $w_K(\Gamma)$ is then defined as
\begin{equation}
\label{kodef6}
w_K(\Gamma) \dfequal \frac{1}{n!(2\pi)^{2n}} \int_{\ech_n} \bigwedge_{i=1}^{n} 
(d\phi_{e_k^1}^h \wedge d\phi_{e_k^2}^h) \; .
\end{equation}

\vv

\no Note that in the restricted situations where we will be able to compute weights, we will be setting up iterated integrals by fixing all but one $p_k$ and integrating out its contribution. In this setting we will abuse notation and still use $\phi_e^h$ to denote the function from $\ech \smallsetminus \{q\}$ to $\br/2\pi\bz$ given by $p_k \mapsto \phi^h(p_k,q)$.

\begin{dfn}{\bf Kontsevich's \stp.}

\no Given the above weights $w_K(\Gamma)$ and a Poisson bi-vector field $\bar{\alpha}$ on $\rd$, the formula
\begin{equation}
\label{kodef7}
\notag
f \star g \dfequal \smo \sum_{\Gamma \in G_j} w_K(\Gamma) B_{\Gamma,\bar{\alpha}}(f,g)
\end{equation}
defines a \stp. We will refer to this as the Kontsevich quantization.
\end{dfn}

\no On occasion it will be simpler to consider the weights without the $n!$ factor in the denominator. With that in mind we define the modified  weight system
$w_I(\Gamma)$
\begin{equation}
\label{4ee2}
w_I(\Gamma) \dfequal  \frac{1}{(2\pi)^{2n}}  \int_{\ech_n}  \bigwedge_{i=1}^{n} (d\phi_{e_k^1}^h \wedge d\phi_{e_k^2}^h) \; .
\end{equation}
\no In his paper \cite{ko:defq} Kontsevich explains the general setting that leads to the above weights. The punctured upper half plane $\ech \smallsetminus \{q\}$ provides the prototype substrata in the configuration spaces attached to admissible graphs. The harmonic angle introduced above is just a special example of more general angle measures (defined on $\ech \smallsetminus \{q\}$)
that would work equally well to give the weights. In Section \ref{secang} we will briefly describe the nature of this angular measure and rescale it to ease our computations. For now we gather some general facts about admissible graphs and their weights.

\subsection{Complexity Analysis}
\label{seccomp}
\no In this subsection we try to set up some structure theory for the admissible graphs. We do this by dividing them into certain natural groupings. This division is guided by the following considerations.
\newcounter{z}
\begin{list}
{(\arabic{z})}{\usecounter{z}}
\item We wish to have groupings that reflect the level of difficulty involved in computing the weights.
\item We wish to establish that Kontsevich's quantization is also of the form \\ $\exp$(symmetrized deformation of the Poisson structure).
\item Given our preoccupation with the linear setting we wish to concentrate on computing weights that are relevant for linear Poisson structures.
\end{list}

\subsubsection{Prime Graphs}
\label{secprim}
\no We have already defined two subcategories of the admissible graphs, Sym-admissible  and Lie-admisible graphs in Section \ref{freeg}. The first subcategory we introduce here is that of  prime graphs. Recall the graph multiplication that we had defined in Definition \ref{defspo}. It extends naturally to all admissible graphs. Namely, the product of two graphs is the graph obtained by merging their respective $X$ and $Y$ vertices. Pictorially,

\vv 

\begin{figure}[H]
\centering{\epsfig{file=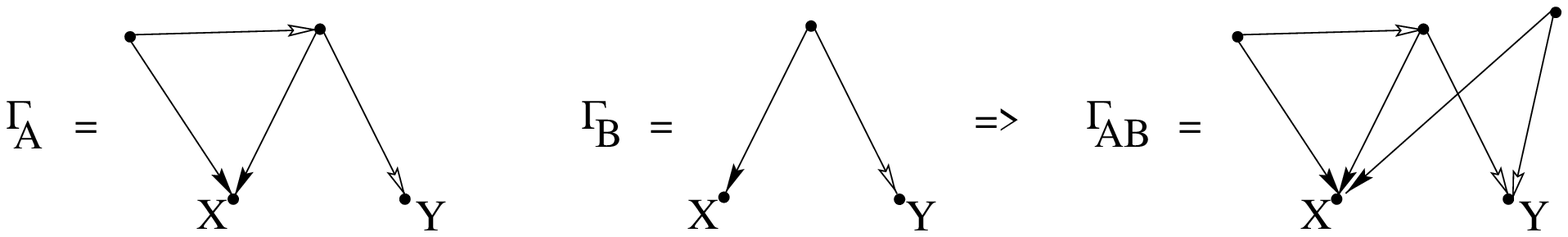 , height=2.2cm , clip= }}
\caption{Graph multiplication revisited}
\label{4p7}
\end{figure}

\no As before, the graph $\Gamma_0$ which has no edges and only the vertices $X$ and $Y$ is the identity for this multiplication. Admissible graphs form a monoid under this multiplication.

\begin{dfn}
\label{4d1}
An admissible graph $\Gamma$ will be called {\bf prime} if it has no factors other than 
$\Gamma_0$ and itself (under the above multiplication).

\no Graphs that are not prime will be called composite.
\end{dfn}

\no The notion of being prime is more or less equivalent to being connected. \ie, if one separates out all the edges (feet) of a graph $\Gamma$ that land at $X$ and $Y$, and if $\Gamma$ is still connected, then it is prime. If it becomes disconnected it was composite.  Note that the Lie-admissible graphs are all prime.

\vv

\no The reason behind defining this category  is that the (modified) weights $w_I$ are multiplicative. \ie, $w_I(\Gamma_{AB}) = w_I(\Gamma_A) w_I(\Gamma_B)$. This follows from the following rudimentary property of iterated integrals.
\begin{equation}
\label{4e3}
\iint f(x) g(y) dy dx = \int f(x) dx \int g(y) dy \; .
\end{equation} 

\no Thus knowing the weights of the prime graphs determines all the weights. Further, when we work with the $w_K$(\ie introduce the $n!$ in the denominator) and take all contributions into account, then we will see that composite graphs in Kontsevich's formula arise precisely as mandated by  the expansion of 

\begin{center}
$\exp$ (sum of prime graphs). 
\end{center}
\subsubsection{Loop Graphs}

\no Next we introduce the subcategory of loop graphs. 
\begin{dfn}
Graphs which contain a valid loop as a subgraph will be called {\bf loop graphs}. 
\end{dfn}

\no By valid loop we mean that one is only allowed to consider loops formed when one moves in accordance with the {\em directed }
edges of the graph. Kontsevich uses the term `wheels' for loops. Here are examples of loop graphs. We make the trivial observation that they may be prime or composite.
\begin{figure}[H]
\centering{\epsfig{file=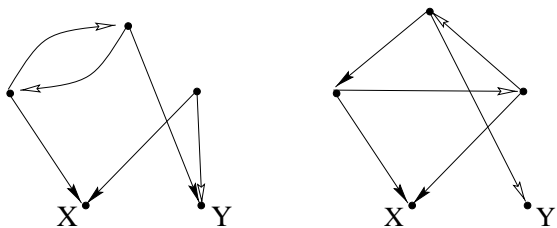 , height=3cm , clip= }}
\caption{Loop graphs}
\label{4p8}
\end{figure}

\no Here are graphs that are {\bf not} loop graphs.
\begin{figure}[H]
\centering{\epsfig{file=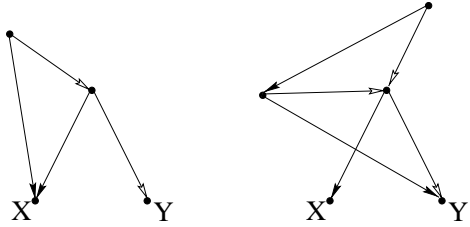 , height=3cm , clip= }}
\caption{Non-loop graphs}
\label{4p9}
\end{figure}

\no We separate out the loop graphs since we are unable to set up computable iterated integrals for their weights or compute them by other means. 

\no {\bf Note:} We do not claim that this cannot be done, just that we are unable to do so at present. 

\no We are still able to relate Kontsevich's product to the CBH-quantization because irrespective of their weights, in the linear setting,  the bi-differential operators corresponding to the loop graphs turn out to be unnecessary for associativity and can be subtracted out. We will address this issue in the next section.
For now we simply choose to ignore the presence of the loop graphs. For the remainder of this section, unless stated otherwise, graphs will mean non-loop graphs.

\no Even for the non-loop prime graphs we are unable to explicitly compute all the weights. While one can categorize the varying levels of complexity, we start with the simplest category, one we will call the W-computable graphs (W for weight).

\subsubsection{W-computable Graphs}

\no We start out by displaying the admissible graphs whose weights we can compute in a straightforward manner by actual integration. They are of the form

\begin{figure}[H]
\centering{\epsfig{file=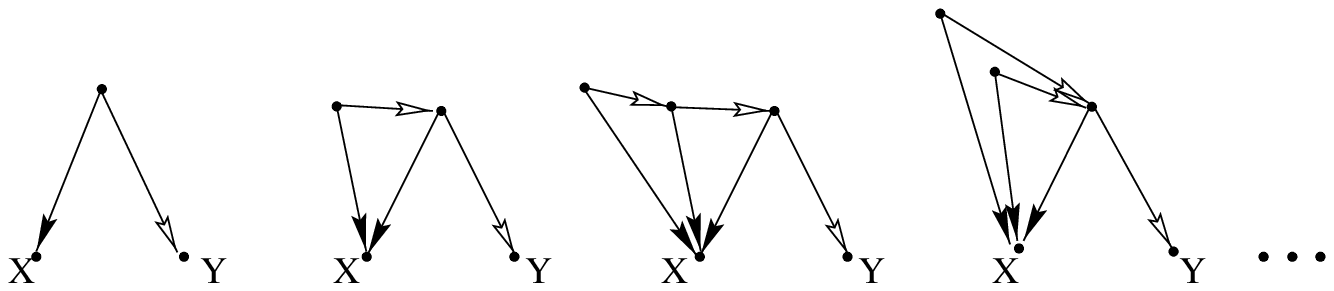 , height=3cm , clip= }}
\caption{Some W-Computable graphs}
\label{4cpic1}
\end{figure}

\no Though the general form may be clear from the above picture we will provide a rigorous definition below. Note that although we have drawn all the above graphs with the black feet at $X$, graphs obtained by reordering the edges and vertices (while retaining the topological structure of the graph) have equally computable weights. In fact their weights, if at all different, differ in sign alone. This is because reordering edges simply changes the ordering of the 1-forms in the weight integral (Equation \eqref{kodef6}). Furthermore weights of graphs obtained by mirror reflecting those in Figure \ref{4cpic1}  about a vertical line ($x=\frac{1}{2}\; \mbox{in}\; \ech$) are equally computable. Pictorially we are referring to the following graphs.

\begin{figure}[H]
\centering{\epsfig{file=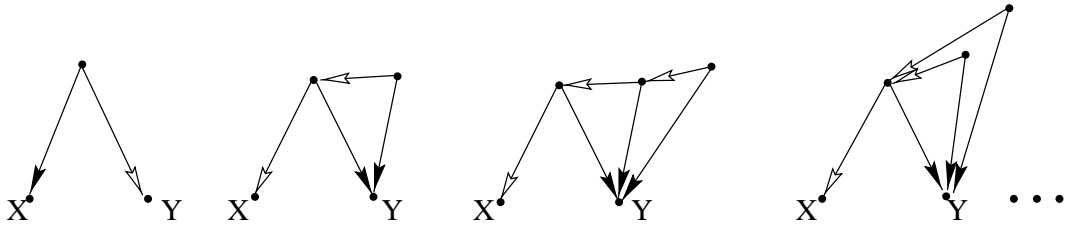 , height=3cm , clip= }}
\caption{Some more W-computable graphs}
\label{4cpic2}
\end{figure}

\no The above equivalence of weights follows from the (anti)symmetry of the angle measure about a vertical line. We now give a concrete description for these types of graphs. Given the preceding statements we will only concern ourselves with graphs such as those displayed in Figure \ref{4cpic1}, the mirrored ones being obtainable in a similar fashion.

\no Let us denote the graph \epsfig{file=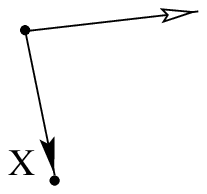 , height=6mm , clip= } by $\frac{1}{2}{\rm ad}_X$. Note that $\frac{1}{2}{\rm ad}_X$ is not an admissible graph but we shall use it to construct admissible ones. 
We define {\bf W-computable} graphs to be those admissible graphs that can be constructed
as follows:

\newcounter{pk}
\begin{list}
{W-(\arabic{pk})}{\usecounter{pk}}
\item The graph \epsfig{file=4p15.eps , height=8mm , clip= } is W-computable.
\item Admissible graphs formed by concatenating \epsfig{file=4cpic3.eps , height=8mm , clip=} with  a W-computable graph are W-computable.
\end{list}  

\no Here by concatenation of \epsfig{file=4cpic3.eps , height=6mm , clip= } with $\Gamma$ we mean the graph formed by attaching the $X$ vertex of \epsfig{file=4cpic3.eps , height=6mm , clip= } to the $X$ vertex of $\Gamma$, and allowing the other (free) foot of \epsfig{file=4cpic3.eps , height=6mm , clip= } to land on any {\bf aerial} vertex of $\Gamma$. \eg,

\begin{figure}[H]
\centering{\epsfig{file=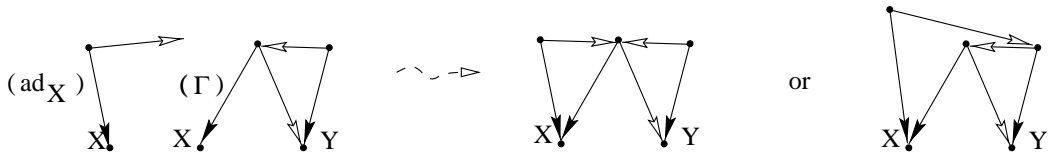 , height=2cm , clip= }}
\caption{Concatenation in action}
\label{4spic5}
\end{figure}
\no The graphs in Figure \ref{4cpic1} are precisely the first four W-computable graphs. The restriction that the free foot of  \epsfig{file=4cpic3.eps , height=8mm , clip= } only land on aerial vertices is an artificial one. \eg, the weight of the graph 
 \epsfig{file=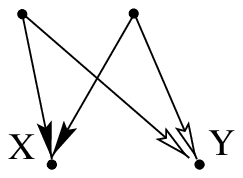 , height=8mm , clip= }
formed by 
 \epsfig{file=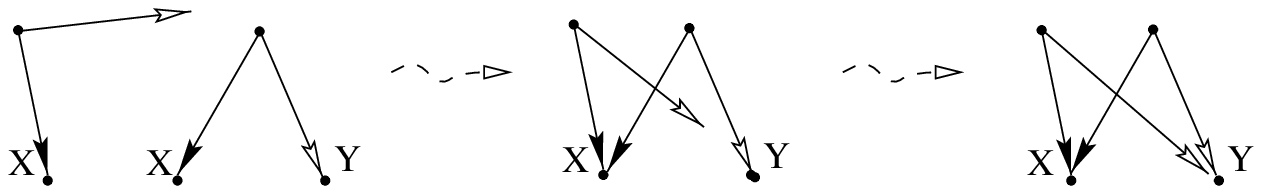 , height=8mm , clip= }
is equally easy to compute. The reason we have imposed this restriction is that under this restriction our procedure generates prime graphs alone. Allowing for the more general concatenation will not yielded any graphs whose weights cannot be obtained from the prime W-computable ones.

\vv
\no We wish to mention that there is a larger category that contains the W-computable graphs for which it is still straightforward to set up iterated integrals, though the computations may require symmetry considerations and/or recursive computations. One could call these the {\bf tractable} graphs. They are the ones formed by concatenating not just $\frac{1}{2}{\rm ad}_X$, but also $\frac{1}{2}{\rm ad}_Y$ (\epsfig{file=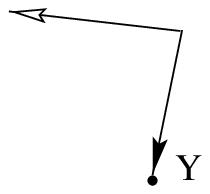 , height=8mm , clip= }). Here are some tractable graphs.
\begin{figure}[H]
\centering{\epsfig{file=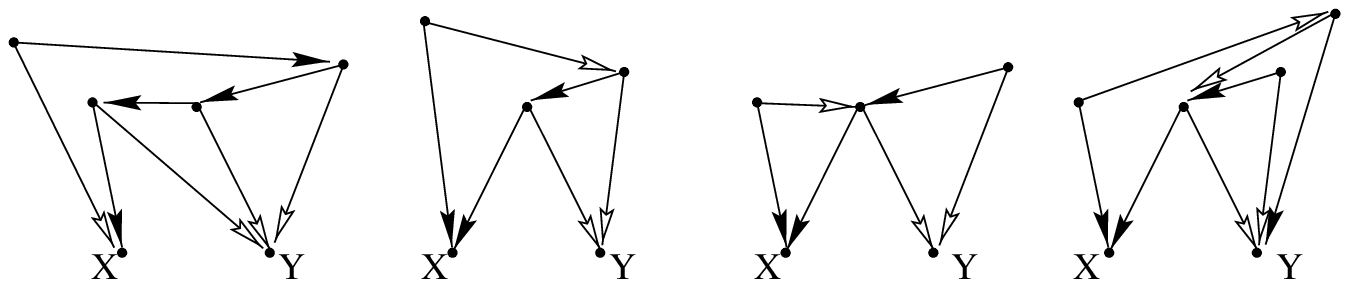 , height=3cm , clip= }}
\caption{Some tractable graphs}
\label{4p13}
\end{figure}

\no The basic point is that these graphs are made solely with grounded wedges, \ie , each constituent wedge of the graph has atleast one foot on the ground. We will not delve too much into the tractable graphs. There are  two reasons for this. One, the W-computable graphs will prove to be adequate for the purposes of relating the Kontsevich and CBH quantizations. Secondly, though we are able to recursively compute the low order tractable graphs, we have not been able to formalize the procedure. 

\no We now go about the task of computing the integrals for the W-computable graphs. For that we first make some modifications to the angle measure which provides the basis for the weights.

 \subsection{Angular Measure}
\label{secang}
The angular measure $\phi^h(p,q)$ can be replaced by any other $\phi(p,q)$ as long as it satisfies the following properties.
\newcounter{y}
\begin{list}
{A-\arabic{y}.}{\usecounter{y}}
\item $\phi(p,q)$ considered as a function of $p$ is a smooth $S^1$ valued function on $\bech \smallsetminus \{q\}$ for all $q \in \bech$;
\item $\phi (p,q), \quad p,q \in \bech$, should give a sensible notion for ``the angle subtended by $p$ at $q$,'' angles being measured counterclockwise from a vertical line through $q$. There is the additional requirement that points $x \in \br \subset \bech$ subtend a zero angle at any $q \in \ech$. It is this requirement that does not allow us to use the standard Euclidean angle measure.
\end{list}

\no For a precise formulation of the angle measure we refer the reader to Section 6.2 of \cite{ko:defq}. For our purposes it will be adequate to work with just a rescaled version of $\phi^h(p,q)$. The easiest way to see that ${\phi^{h}} (p,q)$ satisfies the above conditions is by drawing
level sets for ${\phi^{h}} (p,q)$, $q$ being held fixed. There are two typical cases. $q \in \ech$ and $q=x \in \br$.
\begin{figure}[H]
\centering{\epsfig{file=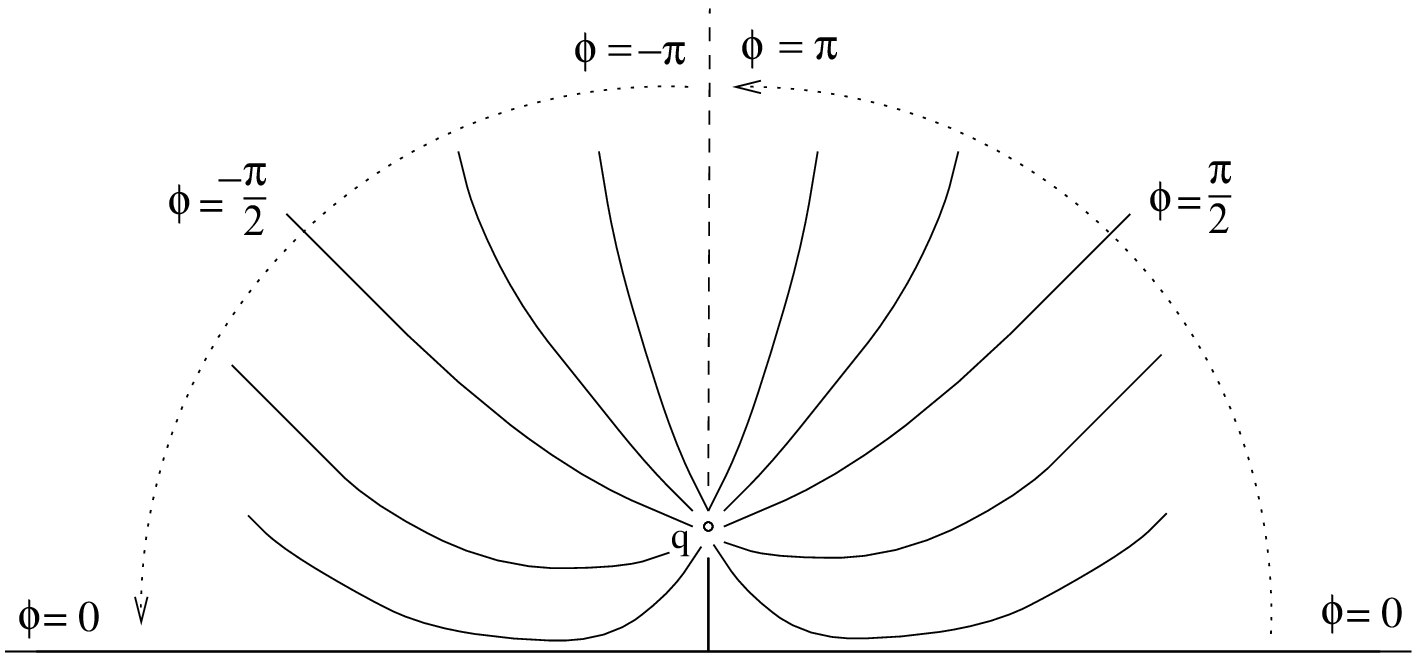 , height=4.5cm , clip= }}
\caption{Level curves for  ${\phi^{h}} (p,q)$, $q \in \ech$ fixed}
\label{4p1}
\end{figure}

\vv
\begin{figure}[H]
\centering{\epsfig{file=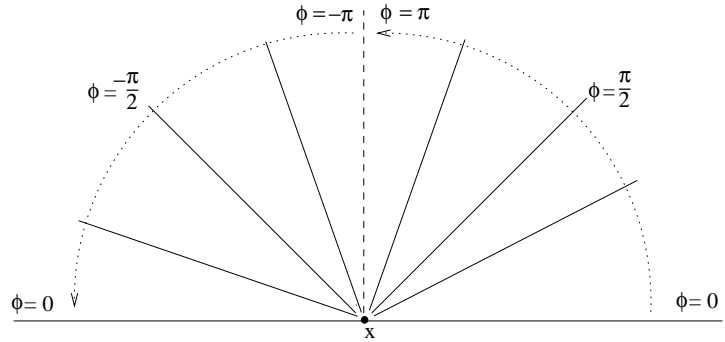 , height=4.5cm , clip= }}
\caption{Level curves for  ${\phi^{h}} (p,x)$, $x \in \br$ fixed}
 \label{4p2}
\end{figure}

\no The curves in Figure \ref{4p1} are hyperbolae. This follows in a straightforward manner from the fact that geodesics in \ech are given by semi-circles centered on the real axis. Observe that if we were to consider $\phi^h (p,q)$ as taking values in $\br$ then it is not defined on the vertical dashed line above $q$, where the function then changes value from $\pi$ to $-\pi$. We point this out because to set up our iterated integrals we will have to choose a specific {\em branch} for $\phi^h (p,q)$ and the integrals are quite complicated if we work with $\phi^h (p,q)$ taking values in $\br$ as described in Figures \ref{4p1} and \ref{4p2}. This can be easily remedied by taking the branch $\varphi(p,q)$ for  $\phi^h (p,q), \; \varphi(p,q)$ being  defined as follows.

\no Instead of going from $0$ to $\pi$ and $-\pi$ to $0$ as we move counterclockwise, $\varphi(p,q)$ will vary from $0$ to $2\pi$. We display this below.

\begin{figure}[H]
\centering{\epsfig{file=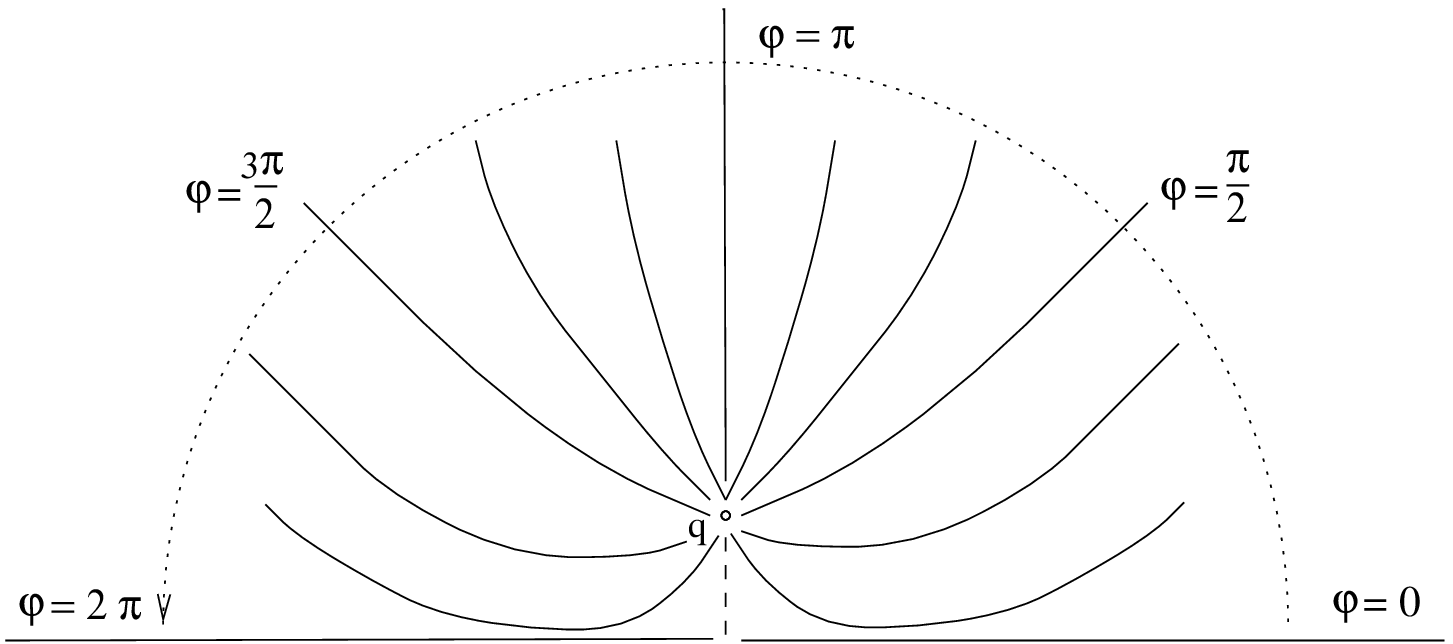 , height=4.5cm , clip= }}
\caption{Level curves for  $\varphi (p,q)$, $q \in \ech$ fixed}
\label{4p3}
\end{figure}

\no The corresponding alteration for Figure \ref{4p2} looks thus.

\begin{figure}[H]
\centering{\epsfig{file=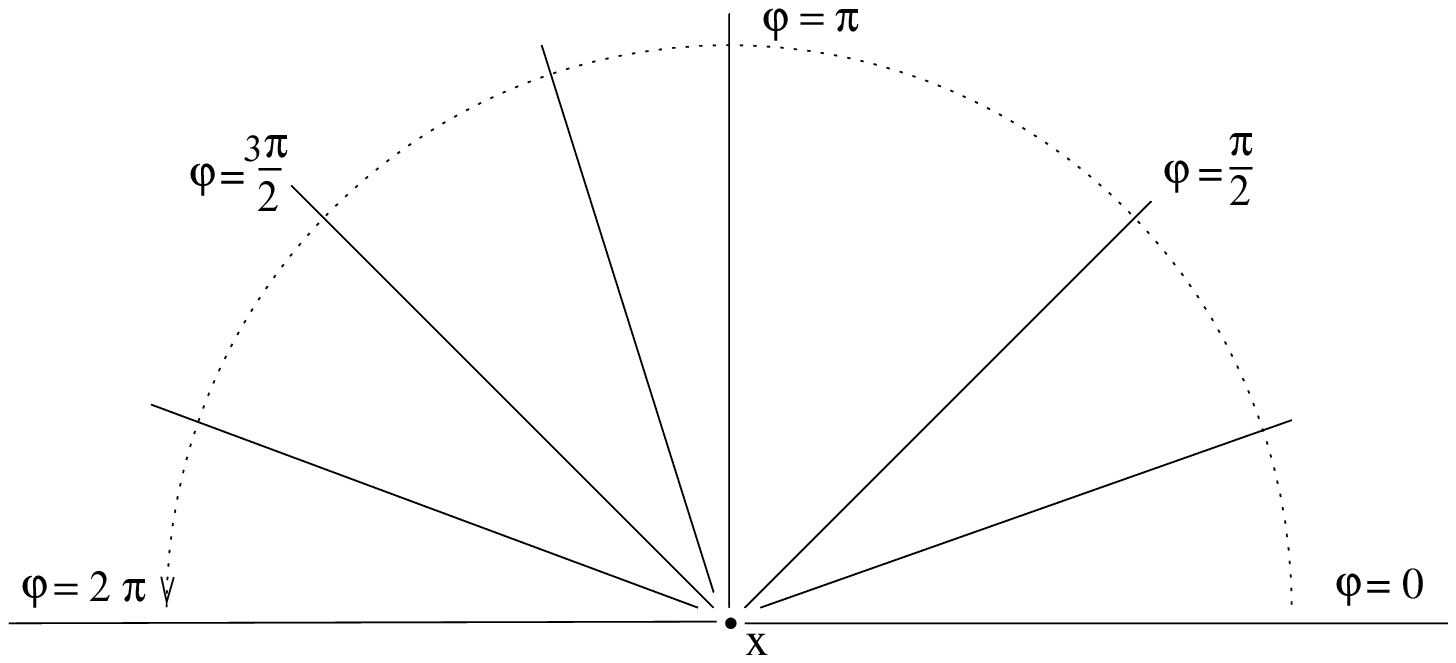 , height=4.5cm , clip= }}
\caption{Level curves for  $\varphi (p,x)$, $x \in \br$ fixed}
\label{4p4}
\end{figure}

\no Note that the discontinuity in Figure \ref{4p3} is now just a finite segment and there is none in Figure \ref{4p4}.

\no We make one more alteration to the functions describing the weight integrals before starting with our computations. Notice that in the definition of the weights (Equation \eqref{kodef6}) there is a factor of $(2\pi)^2$ in the denominator for every integral of type $ \int_{\ech} (d{\phi_{e_k^1}} \wedge d{\phi_{e_k^2}})$. This normalizing factor scales the angular measure from $2\pi$ to $1$. But we can effect this by simply rescaling $S^1$. \ie, requiring that our angle measure take values in $S^1 \cong \br/\bz$. This coordinate transformation will free us from carrying a whole bunch of $\pi$'s around and make it easier to recognize identities satisfied by the weights. We represent our final angle measuring function by $\lambda(p,q)$. For the sake of consolidation here is its  level curve diagram (we skip the second one for brevity).
\begin{figure}[H]
\centering{\epsfig{file=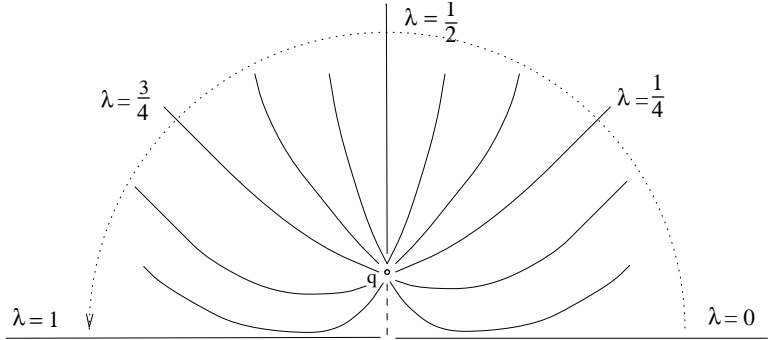 , height=4.5cm , clip= }}
\caption{Level curves for  $\lambda (p,q)$, $q \in \ech$ fixed}
\label{4p5}
\end{figure}

\no  Here is what the Kontsevich weight integrals look like when defined in terms of $\lambda \quad (\Gamma$ being an element of $G_n$).


\begin{equation}
\label{4e1}
w_K(\Gamma) \dfequal \frac{1}{n!} \int_{\ech_n} \bigwedge_{i=1}^{n} 
(d{\lambda_{e_k^1}} \wedge d{\lambda_{e_k^2}}) \; .
\end{equation}

\vv

\no The modified  weights $w_I(\Gamma)$ are given by
\begin{equation}
\label{4e2}
w_I(\Gamma) \dfequal  \int_{\ech_n} \bigwedge_{i=1}^{n} 
(d{\lambda_{e_k^1}} \wedge d{\lambda_{e_k^2}}) \; .
\end{equation}

\no In the next subsection we start integrating the weights for the W-computable graphs. The computations will use no structure beyond of basic multivariable calculus (we've already encountered a rescaling of coordinates).

\subsection{Weight Computation}
\label{4wc}
\no Let us start by computing the weight for the graph {\epsfig{file=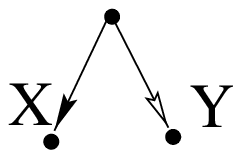 , height=8mm , clip= }}, a totally grounded wedge. We denote its weight $w_I({\epsfig{file=4p7a.eps , height=6mm , clip= }})$ by $w_1$. Along the way we will modify our notation and build up some further machinery to facilitate the integration of the weights. An imbedding of  {\epsfig{file=4p7a.eps , height=8mm , clip= }} in \ech is given by

\begin{figure}[H]
\centering{\epsfig{file=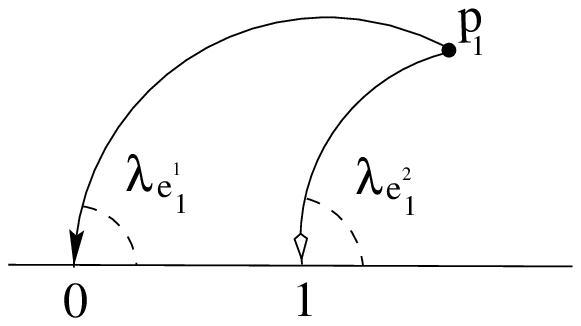 , height=3cm , clip= }}
\caption{A wedge in  $\ech$ }
\label{4wp2}
\end{figure}

\no Using Equation \eqref{4e2} we have 
\begin{equation}
\label{4w1}
w_1 =  \int_{\ech} d{\lambda_{e_1^1}} \wedge d{\lambda_{e_1^2}} \; ,
\end{equation}

\no where, as displayed in Figure \ref{4wp2}, $\lambda_{e_1^1}$ and $\lambda_{e_1^2}$ are the angles subtended by the vertex $p_1 \in \ech$ at $0$ and $1$ respectively. Though we keep using the word angle, we remind ourselves that the modified angle $\lambda$ varies between $0$ and  $1$.

\no {\bf Some notational modification:} Currently we are denoting the angle measure for the edges starting at vertex $p_k$ by $\lambda_{e_k^1}$ and $\lambda_{e_k^2}$. From now on, if an edge starting at $p_k$ ends at $p_i$ we will denote the corresponding angle measure (the angle subtended by $p_k$ at $p_i$) by \lak{i}. Further we set $\lak{X} = \theta_k$ and $\lak{Y} = \psi_k$. Thus for generic pictures we get 

\begin{figure}[H]
\centering{\epsfig{file=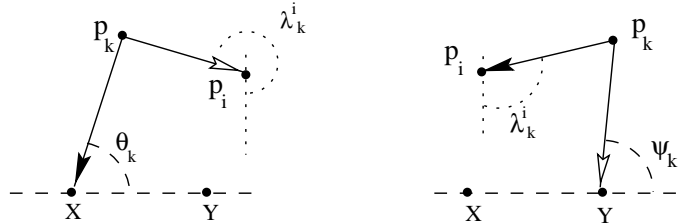 , height=3cm , clip= }}
\caption{New notation for angles}
\label{4p16}
\end{figure}

\no Note also that we revert to using $X$ and $Y$ for $0$ and $1$ on $\br$ so as to avoid confusing the $1$ with the vertex $p_1$. For our specific example we now have the picture
\begin{figure}[H]
\centering{\epsfig{file=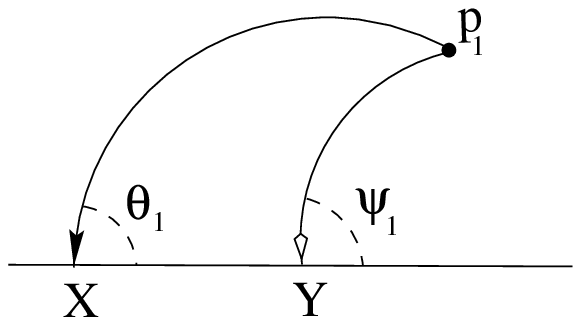 , height=3cm , clip= }}
\caption{New notation for wedge in $\ech$}
\label{4wp3}
\end{figure}

\no and in our current notation
\begin{equation}
\label{4en3}
w_1 = \int_{\ech} d\theta_1 \wedge d\psi_1 \; .
\end{equation}

\no There is some apparent loss of information with this naming procedure. It does not indicate which edge is the first one (this was earlier encoded in the $e_1^1$) and which is the second. But this is not really a problem since in the figures we are able to use the black leg,{\epsfig{file=leg1.eps , height=6mm , clip= }}, to indicate the first edge, and in the integrals it is encoded by the ordering in the form $d\theta_1 \wedge d\psi_1$.

\vv

\no Before we set up our integrals we would also like to address the issue of orientation. There are quite a few inputs which go into deciding the sign of our weight integrals. The standard orientation of \ech, the ordering of the wedges' feet, and the assigning of limits for our iterated integrals ($\int_a^b$ versus $\int_b^a$). All of these choices are {\em alternating}. \ie a single flip of any of these changes the sign of our weight and two flips get us back. So to fix things we do the following. We set up an unambiguous scheme that respects the above determinations and gives the right weight for some standard graph. That standardization is given by requiring that  $w_1 = \frac{1}{2}$. This follows because $G_1$ consists of the two graphs {\epsfig{file=cbga1.eps , height=8mm , clip= }} and {\epsfig{file=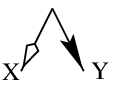 , height=8mm , clip= }} (Notice the difference!), where {\epsfig{file=cbga1.eps , height=8mm , clip= }} represents  $\frac{1}{2}\{\ ,\ \}$ .
Since the two graphs are obtained from each other by flipping edges ${w_I}$({\epsfig{file=cbga1.eps , height=6mm , clip= }})
= $-{w_I}$({\epsfig{file=cbga1m.eps , height=6mm , clip= }}).

\no Now the $\h$ term of $f \star g$ has to be $\frac{1}{2}\{\ f,g \}$. This implies

\begin{equation}
\label{4en2}
\frac{1}{2}\{\ f,g \} = w_1(\frac{1}{2}\{\ f,g \}) -w_1(\frac{1}{2}\{\ g,f \}) \; .
\end{equation}

\no But then $w_1$ must be $\frac{1}{2}$.
\vv
\no We now set up the iterated integral for $w_1 = \int_{\ech} d\theta_1 \wedge d\psi_1 $. 

\no Here is our scheme for assigning an iterated integral for $w_1$. We  use the two polar coordinates ($\theta$ and $\psi$) to represent $p_1$ ranging over \ech, one of them replacing the radial coordinate. Given the ordering of the 2-form as in Equation \eqref{4en3} we will take the second 1-form ($d\psi_1$ here) as the polar variable and it will give the outer integral. We use $\theta_1$ to represent the radial coordinate for a fixed $\psi_1$, and this gives the inner integral. Hopefully the following picture will make this clearer.

\vv

\begin{figure}[H]
\centering{\epsfig{file=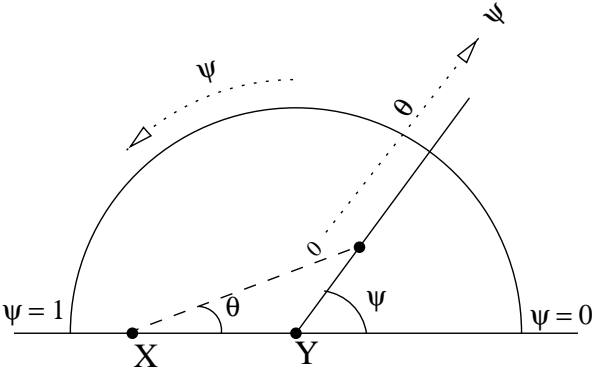 , height=5cm , clip= }}
\caption{Determining limits for $\int_{\ech} d\theta_1 \wedge d\psi_1$}
\label{4p19}
\end{figure}

\no  The inner limits are determined as follows. As $r$, the radial parameter, goes from $0$ to $\infty$ along a fixed $\psi_1$ ray, $\theta_1$ goes from $0$ to $\psi_1$. Written as an iterated integral
\begin{equation}
\label{4enn3}
\int_{\ech} d\theta_1 \wedge d\psi_1 = \int_0^1 ( \int_0^{\psi_1} d\theta_1 ) d\psi_1 \; .
\end{equation}

\no The actual integration is straightforward.

\begin{equation}
\label{4en4}
\int_0^1 ( \int_0^{\psi_1} d\theta_1 ) d\psi_1 \; = \; \int_0^1 \psi_1 d\psi_1 \; = \; [\frac{{\psi_1}^2}{2}]_0^1 \; = \; \frac{1}{2} \; .
\end{equation}

\no To reinforce our method we compute  ${w_I}$({\epsfig{file=cbga1m.eps , height=8mm , clip= }}) directly. We have,

\begin{equation}
\label{4en5}
\int_{\ech} d\psi_1 \wedge d\theta_1 = \int_0^1 ( \int_1^{\theta_1} d\psi_1 ) d\theta_1 \; .
\end{equation}

\no The limits are obtained from the following figure.
\begin{figure}[H]
\centering{\epsfig{file=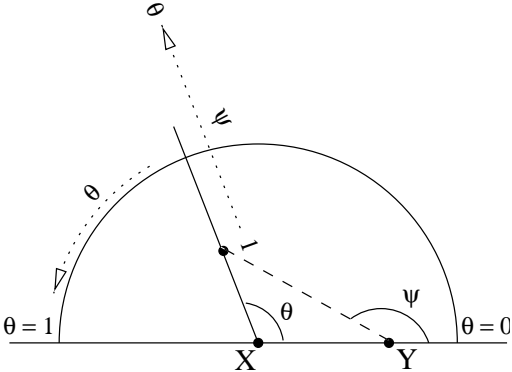 , height=5cm , clip= }}
\caption{Determining limits for $\int_{\ech} d\psi_1 \wedge d\theta_1$}
\label{4pn19}
\end{figure}

\vv

\no Once again the integration is straightforward and we skip the details. The value, as expected, is $-\frac{1}{2}$. Note that this also establishes the compatibility of our iterated integral set-up with the sign change associated with flipping the  2-form $ d\theta_1 \wedge d\psi_1$. We now explain the iterated integral setting for arbitrary W-computable graphs. 

\no Say we were interested in finding the weight for 
\begin{figure}[H]
\centering{\epsfig{file=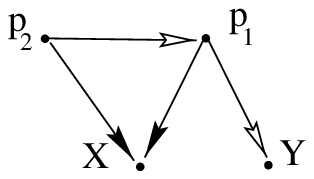 , height=2cm , clip= }} 
\caption{The simplest non-trivial W-computable graph}
\label{simp}
\end{figure}

\no The iterated integral  for the above graph's weight is a (double) integral over two copies of $\ech$. For a fixed $p_1$ we let $p_2$ range over $\ech$. This gives the inner integral and will lead to a function of $p_1$. We then integrate this function as $p_1$ ranges over $\ech$. The setting for more general W-computable graphs is no more complicated than this. One iteratively integrates out the contribution of the concatenating  grounded wedges, \epsfig{file=4cpic3.eps , height=6mm , clip= }, that generate a given graph, the outermost vertex corresponding to the innermost integral. To show that this is a manageable process we study the iterated integral set up for the  following $\frac{1}{2}{\rm ad}_X$ type  subgraph.
\begin{figure}[H]
\centering{\epsfig{file=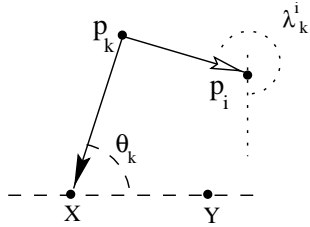 , height=3cm , clip= }}
\caption{Setting the integral for  ${\rm ad}_X$ }
\label{4p20}
\end{figure}

\no We introduce another bit of notation at this point. We denote the integral where $p_k$ ranges over \ech by $\int_{p_k}$ (instead of $\int_{\ech}$).
We claim that the  weight integral for the above graph is given by 
\begin{equation}
\label{4en6}
\int_{p_k} d\theta_k \wedge d\lak{i} = \int_0^1 ( \int_{\theta_i}^{\lak{i}} d\theta_k ) d\lak{i} \; .
\end{equation}

\no To see this we simply mimic the procedure we used to set up Equation \eqref{4enn3}. We let the angle $\lak{i}$ subtended by $p_k$ at $p_i$ vary from $0$ to $1$ for the outer integral, and for the inner one we determine the variation of $\theta_k$ as the point $p_k$ moves from the fixed $p_i$ to $\infty$ along a fixed  $\lak{i}$ level curve. Here is the corresponding figure.

\begin{figure}[H]
\centering{\epsfig{file=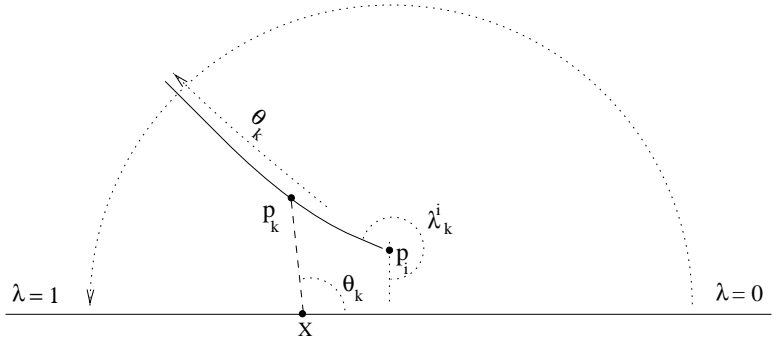 , height=4.5cm , clip= }}
\caption{Determining limits for $\int_{p_k} d\theta_k \wedge d\lak{i}$ }
\label{4p21}
\end{figure}

\no Thus far we have
\begin{equation}
\label{4enn6}
\int_{p_k} d\theta_k \wedge d\lak{i} = \int_0^1 ( \int_{?}^{??} d\theta_k ) d\lak{i} \; .
\end{equation}

\no Observe that for the inner integral's lower limit we start with $p_k$ at the same spot as $p_i$. Thus $\theta_k$ starts at $\theta_i$. Secondly, as $p_k$ travels to $\infty$, the angle it subtends at $X$, namely $\theta_k$, approaches the angle it subtends at $p_i$, namely $\lak{i}$. But since $p_k$ is travelling on a level curve for this latter angle $\lak{i}$, the upper limit for $\theta_k$ is $\lak{i}$. We now compute the integral in equation \eqref{4en6}.
\begin{align}
\label{4en7}
\int_{\ech} d\theta_k \wedge d\lak{i} &= \int_0^1 ( \int_{\theta_i}^{\lak{i}} d\theta_k ) d\lak{i} \notag  \\
&=  \int_0^1 ({\lak{i}} - \theta_i ) d\lak{i} \notag  \\
&= \biggl[ \frac{{\lak{i}}^2}{2} - \theta_i \lak{i}{\biggr]}_0^1 \notag  \\
&= \frac{1}{2} - \theta_i \; .
\end{align}

\no First of all note that Equation \ref{4en7} affords a recomputation of $w_1$when we set $\theta_i = 0$. This corresponds to the point $p_i$ in Figure \ref{4p20} being placed at $Y$, whose $\theta$ coordinate is $0$.
\begin{figure}[H]
\centering{\epsfig{file=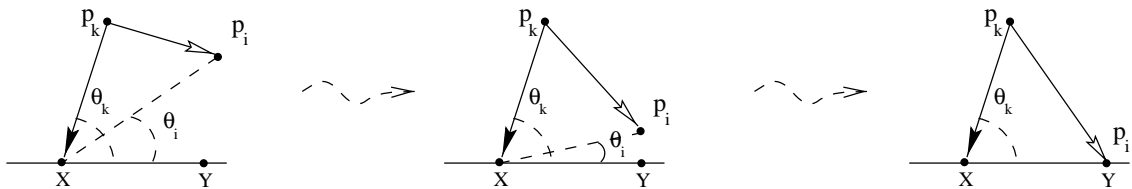 , height=25mm , clip= }}
\caption{Recomputing $w_1$}
\label{4pp20}
\end{figure}

\no Secondly, if the 1-forms in the above integral had been flipped, \ie, if we wished to compute $\int_{p_k} d\lak{i} \wedge  d\theta_k $, we could set up the corresponding integrals with the $\theta_k$ as the outer integral variable. But that turns out to be a lot more complicated (though still computable). The reason underlying this is that while $\theta$ (and $\psi$) are well-behaved along the level curves for $\lak{i}$, the converse is not true. Our branch for $\lak{i}$ has a jump discontinuity along a fixed vertical segment. Taking $\lak{i}$ as the inner integral variable would necessitate the subdivision of the integral. Earlier we had mentioned that we would be invoking basic calculus techniques in the computation of the weight integrals. This choice of a preferred order of integration is another instance of that. The moral of the above discussion is that we always integrate over the {\em ground} angles $\theta$ (and $\psi$) first.

\vv

\no Returning to Equation \eqref{4en7}, it can be interpreted verbally as follows. If the wedge with vertex at $p_k$ has one foot landing on $p_i$ and the other on $X$ (or $Y$), then integrating over the $p_k$ yields a linear function of
$\theta_i$ (or $\psi_i$). Moreover we would like to observe that the function $ \frac{1}{2} - \theta_i$ which  we obtained in  Equation \eqref{4en7} is $\Bb{1}(-\theta_i)$, where $\Bb{1}(x)$ is the first modified $\Bn{}$ernoulli polynomial that we defined in Definition \ref{bpold2}  at the end of Section \ref{chap2}. There are two implications of the above analysis that we would like to discuss.

\vv

\no For one, as we had discussed in the paragraph preceding Figure \ref{4p20}, it will be straightforward to compute weights for W-computable graphs. We first explain this for graphs of the form
\begin{figure}[H]
\centering{\epsfig{file=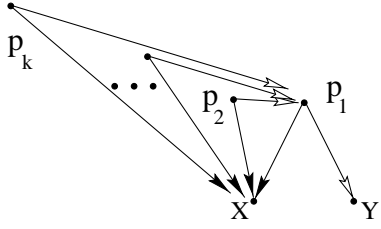 , height=3cm , clip= }}
\caption{Readily computable graphs - I}
\label{4pn1}
\end{figure}

\no Here each vertex $p_i, \; i>1$ has one foot on $X$ and the other on $p_1$.
Each of the $p_i, \; i>1$ can then be integrated independently, and one has to finally integrate $(\frac{1}{2} - \theta_1)^{k-1}$ as $p_1$ ranges over \ech.  For graphs of the form
\begin{figure}[H]
\centering{\epsfig{file=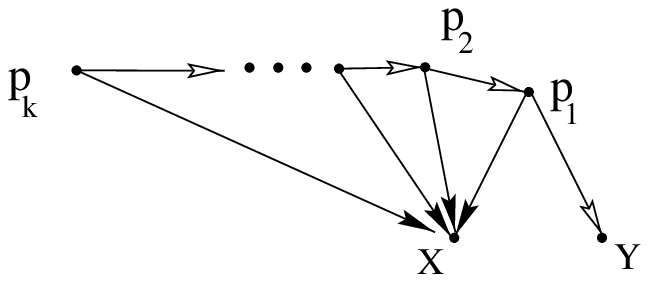 , height=3cm , clip= }}
\caption{Readily computable graphs - II}
\label{4pn2}
\end{figure}
 
\no integrating over the $p_k$ will yield a linear function of $\theta_{k-1}$.
Next, integrating over $p_{k-1}$ will  give a quadratic function of $\theta_{k-2}$. Proceeding in this fashion, finally we will be integrating a polynomial of degree $k-1$ in $\theta_1$ as $p_1$ varies over $\ech$ - this yielding the constant weight. 

\no The second implication of Equation \eqref{4en7} is that the weights that arise are intimately linked with the Bernoulli numbers. To exhibit this we return to computing the weight $w_I$ of {\epsfig{file=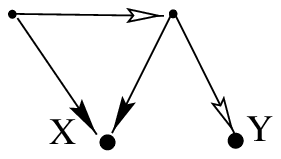 , height=6mm , clip= }}. Let us set $w_I({\epsfig{file=42wed.eps , height=6mm , clip= }}) = w_2$. As an iterated integral this is given by 
\begin{equation}
\label{4en9}
w_2 =   \int_0^1 \int_0^{\psi_1} \; (\int_0^1 \int_{\theta_1}^{\la{2}{1}} d\theta_2  d\la{2}{1}) \; d\theta_1 d\psi_1 \; .
\end{equation}

\no The inner integral giving the integral over $p_2$, the outer one over $p_1$. As in Equation \eqref{4en7} integrating out the $p_2$ yields $ (\frac{1}{2} - \theta_1)$. Integrating this as $p_1$ varies over \ech yields
\begin{align}
\label{4en10}
w_2 &=  \int_0^1 \int_0^{\psi_1}  (\frac{1}{2} - \theta_1) d\theta_1 d\psi_1 \notag \\
&=  \int_0^1 [\frac{\theta_1}{2} - \frac{{\theta_1}^2}{2}]_0^{\psi_1} d\psi_1  \notag \\
&=  \int_0^1 (\frac{\psi_1}{2} - \frac{{\psi_1}^2}{2}) d\psi_1  \notag \\
&= \biggl[\frac{{\psi_1}^2}{4} - \frac{{\psi_1}^3}{6}{\biggr]}_0^1 \notag \\
&= \frac{1}{12} 
\end{align}

\no which is $\frac{\Bn{2}}{2!}$, just as $w_1$ was $\frac{\Bn{1}}{1!} \; = \; \frac{1}{2}$. Here $\Bn{k}$ are the modified Bernoulli numbers that arose in Section \ref{chap2} (\eg refer to Equations \eqref{chp2s8} through \eqref{chp2s10}). We are now in a position to state the main result of this subsection.
\begin{prop}
\label{4bprop}
$${\mbox Let} \; \hat{\Gamma}_k \; \dfequal \epsfig{file=4pn2.eps , height=16mm , clip= } \; .\quad {\mbox Then,}
\;  w_k \; \dfequal \; w_I(\hat{\Gamma}_k ) \; = \; \frac{\Bn{k}}{k!} \; . \notag$$
\end{prop}

\no Proposition \ref{4bprop} will follow from a more general proposition that we shall establish below. The proof hinges on the fact that the iterated integrals for $w_k$ yield a modification of the Bernoulli polynomials, and this is proved by induction. Our first step is to further specialize our notation in order to ease the mechanics of the proof.

\no Given the form of the integrals at hand, with the first foot of vertex $p_k$ landing on $X$ and the other on $p_{k-1}$, we specialize the definition of $\int_{p_k}$ in Equation \eqref{4en6} to 
\begin{equation}
\label{4end}
\int_{p_k} \dfequal \int_0^1 ( \int_{\theta_{k-1}}^{\lambda_k^{k-1}} d\theta_k ) d\lambda_k^{k-1} \; .
\end{equation}
\no $\int_{p_k}$ can be interpreted to be the polynomial  in $\theta_{k-1}$ that is obtained when one carries out the integration on the right hand side of Equation \eqref{4end}. It follows then that $w_n$ is given by 
\begin{equation}
\label{4endy}
w_n = \int_{\ech} (\int_{p_{2}} \ldots \int_{p_n} )\quad d\theta_1 \wedge d\psi_1 \; , 
\end{equation}
\no the iterated integrals having the same general form up until the integration over $p_1$ (all earlier wedges are of the form \epsfig{file=4cpic3.eps , height=6mm , clip= }). There seems to be a slight hitch in incorporating the final $p_1$ integral in the form given by Equation \eqref{4end}, as in  this case the wedge is a totally grounded one and one gets a constant output for the integration, not a function. We shall see that this obstruction is easily surmountable. For this we would like to recall the discussion centered around Figure \ref{4pp20}. The algebraic version of that discussion is the following. One can think of $\int_{\ech}$, with $p_1$ ranging over $\ech$, as being given, not by 
\begin{equation}
\label{4end3}
\int_{\ech} d\theta_1 \wedge d\psi_1 = \int_0^1 ( \int_0^{\psi_1} d\theta_1 ) d\psi_1 
\end{equation}
\no as in Equation \eqref{4enn3}, but rather by
\begin{equation}
\label{4end4}
\int_{p_1} d\theta_1 \wedge d\psi_1 = \int_0^1 ( \int_{\theta_0}^{\psi_1} d\theta_1 ) d\psi_1 \; , 
\end{equation}

\no where $\theta_0$ is just a dummy variable which we then set equal to  $0$. Notice that Equation \eqref{4end4} does represent the integral over $p_1$  in the form of Equation \eqref{4end} since in this case we have 
$\lambda_k^{k-1} = \la{1}{Y} = \psi_1$. 
\begin{dfn}
\label{4definn}
We define the polynomial $\Pol{n}{\theta_0}$ by,
$$\Pol{n}{\theta_0} \dfequal \int_{p_1}\quad\int_{p_2}\;\ldots\;\int_{p_n} \notag$$
where all the $\int_{p_k}, \; 1 \leq k \leq n$, are given by equation \eqref{4end}.
\end{dfn}

\no The conclusion of the above formalism is that 
\begin{equation}
\label{4end5}
w_n \; = \; (\int_{p_1}\;\int_{p_2}\;\ldots\;\int_{p_n}) \bigg|_{\theta_0 = 0} \; = \; \Pol{n}{0}
\end{equation}
\no \ie, $w_n$ is nothing but the polynomial $\Pol{n}{\theta_0}$ evaluated at 
$\theta_0 = 0$. For notational ease we will denote the final dummy variable $\theta_0$ by $\theta$.
\begin{prop}
\label{4bprop2}
$\Pol{n}{\theta}
= \frac{1}{n!}\Bb{n}(-\theta)
\; .$
\end{prop}
\no {\bf Reminder!} Script letters \eg $\Pol{n}{\theta}$ denote a polynomial in $\theta$, not a scalar times $\theta$. Recall that 
\begin{equation}
\label{4pold2} \Bb{n}(x) = \sum_{k=0}^n \binom{n}{k} \Bn{n-k} x^{k} 
\end{equation}
\no is the modified $\Bn{}$ernoulli polynomial that we defined in Section \ref{chap2} (Definition \ref{bpold2}). In passing we would like to mention the following (readily verifiable) relation between the modified $\Bn{}$ernoulli and Bernoulli polynomials. 
\begin{equation}
\label{bpobs}
 \Bb{n}(-\theta) =  (-1)^{n} \Be{n}(\theta) 
\end{equation}

\vv

\no Now $\Bb{n}(0) = \Bn{n}$, and Proposition \ref{4bprop} follows from Equation \eqref{4end5} and Proposition \ref{4bprop2}. We now prove Proposition \ref{4bprop2} by induction.

\begin{proof}
Suppose we do have $\Pol{n-1}{\theta} = \frac{1}{(n-1)!}\Bb{n-1}(-\theta)$. We would like to show that it also holds for  $\Pol{n}{\theta}$. Now, using the definition of the  $\Pol{n}{\theta}$ (Definition \ref{4definn}),

\begin{align}
\Pol{n}{\theta} &= \int_{p_1}\Pol{n-1}{\theta_1} \notag \\
&= \int_0^1 \bigl( \int_{\theta}^{\psi_1}\Pol{n-1}{\theta_1} d\theta_1 \bigr) d\psi_1 \notag  \\
& \quad \mbox{where we used  Equation \eqref{4end4} to obtain the above line.} \notag \\
& \quad \mbox{We now invoke the induction hypothesis and replace $\Pol{n-1}{\theta_1}$. } \notag \\
\Pol{n}{\theta} &= \int_0^1 \bigl( \int_{\theta}^{\psi_1}\frac{1}{(n-1)!}\Bb{n-1}{(-\theta_1)} d\theta_1 \bigr) d\psi_1 \notag \\
& \quad \mbox{Next we expand the $\Bb{n-1}{(-\theta_1)}$ using Equation \eqref{4pold2}} \notag \\
\Pol{n}{\theta} &=\int_0^1 \bigl( \int_{\theta}^{\psi_1}\frac{1}{(n-1)!}\sum_{k=0}^{n-1} \binom{n-1}{k} \Bn{n-1-k}(-\theta_1)^k d\theta_1 \bigr) d\psi_1  \notag 
\end{align}
\no At this stage we interchange the order of the summation and the integration  and integrate out the $d\theta_1$ and the  $d\psi_1$.
\begin{align}
\Pol{n}{\theta} &=\frac{1}{(n-1)!}\sum_{k=0}^{n-1} (-1)^k \binom{n-1}{k} \Bn{n-1-k}\int_0^1 ( \int_{\theta}^{\psi_1}{\theta_1}^k d\theta_1 ) d\psi_1 \notag \\
&=\sum_{k=0}^{n-1}(-1)^k \frac{1}{k!(n-1-k)!} \Bn{n-1-k}\int_0^1 \biggl[ \frac{{\theta_1}^{k+1}}{k+1} \biggr]_{\theta_1 = \theta}^{\theta_1 = \psi_1} \; d\psi_1 \notag \\
&= \sum_{k=0}^{n-1}(-1)^k \frac{1}{k!(n-1-k)!} \Bn{n-1-k}\int_0^1( \frac{{\psi_1}^{k+1}}{k+1} - \frac{{\theta}^{k+1}}{k+1}) d\psi_1 \notag \\
&= \sum_{k=0}^{n-1}(-1)^k \frac{1}{k!(n-1-k)!} \Bn{n-1-k} \biggl[ \frac{{\psi_1}^{k+2}}{(k+1)(k+2)} - \frac{{\theta}^{k+1}\psi_1}{k+1}\biggr]_{\psi_1 = 0}^{\psi_1 = 1}   \notag \\
&= \sum_{k=0}^{n-1}(-1)^k \frac{1}{k!(n-1-k)!} \Bn{n-1-k}[\frac{1}{(k+1)(k+2)} - \frac{{\theta}^{k+1}}{k+1}] \notag
\end{align}
\no We now gather terms so as to be in a position to prove the induction hypothesis.
\begin{align}
\Pol{n}{\theta} &= \sum_{k=0}^{n-1}(-1)^k \frac{1}{(k+2)!(n-1-k)!} \Bn{n-1-k} \notag \\
&\quad \quad  + \quad
\sum_{k=0}^{n-1}(-1)^{k+1} \frac{1}{(k+1)!(n-1-k)!} \Bn{n-1-k}{\theta}^{k+1}
\notag \\
&= \sum_{k+1=1}^{n}(-1)^k \frac{1}{((k+1)+1)!(n-(k+1))!} \Bn{n-(k+1)} \notag \\& \quad \quad + \quad
\sum_{k+1=1}^{n}(-1)^{k+1} \frac{1}{(k+1)!(n-(k+1))!} \Bn{n-(k+1)}{\theta}^{k+1} 
\notag \\
&= \sum_{l=1}^{n}(-1)^k \frac{1}{(l+1)!(n-l)!} \Bn{n-l} \quad + \quad
\sum_{l=1}^{n}(-1)^{l} \frac{1}{l!(n-l)!} \Bn{n-l}{\theta}^{l} \notag
\end{align}
\no where in the last line we set $l = k+1$. Finally we multiply and divide the right hand side of the last line by $n!$ and collect the appropriate terms to obtain
\begin{equation}
\label{4proof}
\Pol{n}{\theta} = \frac{1}{n!}\sum_{l=1}^{n}\frac{(-1)^k}{l+1} \binom{n}{l} \Bn{n-l} \; + \;
\frac{1}{n!}\sum_{l=1}^{n} \binom{n}{l} \Bn{n-l}{(-\theta)}^l
\end{equation}

\no Comparing Equation \eqref{4proof} and Equation \eqref{4pold2}, we would be done if $$ \frac{1}{n!}\sum_{l=1}^{n}\frac{(-1)^{l-1}}{l+1} \binom{n}{l} \Bn{n-l}
= \frac{1}{n!}\Bn{n} \quad .$$

\no Or, pulling all the terms to one side, if
$$\frac{1}{n!}\sum_{l=0}^{n}\frac{(-1)^{l-1}}{l+1} \binom{n}{l} \Bn{n-l}
= 0 \quad .$$
\no Or, setting $n-l = k$, if
$$\frac{1}{n!}\sum_{k=0}^{n}\frac{(-1)^{n-k-1}}{n-k+1} \binom{n}{k} \Bn{k}
= 0 \quad .$$
\no Or, pulling a constant ($\frac{(-1)^{n-2k+1}}{n!}$) out, if
$$\sum_{k=0}^{n}\frac{(-1)^{k}}{n-k+1} \binom{n}{k} \Bn{k}
= 0 \quad .$$
\no But this is exactly the condition satisfied by the $\Bn{k}$ in Equation \eqref{chp2s9}. And so we are done.
\end{proof}

\no One could try presenting a slicker proof, using other identities available in the literature, but we felt it was best to relate the equations here to those that arose earlier in our work. We would like to point out though that proofs of the Euler summation formula (see \cite{gkp} for example) bear a resemblance to the above proof.

\no In the next section we will use  this determination of the $w_n$ to prove  that bi-differential operators that are common to the CBH and Kontsevich quantizations have equal weights. In fact the  above valuation for the $w_n$ gives a direct equivalence for operators of the form $\dell{i}\delr{{j_1}}\delr{{j_2}} \ldots \delr{{j_n}}$, and as in Section \ref{chap2} we use that to establish the rest of the equivalence.

\subsection{Loop vs Non-loop}
\label{chp4f}
\no In this final section we would like to make some comments about graphs with loops and those without. First we would like to state and prove the following theorem.

\begin{thm} 
\label{timp}
A non-loop graph $\Gamma_n \in G_n$ is formed by the concatenation of a single wedge and a non-loop graph $\Gamma_{n-1} \in G_{n-1}$, where by concatenation we mean that the feet of the wedge are allowed to land on any of the vertices of $\Gamma_{n-1}$, not simply the ground vertices.
\end{thm}
\no 
\begin{proof}\footnote{We would like to thank Professor Alan Weinstein for the idea underlying this proof.} Let $\Gamma$ be a non-loop graph made up of $n$ wedges. Since it has no loops, following the edges out of any vertex always gives a chain (\ie a set homeomorphic to an interval). We define the length of a chain to be the number of edges in the chain, and the height of a vertex to be the length of the longest chain emanating from it. The fact that $\Gamma$ is finite and has no loops implies that these notions are well-defined and not vacuous. The finiteness of $\Gamma$ also implies that there is at least one vertex, say $k$, of maximal height. Consequently this vertex must be free (\ie no foot lands on it) because otherwise there would be a vertex with greater height. But then erasing the $k$th wedge gives a non-loop graph $\Gamma_{n-1}$ such that $\Gamma$ is obtained by the concatenation of the $k$th wedge with $\Gamma_{n-1}$.
\end{proof}

\no There is a `physical' metaphor that one can use to try and explain the difference between the arising of loop and non-loop graphs (I apologize in advance for the over-simplified nature of the following imagery). Let us think of $X$ and $Y$ as two iron tabs, and the wedges as being flexible horseshoe magnets. Say one wants to list all the possible agglomeration of wedges as one throws the wedges down to the $X$, $Y$ tabs. One way is to throw them one at a time, not throwing a wedge until the earlier one has settled. This procedure would lead to non-loop graphs only, for it mimics the concatenation alluded to in the above theorem.

\no On the other hand one could throw more than one wedge simultaneously. These could then interact and pair up (or triple, or $\ldots$) while in the air, forming loops before they landed. There are definitely more bi-differential operators arising via the second method, but it is not clear if the extra freedom is required for associativity. 

\no As we have seen via the CBH-quantization, the non-loop graphs were adequate for providing a \stp in the linear case. This may lead one to conjecture that even in the Kontsevich quantization the weights for the loop graphs will be zero. But Kontsevich's own computations \cite{ko2} yield non-zero weights for the loop graphs. Also, work by Penkava and Vanhaecke \cite{pv} on generalizing the universal enveloping algebra approach  to polynomial  Poisson structures suggests that the loop graphs have non-zero weight and that they could  play a non-trivial role  for quantizing cubic and higher degree Poisson structures. They  show that the natural extension of the enveloping algebra approach does give a deformation quantization up to order ${\h}^3$, but there are obstructions to extending it to order ${\h}^4$ (an explicit example of a cubic Poisson structure is given). They compute a correction term which precisely cancels the obstruction and the graphical description of the corresponding bi-differential operator needs loop graphs. 

\vv

\no Though we cannot figure out the weights for the loop graphs, using the machinery developed here and in Section \ref{chap2} we can show that all the coefficients (by default for non-loop graphs) in the CBH-quantization agree with the corresponding coefficients in Kontsevich's quantization. Thus we could set the weights of all the loop graphs to zero and still have an associative product in the linear situation. In the next section we pull together the various ingredients needed to prove the preceding contention and examine some of the consequences.

\section{Summing Up}
\label{chapfin}
\no In this section we study Kontsevich's quantization for linear Poisson structures and relate it to the CBH-quantization. We establish the following theorem.
\begin{thm}
\label{finthm}
Apart from bi-differential operators that arise from loop graphs, the Kontsevich and CBH quantizations for $\fgs$ are identical.
\end{thm}
\no Theorem \ref{finthm} has the following symbolic interpretation. If the bi-differential operator for the CBH-quantization is given by 
\begin{equation}
\label{5new1}
\exp(C), \quad {\mbox{where}} \quad C \; = \; \frac{\h}{2}({\epsfig{file=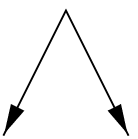 , height=8mm , clip= }}) \; + \; \frac{\h^2}{12}({\epsfig{file=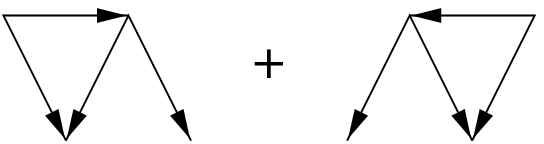 , height=8mm , clip= }}) \; + \; \ldots ,
\end{equation}
then the bi-differential operator for the Kontsevich quantization for $\fgs$ is given by $\exp(C + L)$, where $L$ is a bi-differential operator encoded by prime {\bf loop} graphs.

\no The proof is based on the fact that the two quantizations match exactly for duals of arbitrary nilpotent Lie algebras. We outline our approach below. Except for the last subsection, by quantization we will mean a quantization for $\fgs$.
\newcounter{u}
\begin{list}
{(\arabic{u})}{\usecounter{u}}
\item We show that the bi-differential operators in Kontsevich's quantization are those that arise in the CBH quantization (from Sym-admissible graphs) and those that arise from loop graphs (absent in the CBH quantization).
\item The bi-differential operators corresponding to loop graphs generalize the Killing form and vanish on all nilpotent Lie algebras.
\item As a consequence of (2) above and our computations in Section \ref{4wc},
we are able to show that for nilpotent Lie algebras $\{(X)^n \star Y : \; X,Y \in \fg, \; n \in \bn \}$ are identical for the CBH and Kontsevich quantizations. Recall that by $(X)^n$ we mean the monomial given by the $n$th power of $X$, where we think of $X$ as a coordinate function on $\rd \cong \fgs$. As in Section \ref{c2s1},  we are then able to show that the Kontsevich and CBH quantizations match on nilpotent Lie algebras since they both satisfy Property (P1) introduced in Section \ref{c2s1}. 
\item Given the universality of the coefficients in the CBH and Kontsevich schemes, and the fact that by taking larger and larger dimensional nilpotent Lie algebras we can account for any non-zero coefficient in the CBH formula, we finally obtain a proof of Theorem \ref{finthm}.
\end{list}

\no The various subsections of this section address the above points. Finally, to end this paper, we indicate some avenues for further research.

\subsection{Admissible Graphs in the Linear Setting}
\label{finsec1}
\no The key restriction that makes Kontsevich's \stp amenable to analysis 
in the linear setting is that we only have to consider graphs where no more than one edge can land on an aerial vertex. This is because the Poisson structure being linear, all graphs containing the subgraph {\epsfig{file=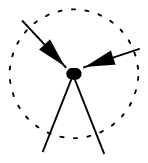 , height=6mm , clip= }} contribute nothing. Consequently aerial vertices of interest are either of the form {\epsfig{file=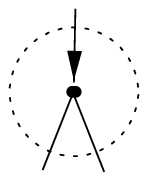 , height=8mm , clip= }},  or they are free vertices (roots) {\epsfig{file=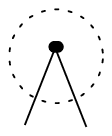 , height=8mm , clip= }}. 
\begin{thm}
\label{thm511}
If a graph contains no loops, and its vertices are of type {\epsfig{file=5p2.eps , height=8mm , clip= }} or {\epsfig{file=5p3.eps , height=8mm , clip= }} then it corresponds to a (possibly multiple rooted) binary tree.
\end{thm}
\no Note that when we use the word tree we ignore the coalescing of edges at the $X$ and $Y$ vertices. \ie, we think of graphs such as  

\no {\epsfig{file=cbg05.eps , height=16mm , clip= }} as being equivalent to {\epsfig{file=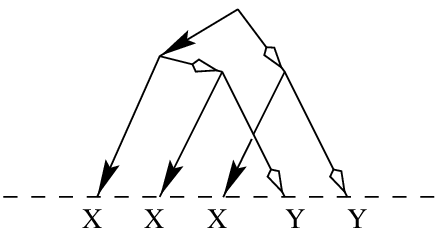 , height=16mm , clip= }}.

\begin{proof} (for Theorem \ref{thm511})
From Theorem \ref{timp} it follows that non-loop graphs $\Gamma_n$ are formed by concatenation of a wedge with another non-loop graph $\Gamma_{n-1}$. But now we only allow for concatenations that result in nodes of the form  {\epsfig{file=5p2.eps , height=6mm , clip= }} and  {\epsfig{file=5p3.eps , height=6mm , clip= }}. This means that a new wedge can land only on $X$, $Y$, or a {\bf free} vertex. Suppose now that at the $n$th stage this procedure has only led to binary trees. Then at the $n+1$st stage we still get binary trees, since all we have done is introduced a new root whose branches are from among the earlier trees, $X$, or $Y$. The fact that we start with a single wedge completes the proof.
\end{proof}

\no Now in Section \ref{freeg} (on page~\pageref{cbg3m}) we indicated that binary, rooted trees (where we allow for multiple roots) correspond exactly to the Sym-admissible graphs. Recall that these account for the bi-differential operators that arise in the CBH-quantization. Thus we can reach the following conclusion. 
\begin{prop}
\label{finprop1}
Bi-differential operators occurring in the Kontsevich quantization for $\fgs$ are of two types:
\newcounter{uv}
\begin{list}
{(\alph{uv})}{\usecounter{uv}}
\item Those corresponding to loop graphs;
\item Those corresponding to Sym-admissible graphs.
\end{list}
\end{prop}

\no We now move on to analyzing the loop graphs. In particular we show that the
corresponding bi-differential operators vanish for nilpotent Lie algebras.

\subsection{Loop Graph Contributions}
The simplest loop graph $\Gamma_l$ is the following one, which is constructible using two wedges.
\begin{figure}[H]
\centering{\epsfig{file=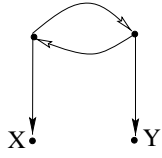 , height=2cm , clip= }}
\caption{The simplest loop graph -  $\Gamma_l$}
\label{5p5}
\end{figure}
\no Note that the definition of admissible graphs (Definition \ref{defff}) precludes the possibility of the loops of the type  {\epsfig{file=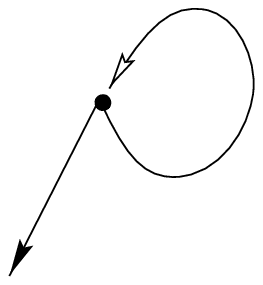 , height=8mm , clip= }}. Let us write out the bi-differential operator corresponding to  $\Gamma_l$. For this we first color the graph appropriately and then follow the procedure given in Equation \eqref{kodef3}.
\begin{figure}[H]
\centering{\epsfig{file=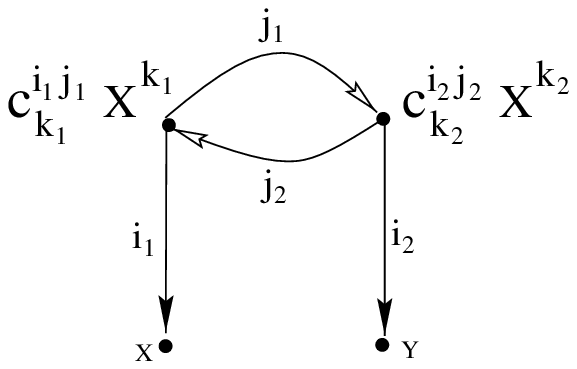 , height=3cm , clip= }}
\caption{ $\Gamma_l$ with coloring}
\label{5p6}
\end{figure}
\no We have omitted factors of $\frac{1}{2}$ in the interest of clarity. Carrying the factors is not a crucial concern since we will show that the corresponding operator vanishes in the case of interest.  The corresponding bi-differential operator is given by 
\begin{align}
\label{defcol}
B_{\Gamma_l} &=\del{{j_2}}(c_{k_1}^{{i_1}{j_1}} X^{k_1})\; \del{{j_1}}(c_{k_2}^{{i_2}{j_2}} X^{k_2})  \; \dell{{i_1}}\;\delr{{i_2}} & \notag \\
&=  c_{j_2}^{{i_1}{j_1}}  c_{j_1}^{{i_2}{j_2}} \; \dell{{i_1}} \; \delr{{i_2}} \quad & ({\because \del{i} X^k = \delta_i^k})
\end{align}
\no where we are using the Einstein summation convention. We wish to point out that applying $B_{\Gamma_l}$ to linear functions (Lie algebra elements) gives the Killing form. This follows because
\begin{align}
\label{killf}
B_{\Gamma_l}(X^i,X^j) &=   c_{j_2}^{{i_1}{j_1}}  c_{j_1}^{{i_2}{j_2}} \; \dell{{i_1}}(X^i) \; \delr{{i_2}}(X^j) \notag \\
&=   c_{j_2}^{{i}{j_1}}  c_{j_1}^{{j}{j_2}} \quad .
\end{align}
\no T he end result of Equation \eqref{killf} gives the matrix coefficients for the Killing form with respect to the basis $\{X^1,X^2,\ldots,X^d\}$. We also wish to note that the decoration of the feet of graph $\Gamma_l$, that we use to indicate the ordering of edges, was inconsequential for our argument. Any other ordering would have led at most to a sign variation\footnote{Even this sign change does not filter through to the final bi-differential operator in the \stp because the corresponding sign change in the weight cancels it out.}. 

\no Next we move on to more general loop graphs. We are interested in the coefficient generated by the subgraph that gives the loop. Here by coefficient we mean the function resulting from the contraction of vector fields  and the coefficients of the Poisson structure at the vertices. \eg, in the case of $\Gamma_l$ above, we obtained the function $ c_{j_2}^{{i_1}{j_1}}  c_{j_1}^{{i_2}{j_2}}$.

\no Consider a loop with $n$ vertices. Since we are only interested in the loop we omit drawing the other edges emanating from a vertex.
\begin{figure}[H]
\centering{\epsfig{file=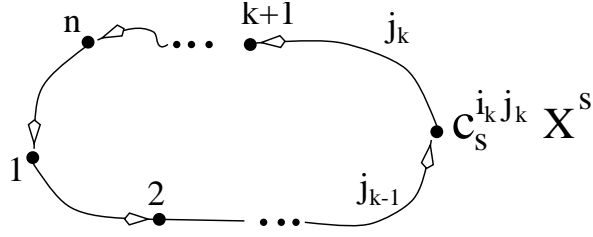 , height=3cm , clip= }}
\caption{A generic loop}
\label{5p8}
\end{figure}
\no In Figure \ref{5p8} we have numbered the vertices other than the $k$th one. At the $k$th vertex we have chosen to display the coloring of the vertex and the relevant incoming and outgoing edges. It follows that the resulting coefficient at vertex $k$ will be 
$c_{j_{k-1}}^{{i_k}{j_k}}$. The complete coefficient due to the entire loop will be the (cyclic) product 
\begin{equation}
\label{cycpro}
c_{j_n}^{{i_1}{j_1}}  c_{j_1}^{{i_2}{j_2}} \ldots
c_{j_{n-1}}^{{i_n}{j_n}} \quad .
\end{equation}
We now specialize to nilpotent Lie algebras, primarily motivated by the following theorem. A proof is available in Varadarajan's text on Lie groups \cite{vr}.
\begin{thm}
\label{thmstruc}
Let $\fg$ be a $d$-dimensional nilpotent Lie algebra over a field of characteristic $0$. Then there is a basis $\{X^1, X^2, \ldots, X^d\}$ for $\fg$ such that the structure constants $c_k^{ij}$ defined by $[X^i,X^j] = c_k^{ij} X^k$
satisfy
\begin{equation}
c_r^{ij} = 0 \quad r \geq \mbox{min}(i,j) \; . \notag
\end{equation}
\end{thm}
\begin{cor}
\label{corny}
For a nilpotent Lie algebra, any cyclic product of the form
$$c_{j_n}^{{i_1}{j_1}}  c_{j_1}^{{i_2}{j_2}} \ldots
c_{j_{n-1}}^{{i_n}{j_n}} \;(n>1)\; \mbox{equals \; {\bf ZERO}}. \notag $$
\end{cor}
\begin{proof}
The above product would be non-zero if and only if all the constituent $c_k^{ij}$ where non-zero. Given Theorem \ref{thmstruc} (at the very least) it would require
$$ j_n \; <  \;  j_1 \; < \;  j_2 \; < \ldots < \; j_{n-1} \; < \; j_n \;. \notag$$
which is impossible for $n > 1$.
\end{proof}
\no Note that the vanishing of the Killing form for nilpotent Lie algebras follows when $n=2$.
\no Using Proposition \ref{finprop1} and Corollary \ref{corny} we observe that the non-zero bi-differential operators in the Kontsevich quantization for the dual of a nilpotent Lie algebra are a subset of the ones obtained from Sym-admissible graphs. We now proceed to show that the Sym-admissible graph contribution
matches with the CBH-quantization.

\subsection{Equivalence of {\stp}s}
Recall from Section \ref{c2s1} that two {\stp}s on $\fgs$ are equal if they satisfy the following two properties:
 \newcounter{b2}
\begin{list}
{(P\arabic{b2})}{\usecounter{b2}}
\item For any two polynomials $p_{n}$ and $q_{m}$ on $\fgs$ (or 
equivalently $p_{n}, q_{m} \in S(\fg)$) of degree
$n$ and $m$ respectively, 
\begin{equation}
\label{stformy}
p_{n} \star q_{m} = p_{n} q_{m} + \mbox{terms of degree $< m+n$}. \notag
\end{equation}
\item The two products are equal for products of the form $(X)^n \star Y, \; n \in \bn$, where $X$ and $Y$ are arbitrary elements of $\fg$ and $(X)^n$ stands for the $n$th power of $X$ in $S(\fg)$.
\end{list}
\no We have already established property (P1) for the CBH-quantization. We now prove that it holds for the Kontsevich quantization in the linear setting. From the  discussion in Section \ref{finsec1} it follows that a graph will yield a non-zero bi-differential operator only if there is no more than one incident edge at each of the graph's aerial vertices. But there are two edges emanating from each such vertex. Hence for a graph with $n$ aerial vertices, at least $n$ edges have to descend earthward. Since these act on the functions whose \stp we are taking, all bi-differential operators in the Kontsevich product that arise from graphs with more than one aerial vertex, lower the  total degree. The only graph that is excluded from this collection is $\Gamma_0$, the graph that yields pointwise multiplication. This establishes Property (P1) for the Kontsevich product (on $\fgs$).
\no Now all that remains to be shown is that $(X)^n \star Y, \; n \in \bn$, is the same in the CBH and Kontsevich quantizations. Let us first list the relevant
bi-differential operators (ones that contribute to  $(X)^n \star Y$) in the CBH-quantization. These were given in Equation \eqref{bersym1} in Section \ref{cbhsec}. They are
\begin{equation}
\label{bersym11}
I + \sum_{k=1}^{\infty}\frac{\Bn{k}}{k!}(({\rm ad}_X)^k(Y)) (\dell{X})^k \delr{Y} .
\end{equation}

\no where the $I$ as usual represents pointwise multiplication. We now gather the graphs that will contribute to  $(X)^n \star Y$ in the Kontsevich product. Since we are restricting ourselves to the nilpotent Lie algebras, we need not concern ourselves with loop graphs. So we only have to find relevant graphs that can be built by concatenating wedges. Given that we are interested in  $(X)^n \star Y$ there is one more restriction in addition to the `one incident edge per vertex' rule - we need not consider any graphs where more than one edge lands on the linear function $Y$. This leads to the following simple picture. After the first wedge has been placed as a grounded wedge (\epsfig{file=cbga1.eps , height=6mm , clip= }), subsequent wedges have a single set of vertices to land on - the $X$ vertex and the lone free vertex available.
\begin{figure}[H]
\centering{\epsfig{file=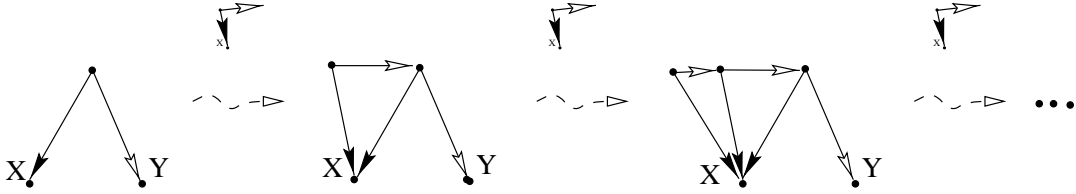 , height=2.5cm , clip= }}
\caption{Relevant graphs for  $(X)^n \star Y$ }
\label{5p10}
\end{figure}
\no But these are none other than the family of graphs $\{\hat{\Gamma}_m \}$\footnote{Actually the $\hat{\Gamma}_m$ do not account for all the relevant graphs but only determine the relevant topological type. We will address this issue in the following paragraphs.} whose weights $w_I(\hat{\Gamma}_m) = w_m$ were given in Proposition \ref{4bprop} in Section \ref{4wc}. It was shown there that $w_m = \frac{\Bn{m}}{m!}$. Notice the equality of this weight with the numerical coefficient in Equation \eqref{bersym11}. To show that this correspondence is actually an equality of bi-differential operators will be a matter of some straightforward book-keeping. Here are the factors that need to be balanced:

\newcounter{ac}
\begin{list}
{(\arabic{ac})}{\usecounter{ac}}
\item As mentioned in the earlier footnote, we need to account for all the graphs that can be obtained by re-labelling and reordering the vertices and edges of the $\hat{\Gamma}_m$. 
\item We should be using the weights $w_K$, not $w_I$, where $w_K(\hat{\Gamma}_m) = 
\frac{w_I(\hat{\Gamma}_m)}{m!}$.
\item The wedges that go into the $\hat{\Gamma}_m$ are colored thus {\epsfig{file=3p7b.eps , height=15mm , clip= }}.
\end{list}
\no Amazingly enough, all these factors balance out perfectly to yield the bi-differential operators that we have  in Equation \eqref{bersym11}. First we will account for the numerical contribution. Since $\hat{\Gamma}_m$ has $m$ aerial vertices that are to be labelled from $1$ through $m$, and there are two choices for the ordering of the edges emanating from each vertex, there are actually $m!2^m$ distinguishable graphs that have the same topological type as $\hat{\Gamma}_m$. They all lead to the same bi-differential operator because the labelling of the vertices is irrelevant to the process of assigning bi-differential operators, and the flipping of edges introduces a change of sign (due to the skew-symmetry of the $\alpha^{ij}$) which is exactly compensated by the corresponding change of sign in the graph's weight. Now the $m!$ will be balanced by the $m!$ that needs to be introduced in the denominator to transit from the $w_I$ to the $w_K$. The $2^m$ will be balanced by the $m$ $(\frac{1}{2})$'s that are imbedded in the $\frac{1}{2}\alpha^{ij}$ coloring of the wedges in $\hat{\Gamma}_m$.

\no The only thing that remains to be shown now is that when we specialize to taking just the $X$ derivatives (instead of arbitrary colorings) for the edges of $\hat{\Gamma}_m$ that land at the $X$ ground vertex (since we wish to compute $(X)^n \star Y$), then the coefficient formed by the contraction of the Poisson coefficients gives exactly the $({\rm ad}_X)^m(Y)$ term in Equation  \eqref{bersym11}. To show this we let $X$ be the first coordinate function $X^1$ and $Y$ be $X^2$. Then the colored $\hat{\Gamma}_m$ looks as in the following figure. (Note that we have already extracted the $\frac{1}{2}$ factor of the coloring.)
\begin{figure}[H]
\centering{\epsfig{file=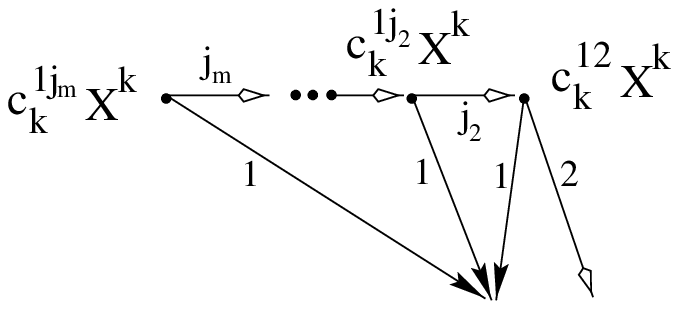 , height=3cm , clip= }}
\caption{$\hat{\Gamma}_m$ suitably colored}
\label{5pfin}
\end{figure}
\no The contraction of functions and vector fields along the horizontal line in Figure \ref{5pfin} precisely leads to   $({\rm ad}_X)^m(Y)$.
Hence on duals of nilpotent Lie algebras the CBH and Kontsevich quantizations are equal.

\no Now given any $n$, one can choose a nilpotent Lie algebra of rank $n$ such that Lie elements occurring in the CBH formula via $n$-fold bracketing are non-zero in this nilpotent Lie algebra. But then the universality of the CBH and Kontsevich quantizations gives the equality of the corresponding weights (for non-loop graphs) for every  $\fgs$, not just those that correspond to  nilpotent $\fg$. Since every term in the CBH formula can be accounted for in this manner, this establishes Theorem \ref{finthm}.

\subsection{Epilogue}
\no Apart from exhibiting the `equivalence' between the Kontsevich and CBH quantizations we hope to have convinced the reader that universal formulae for deformation quantization of $\rd$ seem to have the form $\exp$(deformation of the Poisson bi-vector field). Let us elaborate a bit on this point. 

\no In the earlier subsection we showed how the $m!$ in the order of the symmetry group for a prime graph cancelled out the $m!$ in the denominator of the weights $w_K$. Let us now consider the case of composite graphs. We start with the simplest possible one, $(\hat{\Gamma}_1)^2$.

\vv

\begin{figure}[H]
\centering{\epsfig{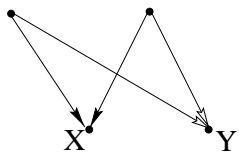}}
\caption{$(\hat{\Gamma}_1)^2$}
\label{expcbg3m}
\end{figure}

\no From our earlier computations it follows that $w_K ((\hat{\Gamma}_1)^2) = \frac{1}{2!} (\frac{1}{2})^2$. But this time the order of the symmetry group (the set of all labelled and oriented graphs with the same topological type) is only $2^2$. This is because the constituent wedges of $(\hat{\Gamma}_1)^2$ are indistinguishable and the re-labelling of vertices does not lead to a differently labelled graph. Thus when one collects all the labelled and ordered graphs of type \epsfig{file=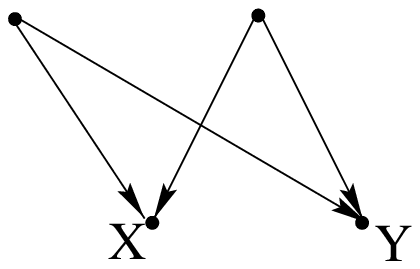 , height=8mm , clip= } one still has
\begin{equation}
\label{expw}
w(\epsfig{file=cbg3m5.eps , height=8mm , clip= })\; = \; \frac{1}{2!} w(\epsfig{file=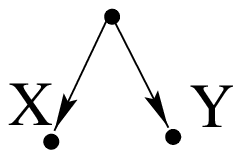 , height=8mm , clip= }) w(\epsfig{file=4p7a5.eps , height=8mm , clip= }) \; .
\end{equation}
\no Here by $w(\Gamma)$ we mean $\sum_{{\Gamma}^{\prime}} w_K({\Gamma}^{\prime})$ where the sum is over all labelled and ordered versions ${\Gamma}^{\prime}$ of $\Gamma$. But this is exactly as it should be, were composite graphs only allowed to occur via exponentiation. One can verify that one has the corresponding result for all composite graphs. This establishes the fact that Kontsevich's \stp is also of the form $\exp$(sum of weighted prime graphs), where the argument of the $\exp$ is a sum that can be viewed as a deformation of the Poisson structure. 

\vv

\no What we haven't done is characterize this deformation of the Poisson structure (the {\em logarithm} of the \stp). Below we list a series of questions that seem to arise from our work and still need to be addressed.
\newcounter{zac}
\begin{list}
{(\arabic{zac})}{\usecounter{zac}}
\item Can one show that all universal formulae for deformation quantization of $\rd$ are
of the type $\exp$(sum of weighted `prime' bi-differential operators)? If yes, or in the cases that are of this form, is there a way to characterize the {\em logarithm} of the \stp? In the CBH and Kontsevich quantizations  we have that the logarithm is symmetric (contributing graphs being symmetric about the vertical axis). Also, in the CBH case the Hausdorff series satisfies the associativity equation $H(X, H(Y,Z)) = H(H(X,Y),Z)$. Do similar characterizations exist for general quantizations? 
\item One can show in the CBH quantization  that if  one adds {\em central} terms to the `logarithm'  then one still gets a \stp. By central terms we mean bi-differential operators that correspond to central elements in  the universal enveloping  algebra of the given Lie algebra, \eg, the bi-differential operator generalizing the Killing form corresponds to the universal Casimir element.  Is there a similar principle for arbitrary deformation quantizations? 
\item Can one find a purely combinatorial approach to computing the weights for the tractable graphs? The fact that the integrations only involve polynomials which result in rational weights provides some hope for this. Though not presented here, there is evidence pointing to  such a principle in the case of the W-computable graphs.
\item How do Kontsevich's weights compare with other computations of the CBH coefficients, in particular the formulae of Dynkin and Goldberg? (See \cite{nt} for both.) Using a solution to Question (3) above, or with a better understanding of the non-loop graphs, one could  try to compare the relative complexity of the various approaches. We hope to address the two preceding questions in a sequel to this paper.
\item We have seen that loop graphs are irrelevant for associativity in the linear case. Can this be generalized any further, or are they essential for quadratic and higher degree Poisson structures?
\item {\bf Connections with K.T. Chen's work}. What follows is highly speculative and based on some similarities between approaches to deformation quantization and certain constructs that arise in the works of K.T. Chen  (see \cite{ch2} for a list of all his collected works). In \cite{ch1} K.T. Chen outlined a theory of iterated path integrals in $\rn$ which allowed him to associate certain formal power series with paths in $\rn$. One of the immediate consequences is a proof of the CBH formula. The iterated integrals have a flavor similar to  the ones arising in Kontsevich's weight computations.  Chen builds on this in his later work and one of the basic building blocks is a non-commutative multiplication on the space of paths on a manifold. This seems to bear some relation to the construction of a \stp. Further, Chen generalizes his iterated integrals from path integrals in $\rn$ to iterated integrals for arbitrary forms on differentiable manifolds. One field of  application is the theory of connections  that take values in formal power series. His work here seems to be related to Fedosov's approach to deformation quantization \cite{fe}. It seems that Chen's formalism could provide a bridge for relating the Kontsevich and Fedosov quantization schemes.
\end{list}

\bibliography{fref}
\bibliographystyle{plain}

\end{document}